\documentclass[reqno,12pt]{amsart}
\usepackage[mathscr]{eucal}
\usepackage{amsmath}
\usepackage{amssymb}
\usepackage{latexsym}
\usepackage{amsthm}
\usepackage{color}
\newtheorem{thm}{\indent Theorem}[section]
\newtheorem{cor}[thm]{\indent Corollary}
\newtheorem{lem}[thm]{\indent Lemma}
\newtheorem{prop}[thm]{\indent Proposition}
\newtheorem{dfn}{{\indent\bf Definition}}[section]
\newtheorem{rmk}{{\indent\bf Remark}}[section]
\newtheorem{expl}{{\indent\bf Example}}[section]

\newcommand{\fr}{\frac}
\newcommand{\lmx}{\left(\begin{matrix}}
\newcommand{\rmx}{\end{matrix}\right)}
\newcommand{\ldt}{\left|\begin{matrix}}
\newcommand{\rdt}{\end{matrix}\right|}

\newcommand{\td}{\tilde}

\newcommand{\lang}{\langle}
\newcommand{\rang}{\rangle}
\newcommand{\nnm}{\nonumber}
\newcommand{\bbr}{{\mathbb R}}

\newcommand{\be}{\begin{equation}}
\newcommand{\ee}{\end{equation}}

\newcommand{\vfi}{\varphi}

\makeatletter
\numberwithin{equation}{section}
\begin{document}
\title [Calabi Hypersurfaces with parallel Fubini-Pick form]{{\large \bf Classification of Calabi Hypersurfaces in $\bbr^{n+1}$ with parallel Fubini-Pick form}}
\author [M.X. Lei and R.W. Xu]{Miaoxin Lei\quad\quad Ruiwei Xu}
\address{School of Mathematics and Information Sciences,
\newline \indent Henan Normal University, Xinxiang 453007, P. R. China
\newline \indent lmiaoxin@126.com;\;rwxu@htu.edu.cn }
\date{}
\footnotetext{$^1$\,The authors are supported by NSFC(Grant No.
 11871197).}
\maketitle
\begin{center} {\emph{ In memory of Prof. Xingxiao Li}}\end{center}
\renewcommand{\baselinestretch}{1}
\renewcommand{\arraystretch}{1.2}
\catcode`@=11 \@addtoreset{equation}{section} \catcode`@=12
\makeatletter
\renewcommand{\theequation}{\thesection.\arabic{equation}}
\@addtoreset{equation}{section} \makeatother
\noindent {\bf Abstract}: The classifications of locally strongly convex equiaffine hypersurfaces (resp. centroaffine hypersurfaces) with parallel Fubini-Pick form with respect to the Levi-Civita connection of the Blaschke-Berwald affine metric (resp. centroaffine metric) have been completed by several geometers in the last decades, see \cite{HLV} and \cite{CHM}. In this paper we define a generalized Calabi product in Calabi geometry and prove decomposition theorems in terms of their Calabi invariants. As the main result, we obtain a complete classification of Calabi hypersurfaces in $\bbr^{n+1}$ with parallel Fubini-Pick form with respect to the Levi-Civita connection of the Calabi metric. This result is a counterpart in Calabi geometry of the classification theorems in equiaffine situation \cite{HLV} and centroaffine situation \cite{CHM}.

\vskip 2pt\noindent {\bf
2000 AMS Classification:} Primary 53A15; Secondary 53C24, 53C42.
 \vskip 0.1pt

\section{Introduction}

In the affine differential geometry of hypersurfaces, there are three well-known normalizations, i.e., the equiaffine normalization, the centroaffine normalization, and the Calabi normalization for graph hypersurfaces. For each of the three normalizations, an affine hypersurface is associated with two fundamental invariants, namely a cubic form called the Fubini-Pick form, and a metric called respectively the equiaffine metric (also called the Blaschke-Berwald affine metric), the centroaffine metric and the Calabi metric. If the affine hypersurface is locally strongly convex, the metric becomes definite and can be assumed positive definite. Then, we have the natural problem of classifying locally strongly convex affine hypersurfaces with parallel Fubini-Pick form with respect to the Levi-Civita connection of the preceding metric. After many efforts of several geometers, the classifications of locally strongly convex affine hypersurfaces with parallel Fubini-Pick form have been finally completed in the equiaffine normalization and the centroaffine normalization.
In the article we obtain a complete classification of Calabi hypersurfaces in $\bbr^{n+1}$ with parallel Fubini-Pick form with respect to the Levi-Civita connection of the Calabi metric.

Classification of affine hypersurfaces with parallel Fubini-Pick form  with respect to the Levi-Civita connection of the Blaschke-Berwald metric have been studied intensively, from the earlier beginning paper by Bokan-Nomizu-Simon \cite{BNS}. For affine surfaces, this condition has been studied by Magid and Nomizu in \cite{MN},
see also Theorem 1 of Li and Penn \cite{LP} (locally strongly convex surfaces). The classifications of 3 and 4-dimensional locally strongly convex affine hypersurfaces with parallel Fubini-Pick form are due to Dillen-Vrancken \cite{DV}; and Dillen-Vrancken-Yaprak \cite{DVY}. Later, Hu-Li-Simon-Vrancken \cite{HLSV} presented some new examples of affine spheres and classified locally strongly convex such affine hypersurfaces of $\bbr^{n+1}(n\leq 7)$. In 2011, Hu-Li-Vrancken \cite{HLV} introduced some typical examples and gave the complete classification of locally strongly convex affine hypersurfaces of $\bbr^{n+1}$ with parallel cubic form.

In centroaffine differential geometry, Cheng-Hu-Moruz \cite{CHM} obtained a complete classification of locally strongly convex centroaffine hypersurfaces with parallel cubic form with respect to the Levi-Civita connection of the centroaffine metric.
On the other hand, Liu and Wang \cite{LW} gave the classification of the centroaffine surfaces with parallel traceless cubic form relative to the Levi-Civita connection.
In \cite{CH}, Cheng and Hu established a general inequality for locally strongly convex centroaffine hypersurfaces in $\bbr^{n+1}$ involving the norm of the covariant differentiation of both the difference tensor and the Tchebychev vector field.
Applying the classification results in \cite{CHM}, Cheng and Hu completely classified locally strongly convex centroaffine hypersurfaces with parallel traceless difference tensor, which realize the equality case of the inequality.

Let $\Omega$ be a domain in $\mathbb{R}^n$ and $f$ be a strictly convex $C^\infty$-function on $\Omega$. Consider the \emph{graph immersion} by $f$
$$  M^n:=\{(x,f(x))\;|\; x\in\Omega\}$$
with the so-called {\em Calabi affine normalization}
$$Y =(0,\cdots,0,1)^t\in \bbr^{n+1}.$$
Here we called such graph immersion as \emph{Calabi hypersurface}.
It relative affine metric is the Calabi metric
$G=\sum \frac{\partial^2f}{\partial x_i\partial x_j}dx_idx_j$.
The Fubini-Pick form $A$ (also called cubic form) measures  the difference of the induced affine connection and Levi-Civita connection of Calabi metric $G$.

A Calabi hypersurface is called \emph{canonical} if its Fubini-Pick form is parallel with respect to the Levi-Civita connection of the Calabi metric and its Calabi metric $G$ is flat. Xu and Li \cite{XL} presented a class of new canonical examples denoted by $Q(c_1,\cdots,c_r;n)$ and classified the Calabi complete Tchebychev affine K\"ahler hypersurfaces with nonnegative Ricci curvature. Recently, in \cite{XL1} and \cite{XL2}, the authors classified  $n(\leq4)$-dimensional Calabi hypersurfaces with parallel Fubini-Pick form.

Motivated by the techniques and methods in equiaffine case \cite{HLV} and centroaffine case \cite{CHM}, here we will complete classify $n$-dimensional Calabi hypersurfaces with parallel Fubini-Pick form. To state the main result of this paper, we need the following generalized Calabi product in Calabi geometry, for details see section 3. Let $\varphi: M_1\rightarrow \bbr^{n_1+1}$ be an $n_1$-dimensional locally strongly convex hyperbolic centroaffine hypersurface, and $\psi: M_2\rightarrow \bbr^{n_2+1}$ be an $n_2$-dimensional Calabi hypersurface. For a nonzero constant $\lambda$, we can define the generalized Calabi product of $M_1$ and $M_2$ in Calabi geometry by
\be x(t,p,q):=(e^t\varphi(p),\psi(q))-\frac{t}{\lambda^2}Y \in \bbr^{n_1+n_2+2}, \quad p\in M_1,\; q\in M_2,\; t\in \bbr.\ee

Similarly, the generalized Calabi product of locally strongly convex hyperbolic centroaffine hypersurface $M_1$ and a point is defined by
\be x(t,p):=(e^t\varphi(p),-\frac{t}{\lambda^2})\in \bbr^{n_1+2}, \quad p\in M_1,\; t\in \bbr.\ee
Moreover, the generalized Calabi product of Calabi hypersurface $M_2$ and a point is defined by
\be x(t,q):=(e^t,\psi(q))-\frac{t}{\lambda^2}Y\in \bbr^{n_2+2}, \quad q\in M_2,\; t\in \bbr.\ee

The main theorem of this paper is the following
\begin{thm}\label{main} Let $f$ be a smooth strictly convex function on a domain $\Omega\subset \bbr^n$. If the Calabi hypersurface $M^n=\{(x,f(x)) \;| \;   x\in\Omega\}$  has  parallel cubic form. Then $M$ is Calabi affine equivalent to an open part of one of the following types of hypersurfaces:
\begin{enumerate}
  \item[(i)] elliptic paraboloid $(A\equiv0)$; or
  \item[(ii)] $M^n$ is obtained as the Calabi product of an $(n-1)$-dimensional locally strongly convex hyperbolic centroaffine hypersurface with parallel cubic form and a point; or
  \item[(iii)] $M^n$ is obtained as the Calabi product of an $(n-1)$-dimensional Calabi hypersurface with parallel cubic form and a point; or
  \item[(iv)] $M^n$ is obtained as the Calabi product of a lower dimensional locally strongly convex hyperbolic centroaffine hypersurface with parallel cubic form and a lower dimensional Calabi hypersurface with parallel cubic form.
\end{enumerate}
\end{thm}

\begin{rmk}
Theorem 1.1 implies that all such Calabi hypersurfaces but elliptic paraboloid  can be decomposed as the Calabi product, which is different to cases in the equiaffine differential geometry \cite{HLV} and centroaffine differential geometry \cite{CHM}.\end{rmk}
\begin{rmk}
These hypersurfaces in Theorem 1.1 are Tchebychev affine K\"{a}hler hypersurfaces and affine extremal hypersurfaces.\end{rmk}
\begin{rmk}
For a locally strongly convex hyperbolic centroaffine hypersurface with parallel cubic form, if it can be decomposed as the Calabi product of the lower dimension centroaffine hypersurfaces, then they must be locally strongly convex and hyperbolic, see \cite{CHM}. Combining the complete classification of the locally strong convex hyperbolic centroaffine hypersurfaces with parallel cubic form \cite{CHM} and classification of Calabi hypersurfaces $(n\leq 4)$ with parallel Fubini-Pick form, we obtain the complete classification of Calabi
Calabi hypersurfaces in $\bbr^{n+1}$ with parallel Fubini-Pick form.
\end{rmk}

The paper is organized as follows. In section 2, we review relevant
materials for Calabi hypersurfaces. In section 3, we study both the generalized Calabi product and their characterizations in Calabi geometry, which are main contribution in this paper. In section 4, properties of Calabi hypersurfaces with parallel cubic form in terms of a typical basis are presented, so that the classification problem of such hypersurfaces is divided into $(n+1)$ cases, namely: $\{\mathfrak C_m\}_{0\leq m\leq n}$, depending on the decomposition of the tangent space into three orthogonal distributions, i.e., $\mathcal D_1$ (of dimension one), $\mathcal D_2$ and $\mathcal D_3$. Recall that Cases $\mathfrak C_0$, $\mathfrak C_n$ and $\mathfrak C_{n-1}$ were proved in \cite{XL1} and \cite{XL2}, respectively. As a corollaries of theorems 3.6 and 3.3, the Cases $\mathfrak C_1$ and $\mathfrak C_{n-1}$ will be settled in this section.

To consider the cases $\{\mathfrak C_m\}_{2\leq m\leq n-2}$, we follow closely the same procedure as in \cite{HLV} and \cite{CHM} to introduce two extremely important operators, i.e., an isotropic
bilinear map $L : \mathcal D_2\times\mathcal D_2 \rightarrow \mathcal D_3$ in section 5.1, and, for any unit vector $v \in\mathcal D_2$, the
symmetric linear map $P_v : \mathcal D_2\rightarrow \mathcal D_2$ in section 5.2. With the help of $L$ and $P_v$, we can give a remarkable decomposition of $\mathcal D_2$ in section 5.3. Then in sections 6-10,
according to the decomposition of $\mathcal D_2$, we calculate these cases in much detail in order to achieve the corresponding conclusion, respectively. Finally in section 11 we complete the proof of Theorem 1.1.

\section{Preliminaries}
\subsection {Calabi geometry.}
In this section, we shall recall some basic facts for the Calabi geometry, see \cite{Ca}, \cite{P} or section 3.3.4 in \cite{LXSJ}. Let $f$ be a strictly convex $C^\infty$-function on a domain $\Omega \subset \mathbb{R}^n.$ Consider the graph hypersurface
\be\label{2.1} M^n:=\{(x,f(x))\;|\; x\in\Omega\}.\ee
This Calabi geometry can also be essentially realised as a very special relative affine geometry of the graph hypersurface (\ref{2.1})
by choosing the so-called {\em Calabi affine normalization} with
$$Y =(0,\cdots,0,1)^t\in \bbr^{n+1}$$
being the fixed relative affine normal vector field, which is called the {\em Calabi affine normal}.

For the position vector $x=(x_1,\cdots,x_n,f(x_1,\cdots,x_n))$ we have the decomposition
\be x_{ij}=\fr{\partial^2 x}{\partial x_i\partial x_j}
=\sum c_{ij}^k x_*\frac{\partial }{\partial x_k}+f_{ij}Y,\ee
with respect to the bundle decomposition $\bbr^{n+1}=x_*TM^n\oplus\bbr\cdot Y$, where the induced affine connection coefficients $c^k_{ij}\equiv 0$. It follows that the relative affine metric is the Calabi metric
$$ G=\sum \frac{\partial^2f}{\partial x_i\partial x_j}dx_idx_j.$$
The Levi-Civita connection with respect to the metric $G$ has the Christoffel symbols
$$\Gamma_{ij}^k=\frac{1}{2}\sum f^{kl}f_{ijl},$$
where and hereafter
$$f_{ijk}=\frac{\partial^3f}{\partial x_i\partial x_j\partial x_k}, \quad (f^{kl})=(f_{ij})^{-1}. $$
Then we can rewrite the Gauss structure equation as follows:\be x_{,ij}=\sum A_{ij}^kx_*\frac{\partial }{\partial x_k}+G_{ij}Y.\ee
The Fubini-Pick tensor (also called cubic form) $A_{ijk}$ and the Weingarten tensor $B_{ij}$ satisfy
\be A_{ijk}=\sum A_{ij}^lG_{kl}=-\frac{1}{2}f_{ijk},\quad B_{ij}=0,\ee
which means that $A_{ijk}$ are symmetric in all indexes.
Classically, the tangent vector field
\be T:=\frac{1}{n}\sum G^{kl}G^{ij}A_{ijk}\frac{\partial}{\partial x_l}=-\frac{1}{2n}\sum f^{kl}f^{ij}f_{ijk}\frac{\partial}{\partial x_l}\ee
is called the \emph{Tchebychev vector field} of the hypersurface $M^n$, and the invariant function
\be J:=\frac{1}{n(n-1)}\sum G^{il}G^{jp}G^{kq}A_{ijk}A_{lpq}=\frac{1}{4n(n-1)}\sum f^{il}f^{jp}f^{kq}f_{ijk}f_{lpq}\ee
is named as the  \emph{relative Pick invariant}  of $M^n$. As in the report \cite{Li},  $M^n$ is called a \emph{Tchebychev affine K\"{a}hler hypersurface},
if the Tchebychev vector field $T$ is parallel with respect to the Calabi metric $G$.
The \emph{Gauss} and \emph{Codazzi equations} read
\be \label{2.7} R_{ijkl}=\sum f^{mh}(A_{jkm}A_{hil}-A_{ikm}A_{hjl}),\ee
\be \label{2.8} A_{ijk,l}=A_{ijl,k}.\ee
From (\ref{2.7}) we get the  \emph{Ricci tensor}
\be\label{2.9} R_{ik}=\sum f^{jl}f^{mh}(A_{jkm}A_{hil}-A_{ikm}A_{hjl}).\ee
Thus, the scalar curvature is given by
\be\label{2.10} R=n(n-1)J-n^2|T|^2.\ee

Let $A(n+1)$ be the group of $(n+1)$-dimensional affine transformations on $\bbr^{n+1}$. Then $A(n+1)=GL(n+1)\ltimes \bbr^{n+1}$, the semi-direct product of the general linear group $GL(n+1)$ and the group $\bbr^{n+1}$ of all the parallel translations on $\bbr^{n+1}$. Define
\be
SA(n+1)=\{(M,b)\in A(n+1)=GL(n+1)\ltimes \bbr^{n+1};\ M(Y)=Y\}
,\ee
where $Y=(0,\cdots,0,1)^t$ is the Calabi affine normal. Then the subgroup $SA$ consists of all the transformation $\phi \in SA(n+1)$ of the following type:
\begin{align}
X:&=(X^j,X^{n+1})^t\equiv(X^1,\cdots, X^n,X^{n+1})^t\nnm\\
&\mapsto \phi(X):=\lmx a^i_j&0\\a^{n+1}_j&1\rmx X+b,\quad \forall\, X\in \bbr^{n+1},
\end{align}
for some $(a^i_j)\in A(n)$, constants $a^{n+1}_j$ ($j=1,\cdots,n$) and some constant vector $b\in \bbr^{n+1}$. Clearly, the Calabi metric $G$ is invariant under the action of $SA(n+1)$ on the graph hypersurfaces or, equivalently, under the induced action of $SA(n+1)$ on the strictly convex functions, which is naturally defined to be the composition of the following maps:
$$f\mapsto (x_i,f(x_i))\mapsto(\td x_i,\td f(\td x_i)):=\phi(x_i,f(x_i))\mapsto \phi(f):=\td f,\quad\forall\,\phi\in SA(n+1).$$
\begin{dfn}\cite{XL1}
Two graph hypersurfaces $(x_i,f(x_i))$ and $(\td x_i,\td f(\td x_i))$, defined respectively in domains $\Omega,\td\Omega\subset \bbr^n$, are called Calabi-affine equivalent if they differ only by an affine transformation $\phi\in SA(n+1)$.
\end{dfn}

Accordingly, we have

\begin{dfn}\cite{XL1}
Two smooth functions $f$ and $\td f$ respectively defined on domains $\Omega,\td\Omega\subset \bbr^n$ are called affine equivalent $($related with an affine transformation $\vfi\in A(n))$ if there exist some constants  $a^{n+1}_1,\cdots,a^{n+1}_n,b^{n+1}\in \bbr$ such that $\vfi(\Omega)=\td\Omega$ and
\be\td f(\vfi(x)) =f(x)+\sum a^{n+1}_jx_j+b^{n+1}\quad
\text{for all}\ x=(x_1,\cdots,x_n)\in\Omega.\ee
\end{dfn}

Clearly, the above two definitions are equivalent to each other.

In \cite{XL} Xu and Li introduced a large class of new canonical Calabi hypersurfaces, which are denoted by $Q(c_1,\cdots,c_r;n)$, see also \cite{XL1}.

\begin{expl} Given the dimension $n$ and $1\leq r\leq n$, we define
$$\Omega_{r,n}=\{(x_1,\cdots,x_n)\in \bbr^n;\ x_1>0,\cdots,x_r>0\}.$$
Then, for any positive numbers $c_1,\cdots, c_r$, we have a Calabi hypersurface given by the function
\begin{align}
f(x_1,\cdots,x_n):=&-\sum_{i=1}^r c_i\ln x_i+\frac12\sum_{j=r+1}^n x_j^2,\qquad (x_1,\cdots,x_n)\in\Omega_{r,n}.\label{2.14}
\end{align}
In follows, we shall denote such Calabi hypersurfaces by $Q(c_1,\cdots,c_r;n)$.
\end{expl}

\section{Generalized Calabi product and Decomposed theorem}

In this section, we should study the generalized Calabi products as defined in (1.1)--(1.3) in Calabi geometry. Firstly, we state some elementary calculations on Calabi product, formulated as Propositions 3.1 and 3.3. Then, considering the converse of these propositions, we shall prove Theorems \ref{thm3.2} and \ref{thm3.5}.

Let $\varphi: M_1\rightarrow \bbr^{n_1+1}$ be a locally strongly convex hyperbolic centroaffine hypersurface of dimension $n_1$. Denote by $g$ the centroaffine metric of $\varphi$ and by $\{u_1,\cdots,u_{n_1}\}$ local coordinates for $M_1$. For a nonzero constant $\lambda$, the generalized Calabi product of hyperbolic centroaffine hypersurface $M_1$ and a point
$$ x: M^{n_1+1}=\bbr\times M_1\rightarrow \bbr^{n_1+2}$$
is defined by
\be x(t,p):=(e^t\varphi(p),-\frac{t}{\lambda^2}), \qquad p\in M_1,\; t\in \bbr.\ee
Then we can state the following
\begin{prop}
The Calabi product of locally strongly convex hyperbolic centroaffine hypersurface $M_1$ and a point
$$x: M^{n_1+1}=\bbr\times M_1\rightarrow \bbr^{n_1+2}$$
defined by (3.1) is a Calabi hypersurface, and the Calabi metric $G$  is given by
\be G=\frac{1}{\lambda^2}(dt^2\oplus g).\ee
 The Fubin-Pick tensor $A$ of $x$ takes the following form:
\be A(\tilde T,\tilde T )= \lambda\tilde T, \quad A(\tilde T,  \frac{\partial x}{\partial u_i})=\lambda \frac{\partial x}{\partial u_i}, \qquad \forall 1 \leq i\leq n_1,\ee
where $\tilde T:= \lambda \frac{\partial x}{\partial t}$ denotes a $G$-unit tangential vector field.

Moreover, $x$ is flat (resp. of parallel cubic form) if and only if $\varphi$ is flat (resp. of parallel cubic form).
\end{prop}

\textbf{Proof.} The following calculations are elementary. \be x_t= (e^t\varphi,-\frac{1}{\lambda^2} ),\quad x_{u_i} = (e^t\varphi_{i},0 ), \quad Y=(0,...,0,1)\in R^{n+1} .\ee
From \be x_{tt}= (e^t\varphi,0)=x_{\ast}(\nabla_{\partial t}\partial t)+G_{tt}Y,\ee
we have \be \nabla_{\partial t}\partial t= \partial t,\qquad G_{tt}=\frac{1}{\lambda^2}.\ee
Similarly, from \be x_{tu_i}= (e^t\varphi_{i},0)=x_{\ast}(\nabla_{\partial t}\partial u_i)+G_{tu_i}Y,\ee
we get\be \nabla_{\partial t}\partial u_i= \partial u_i,\qquad G_{tu_i}=0.\ee
And, from \begin{align} x_{u_iu_j}=& (e^t\varphi_{ij},0)=e^t(\Gamma_{ij}^k\varphi_k+g_{ij}\varphi,0)\notag\\
= &x_{\ast}(\nabla_{\partial u_i}\partial u_j)+G_{u_iu_j}Y\notag\\
=& a_{ij}x_t+ c_{ij}^k x_{u_k}+ G_{u_iu_j}Y,
\end{align}
we have\be a_{ij}=g_{ij},\qquad G_{u_iu_j}=\frac{1}{\lambda^2}g_{ij}.\ee
Then the Calabi metric is given by
\be G=\frac{1}{\lambda^2}(dt^2\oplus g).\ee
Therefore, the difference tensor $A$ of $x$ takes the following form:
\be A(\tilde T,\tilde T )= \lambda\tilde T, \quad A(\tilde T,  x_{u_i})=\lambda x_{u_i}, \qquad \forall 1 \leq i\leq n_1.\ee
The rest conclusions of Proposition is easily to be verified. \hfill $\Box$

By the above method, we can construct Calabi hypersurfaces with parallel cubic form as follows:

\begin{expl}
Assume that $$ \textbf{x} =(x_1,x_2,...,x_n, f(x))=(e^{t}(y_1,y_2,..., y_n),-\frac{t}{\lambda^2}). $$
\begin{enumerate}
\item If $\varphi(M_1)$ is chosen to be $y_1^\alpha y_2^\beta y_3^\gamma=1$ ($\alpha, \beta, \gamma>0$), which is a canonical hyperbolic centroaffine surface, then the Calabi product of $\varphi(M_1)$ and a point is Calabi affine equivalent to
 \be x_4=-\frac{1}{(\alpha+\beta+\gamma)\lambda^2}(\alpha\ln x_1+\beta\ln x_2+\gamma\ln x_3).\ee \\
\item If $\varphi(M_1)$ is chosen to be the hyperboloid $y_1^2+y_2^2-y_3^2=-1$, then the Calabi product of $\varphi(M_1)$ and a point is Calabi affine equivalent to
\be x_4=-\frac{1}{2\lambda^2}\ln(x_3^2-( x_1^2+ x_2^2)). \ee \\
\item If $\varphi(M_1)$ is chosen to be
 $(y_1^2-y_2^2-y_3^2)^\alpha y^\beta_4=1$ ($\alpha, \beta>0$),  which is 3-D locally strongly convex hyperbolic centroaffine hypersurface with $\hat\nabla A=0$ (see (3.2) in \cite{CHM}), then the Calabi product of $\varphi(M_1)$ and a point is Calabi affine equivalent to
\be x_5=-\frac{1}{(2\alpha+\beta)\lambda^2}[\beta\ln x_4+ \alpha\ln(x_1^2-(x_2^2+ x_3^2))]. \ee
 \end{enumerate}
\end{expl}

Next, as the converse of Proposition 3.1, we can prove the following theorem.
\begin{thm}\label{thm3.2}
Let $x:M^n\rightarrow \bbr^{n+1}$ be a Calabi hypersurface. Assume that there exist two distributions $\mathcal {D}_1$ (of dimension 1, spanned by a unit vector filed $\tilde T$), $\mathcal {D}_2$ (of dimension $n-1$) such that
\begin{enumerate}
  \item[(i)] $\mathcal {D}_1$ and $\mathcal {D}_2$ are orthogonal with respect to the Calabi metric $G$;
  \item[(ii)] there exists a nonzero constant $\lambda$ such that $$A(\tilde T,\tilde T)=\lambda\tilde T,\quad A(\tilde T,V)=\lambda V, \quad   \forall V\in \mathcal {D}_2.$$
\end{enumerate}
Then $x:M^n\rightarrow \bbr^{n+1}$ can be locally decomposed as the Calabi product of a locally strongly convex hyperbolic centroaffine hypersurface $\varphi: M_1^{n-1}\rightarrow \bbr^{n}$ and a point.
\end{thm}

\textbf{Proof.} Note that
\begin{align}
(\hat \nabla A)(V,\tilde T,\tilde T)=&\hat \nabla_V A(\tilde T,\tilde T)-2A(\hat \nabla_V\tilde T,\tilde T )\notag\\ = &\lambda\hat \nabla_V  \tilde T -2A(\hat \nabla_V\tilde T,\tilde T ). \end{align}
Similarly,\begin{align} (\hat \nabla A)(\tilde T,V,\tilde T)= &\hat \nabla_{\tilde T} A(V,\tilde T)-A(\hat \nabla_{\tilde T} V,\tilde T )- A(V,\hat \nabla_{\tilde T} \tilde T )\notag\\
= &\lambda\hat \nabla_{\tilde T}V -A(\hat \nabla_{\tilde T} V,\tilde T )- A(V,\hat \nabla_{\tilde T} \tilde T ).\end{align}
Taking inner product of the right hand sides of (3.16) and (3.17) with $\tilde T$, (resp. $\tilde V\in \mathcal {D}_2$) and using the Codazzi equations and the symmetry of Fubini-Pick form, we obtain
\be \hat\nabla_{\tilde T} \tilde T=0,\qquad \langle\hat \nabla_{V}\tilde T,\tilde V\rangle=0.\ee
Thus, for any vector $X\in TM$, the following conclusions hold
\be\hat\nabla_X\tilde T=0,\;\;\hat\nabla_XV\in\mathcal {D}_2.\ee

This is sufficient to conclude that locally $M$ is isometric with $I\times M_1$, where $\tilde T$ is tangent to $I$ and $\mathcal {D}_2$ is tangent to $M_1$. We may assume that $\tilde T=\lambda\frac{\partial}{\partial t}$. Put
\be\varphi:=e^{-t}(\frac{\tilde T}{\lambda}+\frac{Y}{\lambda^2}),\qquad \psi:= - \frac{1}{\lambda}\tilde T+x+\frac{t}{\lambda^2}Y.\ee
where $Y=(0,...,0,1)^t\in \bbr^{n+1}$. It follows (3.20) that
\begin{align}
   D_{\tilde T}\varphi=&-\lambda e^{-t}(\frac{\tilde T}{\lambda}+\frac{Y}{\lambda^2})+ \lambda^{-1} e^{-t}D_{\tilde T}\tilde T\notag\\
   =& -e^{-t}(\tilde T+\frac{Y}{\lambda})+ \lambda^{-1} e^{-t}(\hat\nabla_{\tilde T}\tilde T+ A_{\tilde T}\tilde T+G(\tilde T,\tilde T)Y)\notag\\
   = \;& 0.
\end{align}

Similarly, for any $V\in \mathcal {D}_2$,
\begin{align}
   D_{\tilde T}\psi=\,& D_V \psi=0,\\
    d\varphi(V)=\,& D_V\varphi=e^{-t}V.
\end{align}
The above relations imply that $\varphi$ reduces to a map of $M_1$ in $\bbr^{n+1}$. Whereas $\psi$ is a constant vector in $\bbr^{n+1}$. Moreover, denoting by $\nabla^1$
the $\mathcal {D}_2$ component of $\nabla$ and for any $V,\tilde V\in \mathcal {D}_2$, we find that
\begin{align}
   D_{V}d\varphi(\tilde V)=\;& e^{-t}D_V \tilde V\notag \\
    =\;&  e^{-t}(\nabla^1_V \tilde V+ \langle\nabla_V\tilde V,\tilde T  \rangle\tilde T+G(V,\tilde V)Y )\notag\\
    =\;&   d\varphi(\nabla^1_V \tilde V) + e^{-t} ( \langle\hat\nabla_V\tilde V+A_V\tilde V,\tilde T  \rangle\tilde T+G(V,\tilde V)Y )\notag\\
    =\;&   d\varphi(\nabla^1_V \tilde V) + \lambda^2G(V,\tilde V)\varphi.
\end{align}
Hence $\varphi$ can be interpreted as a hyperbolic centroaffine immersion contained in an $n$-dimensional vector subspace of $\bbr^{n+1}$ with induced connection $\nabla^1$ and centroaffine metric $g=\lambda^2 G$.

Solving (3.20) for the immersion $x$, we have
\be x=e^t\varphi-\frac{t}{\lambda^2}Y+(\psi-\frac{1}{\lambda^2}Y). \ee
Since $M^n$ is nondegenerate, $x$ lies full in $\bbr^{n+1}$. Therefore $\varphi$ lies in the space spanned by the first $n$ coordinates of $\bbr^{n+1}$. By an affine transformation $\phi\in SA(n+1)$, we can assume that
\be x=e^t\varphi-\frac{t}{\lambda^2}Y=(e^t\varphi, -\frac{t}{\lambda^2}). \ee
We have complete the proof of Theorem 3.2.\hfill $\Box$

In Theorem 3.2, if additionally $M^n$ has parallel cubic form $\hat\nabla A=0$, then we can prove the following theorem.

\begin{thm}\label{thm3.3}
Let $x:M^n\rightarrow \bbr^{n+1}$ be a Calabi hypersurface. Assume that $\hat\nabla A=0$ and there exist orthogonal distributions $\mathcal {D}_1$ (of dimension 1, spanned by a unit vector filed $\tilde T$), $\mathcal {D}_2$ (of dimension $n-1$) and  a nonzero constant $\lambda$ such that
$$A(\tilde T,\tilde T)=\lambda\tilde T,\quad A(\tilde T,V)=\lambda V, \qquad
  \forall V\in \mathcal {D}_2.$$
Then $x:M^n\rightarrow \bbr^{n+1}$ can be locally decomposed as the Calabi product of a locally strongly convex hyperbolic centroaffine hypersurface $\varphi: M_1^{n-1}\rightarrow \bbr^{n}$ with parallel cubic form and a point.
\end{thm}

\begin{rmk}  Note that the unit vector field $\tilde T$ of $M^n$ is not necessarily parallel to its Tchebychev vector field $T$. If it happens, then the immersion $\varphi: M_1^{n-1}\rightarrow \bbr^{n}$ is a locally strong convex hyperbolic affine hypersphere in Theorem 3.3.
\end{rmk}

Let $\varphi: M_1\rightarrow \bbr^{n_1+1}$ be a locally strongly convex hyperbolic centroaffine hypersurface of dimension $n_1$. Let $\psi: M_2\rightarrow \bbr^{n_2+1}$ be a Calabi hypersurface of dimension $n_2$. For a nonzero constant $\lambda$, the Calabi product of $M_1$ and $M_2$ is defined by
\be x(t,p,q):=(e^t\varphi(p),\psi(q))-\frac{t}{\lambda^2}Y, \qquad p\in M_1,\quad q\in M_2,\quad t\in \bbr.\ee
Denote by $g^1$ the centroaffine metric of $\varphi$ and by $\{u_1,\cdots,u_{n_1}\}$ local coordinates for $M_1$. Denote by $g^2$ the Calabi metric of $\psi$ and by $\{u_{n_1+1},\cdots,u_{n_1+n_2}\}$ local coordinates for $M_2$. For simplicity, we use the following range of indices:
$$1 \leq i,j,k \leq n_1,\quad n_1+1\leq\alpha,\beta,\gamma\leq n_1+n_2.$$

\begin{prop}
The Calabi product of $M_1$ and $M_2$
$$ x: M^{n_1+n_2+1}=\bbr\times M_1\times M_2\rightarrow \bbr^{n_1+n_2+2}$$
defined by (3.27)  is a Calabi hypersurface, and the Calabi metric $G$  is given by
\be G=\frac{1}{\lambda^2}(dt^2\oplus g^1)\oplus g^2.\ee
The difference tensor $A$ of $x$ takes the following form:
\be A(\tilde T,\tilde T )= \lambda\tilde T, \quad A(\tilde T,  x_{u_i})=\lambda x_{u_i}, \quad  A(\tilde T,  x_{u_\alpha})=0,\quad A(x_{u_i},  x_{u_\alpha})=0.\ee

Moreover, $x$ is flat (resp. of parallel cubic form) if and only if both $\varphi$ and $\psi$ are flat (resp. of parallel cubic form).
\end{prop}

\textbf{Proof.} Denote by $Y=(0,...,0,1)\in \bbr^{n_1+n_2+2}$ the Calabi normalization.
\be x_t= (e^t\varphi, 0)-\frac{1}{\lambda^2}Y,\quad x_{u_i} = (e^t\varphi_{u_i},0 ),\quad x_{u_\alpha} = (0,\psi_{u_\alpha} ).\ee
From\be x_{tt}= (e^t\varphi,0)=x_{\ast}(\nabla_{\partial t}\partial t)+G_{tt}Y=x_t+\frac{1}{\lambda^2}Y,\ee
we have\be \nabla_{\partial t}\partial t= \partial t,\qquad G_{tt}=\frac{1}{\lambda^2}.\ee
Similarly, by\be x_{tu_i}= (e^t\varphi_{u_i},0)=x_{\ast}(\nabla_{\partial t}\partial u_i)+G_{tu_i}Y=x_{u_i},\ee
we get\be \nabla_{\partial t}\partial u_i= \partial u_i,\qquad G_{tu_i}=0.\ee
By $x_{u_iu_\alpha}=0$ and $x_{tu_\alpha}= 0$, we get
\be  \nabla_{\partial u_i}\partial u_\alpha= 0,\quad G_{u_iu_\alpha}=0,\quad
 \nabla_{\partial t}\partial u_\alpha= 0,\quad G_{tu_\alpha}=0.\ee
By
\be\frac{\partial^2 \varphi}{\partial u_i \partial u_j}=\Gamma_{ij}^k\varphi_k +g^1_{ij}\varphi,\ee
we have\begin{align} x_{u_iu_j}=& (e^t\varphi_{ij},0)=e^t(\Gamma_{ij}^k\varphi_k+g_{ij}^1\varphi,0)\notag\\
= &x_{\ast}(\nabla^1_{\partial u_i}\partial u_j)+a_{ij}x_t+D_3^T
+G_{u_iu_j}Y\notag\\
= &x_{\ast}(\nabla^1_{\partial u_i}\partial u_j)+a_{ij}(e^t\varphi,0)
+(G_{u_iu_j}-\frac{a_{ij}}{\lambda^2})Y,
\end{align}
where denote by $\nabla^1$ and $D^T_3$ the tangential part  of $\nabla_{\partial i}{\partial j}$ onto  $M_1$ and $M_2$, respectively.
Then
$$a_{ij}=g^1_{ij},\qquad G_{u_iu_j}=\frac{1}{\lambda^2}g^1_{ij}.$$

Similarly, by
$$(0,\psi_{u_\alpha u_\beta})=(0,\Gamma_{\alpha\beta}^\gamma\psi_\gamma+g^2_{\alpha\beta}\tilde Y)$$$$=x_{u_\alpha u_\beta}=x_{\ast}(\nabla_{\partial u_\alpha}{\partial u_\beta})+G_{u_\alpha u_\beta }Y=x_{\ast}(\nabla^2_\alpha\beta)+G_{u_\alpha u_\beta }Y, $$
we get
$$ G_{u_\alpha u_\beta}=g^2_{\alpha\beta}.$$

Then the Calabi metric is given by
\be G=\frac{1}{\lambda^2}(dt^2\oplus g^1)\oplus g^2.\ee
Therefore, the difference tensor $A$ of $x$ takes the following form:
\be A(\tilde T,\tilde T )= \lambda\tilde T, \quad A(\tilde T,  x_{u_i})=\lambda x_{u_i}, \quad     A(\tilde T,  x_{u_\alpha})=0\quad A(x_{u_i},  x_{u_\alpha})=0..\ee
where $\tilde T=\lambda x_t$ is a unit vector field.\hfill $\Box$

\begin{rmk} Here $n_1$ and $n_2$ may allow to zero  if we define 0-dimensional centroaffine (resp. Calabi) hypersurfaces to be the positive number (resp. arbitrary number) constant maps with vanishing centroaffine (resp. Calabi) metric and Fubini-Pick form. If $n_2=0$, Proposition 3.4 is exactly Proposition 3.1. If $n_1=0$ and the Calabi hypersurface $\psi$ is graph of a strictly convex function $f$, then the Calabi product $x$ is Calabi affine equivalent to graph $$f(x_2,...,x_{n_2+1})-\frac{1}{\lambda^2}\ln x_1.$$
\end{rmk}

By the above method, we can construct Calabi hypersurfaces with parallel cubic form as follows:

\begin{expl}
 \begin{enumerate}
\item If $n_1=0$ and $\psi(M_2)$ is chosen to be $y_2=\frac{1}{2}y_1^2$,
then the Calabi product of $\psi(M_2)$ and a point $$( e^t,y_1,y_2-\frac{t}{\lambda^2})$$
is Calabi affine equivalent to one of $Q(c_1,2)$,
$x_3=\frac{x_2^2}{2}-\frac{1}{\lambda^2}\ln x_1.$
\item
If we choose that $\varphi(M_1)$ is $y_1^\alpha y_2^\beta y_3^\gamma=1$
($\alpha, \beta, \gamma>0$), which is a hyperbolic canonical centroaffine surface, and $\psi(M_2)$ is elliptic paraboloid $z_3=z_1^2+z_2^2$, then the Calabi product of $\varphi(M_1)$ and  $\psi(M_2)$  $$x=(e^t(y_1,y_2,y_3), z_1,z_2,z_3-\frac{t}{\lambda^2})$$
is Calabi affine equivalent to
\be x_6=-\frac{1}{(\alpha+\beta+\gamma)\lambda^2}(\alpha\ln x_1+\beta\ln x_2+\gamma\ln x_3)+x_4^2+x_5^2,\ee which is a canonical Calabi hypersurface in $\bbr^6$.
\end{enumerate}
\end{expl}

Next, as the converse of Proposition 3.4, using the similar proof of Lemmas 1-4 in \cite{HLV}, we can prove the following theorem.

\begin{thm}\label{thm3.5}
Let $x:M^n\rightarrow \bbr^{n+1}$ be a Calabi hypersurface. Assume that there exist distributions $\mathcal {D}_1$ (of dimension 1, spanned by a unit vector filed $\tilde T$), $\mathcal {D}_2$ (of dimension $n_1$),  $\mathcal {D}_3$ (of dimension $n_2$) such that
\begin{enumerate}
  \item[(i)] $1+n_1+n_2=n$,
  \item[(ii)] $\mathcal {D}_1$, $\mathcal {D}_2$ and $\mathcal {D}_3$ are mutually orthogonal with respect to the Calabi metric $G$,
  \item[(iii)] there exists a nonzero constant $\lambda$ such that, for any $V\in \mathcal {D}_2, W\in \mathcal {D}_3$, the following relations hold
  $$A(\tilde T,\tilde T)=\lambda\tilde T,\, A(\tilde T,V)=\lambda V,\, A(\tilde T,W)=0,\, A(V,W)=0.$$
\end{enumerate}
Then $x:M^n\rightarrow \bbr^{n+1}$ can be locally decomposed as the Calabi product of a locally strongly convex hyperbolic centroaffine hypersurface $\varphi: M_1^{n_1}\rightarrow \bbr^{n_1+1}$ and a Calabi hypersurface $\psi:M_2^{n_2}\rightarrow \bbr^{n_2+1}$.
\end{thm}
\textbf{Proof.} Taking inner product of the right hand sides of (3.16) and (3.17) with $\tilde T$, (resp. $\tilde V\in \mathcal {D}_2$) and using the Codazzi equations and the symmetry of Fubini-Pick form, we obtain
\be \hat\nabla_{\tilde T} \tilde T\in\mathcal {D}_2^{\bot},\qquad \hat \nabla_{V} \tilde V\in \mathcal {D}_2\oplus\mathcal {D}_3.\ee
Using the similar method, we obtain
\be \hat\nabla_{\tilde T} \tilde T\in\mathcal {D}_3^{\bot},\qquad \hat \nabla_{W} \tilde W\in \mathcal {D}_2\oplus\mathcal {D}_3.\ee

Similarly,
\begin{align} (\hat \nabla A)(\tilde T,V,W)= &\hat \nabla_{\tilde T} A(V,W)-A(\hat \nabla_{\tilde T} V,W )- A(V,\hat \nabla_{\tilde T} W )\notag\\
= &  -A(\hat \nabla_{\tilde T} V,W )- A(V,\hat \nabla_{\tilde T} W ).\end{align}
\begin{align} (\hat \nabla A)(V,\tilde T,W)= &\hat \nabla_{V} A(\tilde T,W)-A(\hat \nabla_{V}\tilde T,W )- A(\tilde T,\hat \nabla_{V} W )\notag\\
= & -A(\hat \nabla_{V}\tilde T,W )- A(\tilde T,\hat \nabla_{V} W ).\end{align}
Taking inner product of the right hand sides of (3.43) and (3.44) with $\tilde T$, and using the Codazzi equations and the symmetry of Fubini-Pick form, we obtain
\be\langle\hat \nabla_{\tilde T} W,  V\rangle= -\langle\hat \nabla_V\tilde T, W\rangle.\ee
Interchanging $V$ and $W$ in (3.43) and (3.44), we have
\be\langle\hat \nabla_{\tilde T} W,  V\rangle= \langle\hat \nabla_W\tilde T, V\rangle.\ee

Similarly,
\be  (\hat \nabla A)(\tilde T,V,\tilde V)=  \hat \nabla_{\tilde T} A(V,\tilde V )-A(\hat \nabla_{\tilde T} V,\tilde V )- A(V,\hat \nabla_{\tilde T} \tilde V ); \ee
\begin{align} (\hat \nabla A)(V,\tilde T,\tilde V)= &\hat \nabla_{V} A(\tilde T,\tilde V)-A(\hat \nabla_{V}\tilde T,\tilde V )- A(\tilde T,\hat \nabla_{V}\tilde V )\notag\\
= &\lambda\hat\nabla_V \tilde V -A(\hat \nabla_{V}\tilde T,\tilde V )- A(\tilde T,\hat \nabla_{V} \tilde V ).\end{align}
Taking inner product of the right hand sides of (3.47) and (3.48) with $W$, and using (3.46), we obtain
\be -\lambda\langle\hat \nabla_V\tilde V, W\rangle =\langle A(V,\tilde V), \hat \nabla_{\tilde T}W\rangle= \langle A(V,\tilde V), \hat \nabla_W{\tilde T}\rangle.\ee

Similarly,
\be  (\hat \nabla A)(\tilde T,W,\tilde W)=  \hat \nabla_{\tilde T} A(W,\tilde W )-A(\hat \nabla_{\tilde T} W,\tilde W )- A(W,\hat \nabla_{\tilde T} \tilde W ); \ee
\begin{align} (\hat \nabla A)(W,\tilde T,\tilde W)= &\hat \nabla_{W} A(\tilde T,\tilde W)-A(\hat \nabla_{W}\tilde T,\tilde W )- A(\tilde T,\hat \nabla_{W}\tilde W )\notag\\
= & -A(\hat \nabla_{W}\tilde T,\tilde W )- A(\tilde T,\hat \nabla_{W} \tilde W ).\end{align}
Taking inner product of the right hand sides of (3.50) and (3.51) with $V$, and using (3.45), we obtain
\be \lambda\langle\hat \nabla_W\tilde W, V\rangle = \langle A(W,\tilde W), \hat \nabla_{\tilde T}V\rangle=\langle A(W,\tilde W), \hat \nabla_V{\tilde T}\rangle.\ee

\textbf{Claim}: $\hat \nabla_V\tilde T=\hat \nabla_W\tilde T=0,\quad \forall V\in\mathcal {D}_2,W\in\mathcal {D}_3.$

 In fact, using the Gauss equation, we have that
 \begin{align}  \langle R(V,\tilde T)\tilde T,V\rangle= & -\langle A_VA_{\tilde T}\tilde T, V\rangle
 +\langle A_{\tilde T}A_V\tilde T, V \rangle\notag\\
=  & -\langle A_{\tilde T}\tilde T, A_VV\rangle
 +\langle A_V\tilde T, A_{\tilde T}V \rangle\notag\\
 =  & (-\lambda^2
 +\lambda^2)\langle V,V \rangle=0.\end{align}
On the other hand, by a direct computation and using relations (3.45)-(3.46), we have
\begin{align}  \langle R(V,\tilde T)\tilde T,V\rangle= & \langle \hat\nabla_V\hat\nabla_{\tilde T}\tilde T-\hat\nabla_{\tilde T}\hat\nabla_V\tilde T-\hat\nabla_{\hat\nabla_V\tilde T-\hat\nabla_{\tilde T} V   } \tilde T, V\rangle \notag\\
=  &  -\langle\hat\nabla_{\tilde T}\hat\nabla_V\tilde T,V\rangle
 -\sum_{k=1}^{n_2}\langle\hat\nabla_V\tilde T-\hat\nabla_{\tilde T} V,W_k    \rangle \langle\hat\nabla_{W_k}\tilde T,V\rangle \notag\\
 =  & -\langle\hat\nabla_{\tilde T}\hat\nabla_V\tilde T,V\rangle\notag\\
 =  & - \sum_{k=1}^{n_2} \langle\hat\nabla_V\tilde T,W_k   \rangle \langle\hat\nabla_{\tilde T}W_k,V\rangle\notag\\
=  &  \langle\hat\nabla_V\tilde T,\hat\nabla_V\tilde T \rangle.\end{align}
Then we have $\hat \nabla_V\tilde T=0$. Similarly, we can prove $\hat \nabla_W\tilde T=0$.

Thus we obtain that: for any vector $X\in TM$, $V\in \mathcal {D}_2$ and $W\in \mathcal {D}_3$, the following conclusions hold
\be\hat\nabla_X\tilde T=0,\;\;\hat\nabla_XV\in\mathcal {D}_2,\;\;\hat\nabla_XW\in\mathcal {D}_3.\ee

This is sufficient to conclude that $(M,G)$ is locally isometric to $I\times M_1\times M_2$, where $\tilde T$ is tangent to $I$ and $\mathcal {D}_2$ (resp.$\mathcal {D}_3$ ) is tangent to $M_1$ (resp. $M_2$). We may assume that $\tilde T=\lambda \frac{\partial}{\partial t}$. Put
\be\varphi=e^{-t}(\frac{\tilde T}{\lambda}+\frac{Y}{\lambda^2}),\qquad \psi= - \frac{1}{\lambda}\tilde T+x+\frac{t}{\lambda^2}Y.\ee
It follows (3.56) that
\begin{align}
   D_{\tilde T}\varphi=&-\lambda e^{-t}(\frac{\tilde T}{\lambda}+\frac{Y}{\lambda^2})+ \lambda^{-1} e^{-t}D_{\tilde T}\tilde T\notag\\
   =& -\lambda e^{-t}(\frac{\tilde T}{\lambda}+\frac{Y}{\lambda^2})+ \lambda^{-1} e^{-t}(\hat\nabla_{\tilde T}\tilde T+ A_{\tilde T}\tilde T+G(\tilde T,\tilde T)Y)\notag\\
   = & 0.
\end{align}

Similarly, for any $V\in \mathcal {D}_2$ and $W\in \mathcal {D}_3$,
\be D_W\varphi=0,\qquad d\varphi(V)= D_V\varphi=e^{-t}V.\ee
\be D_{\tilde T} \psi=D_V\psi=0,\quad d\psi (W)=D_W\psi=W.\ee
The above relations imply that $\varphi$ (resp. $\psi$) reduces to a map of $M_1$ (resp. $M_2$) in $\bbr^{n+1}$. Moreover both maps $\varphi$ and $\psi$ are actually immersions. Denoting by $\nabla^1$ (resp. $\nabla^2$)
the $\mathcal {D}_2$ (resp. $\mathcal {D}_3$) component of $\nabla$ and for any $V,\tilde V\in \mathcal {D}_2$, we further find that
\begin{align}
   D_{V}d\varphi(\tilde V)=& e^{-t}D_V \tilde V\notag \\
    =&  e^{-t}(\nabla^1_V \tilde V+ \langle\nabla_V\tilde V,\tilde T  \rangle\tilde T+\nabla_V^2\tilde V+G(V,\tilde V)Y )\notag\\
    =&   d\varphi(\nabla^1_V \tilde V) + e^{-t} ( \langle\hat\nabla_V\tilde V+A_V\tilde V,\tilde T  \rangle\tilde T+G(V,\tilde V)Y )\notag\\
    =&   d\varphi(\nabla^1_V \tilde V) + \lambda^2G(V,\tilde V)\varphi.
\end{align}
Hence $\varphi$ can be interpreted as a hyperbolic centroaffine immersion contained in an $(n_1+1)$-dimensional vector subspace of $\bbr^{n+1}$ with induced connection $\nabla^1$ and centroaffine metric $g^1=\lambda^2 G$.
Similarly,
\begin{align}
   D_{W}d\psi(\tilde W)=&  D_W \tilde W \notag\\
    =&\nabla^2_W \tilde W+ \langle\nabla_W\tilde W,\tilde T  \rangle\tilde T+\nabla^1_W\tilde W+G(W,\tilde W)Y \notag\\
    =&d\psi(\nabla^2_W \tilde W) +G(W,\tilde W)Y.
\end{align}
Hence $\psi$ can be interpreted as a Calabi hypersurface contained in an $(n_2+1)$-dimensional vector subspace of $\bbr^{n+1}$ with induced connection $\nabla^2$ and Calabi metric $g^2=G$.

As both subspaces are complementary and $Y$ is the Calabi normalization of $\psi$, we may assume that up to a linear transformation $\phi\in SA(n+1)$, the $(n_1+1)$-dimensional subspace is spanned by the first $(n_1+1)$ coordinates of
$\bbr^{n+1}$, whereas  the $(n_2+1)$-dimensional subspace is spanned by the last $(n_2+1)$ coordinates of $\bbr^{n+1}$.

Solving (3.56) for the immersion $x$, we have
\be x=e^t\varphi+\psi-\frac{t}{\lambda^2}Y-\frac{1}{\lambda^2}Y. \ee
Thus, $x$ is Calabi affine equivalent to
\be x=(e^t\varphi,\psi)-\frac{t}{\lambda^2}Y. \ee
We have complete the proof of Theorem 3.5.\hfill $\Box$

In Theorem 3.5, if additionally $M$ has parallel cubic form $\hat\nabla A=0$, then by
\be
\hat R(V,W)A(\tilde T, \tilde T)= 2\lambda A_ {\tilde T}(A_W V),
\ee
and
\be
  \hat R(V,W)A(\tilde T, \tilde T)= \lambda^2 A_W V,
\ee
we deduce that $A_WV=0$. Now we apply Theorem 3.5 to prove
the following theorem

\begin{thm}\label{thm3.6}
Let $x:M^n\rightarrow \bbr^{n+1}$ be a Calabi hypersurface. Assume that $\hat\nabla A=0$ and there exist orthogonal distributions $\mathcal {D}_1$ (of dimension 1, spanned by a unit vector filed $\tilde T$), $\mathcal {D}_2$ (of dimension $n_1$),  $\mathcal {D}_3$ (of dimension $n_2$) and  a nonzero constant $\lambda$ such that
$$A(\tilde T,\tilde T)=\lambda\tilde T,\quad A(\tilde T,V)=\lambda V, \quad  A(\tilde T,W)=0,\quad
  \forall V\in \mathcal {D}_2,\; W\in \mathcal {D}_3.$$
Then $x:M^n\rightarrow \bbr^{n+1}$ can be locally decomposed as the Calabi product of a locally strongly convex hyperbolic centroaffine hypersurface $\varphi: M_1^{n_1}\rightarrow \bbr^{n_1+1}$ with parallel cubic form and a Calabi hypersurface $\psi:M_2^{n_2}\rightarrow \bbr^{n_2+1}$ with parallel cubic form.
\end{thm}
\section{A typical basis}
\subsection{A typical basis and known results}
Now, we fix a point $p\in M^n$. For subsequent purpose, we will review the well known construction of a typical orthonormal basis for $T_p M^n$, which was introduced by Ejiri and has been widely applied, and proved to be very useful for various situations, see e.g., \cite{HLSV}, \cite{LV} and \cite{CHM}. The idea is to construct from the $(1, 2)$ tensor $A$ a self adjoint operator at a point; then one extends the eigenbasis to a local field.
Let $p \in M^n$ and $U_p M^n = \{v \in T_pM^n \;|\; G(v, v) = 1\}$. Since $M^n$ is locally strong convex, $U_pM^n$ is compact. We define a function $F$ on $U_pM^n$ by $F (v) = A(v,v,v)$.
Then there is an element $e_1\in U_pM^n$ at which the function $F(v)$ attains an absolute maximum, denoted by  $\mu_1$. Then we have the following lemmas

\begin{lem}\label{lemma4.1}\cite{LSZH}  There exists an orthonormal basis $\{e_1,\cdots,e_n\}$ of $T_pM^n$ such that the following hold: \\
(i) $A(e_1,e_i,e_j)=\mu_i\delta_{ij}$, for $i=1,\cdots,n$. \\
(ii) $\mu_1\geq2\mu_i,$ for $i\geq2$. If $\mu_1=2\mu_i$, then $A(e_i,e_i,e_i) =0$.
\end{lem}
\begin{lem}\label{lemma4.2}\cite{XL1} Let $M^n$ be a Calabi hypersurface with parallel Fubini-Pick form. Then, for every point $p\in M^n$, there exists an orthonormal basis $\{e_j\}_{1\leq j\leq n}$ of $T_pM^n$ (if necessary, we rearrange the order), satisfying $A(e_1,e_j)=\mu_je_{j}$, and there exists a number $i$, $0\leq i\leq n$, such that
$$\mu_2=\mu_3=\cdots=\mu_i=\frac{1}{2}\mu_1;\quad \mu_{i+1}=\cdots=\mu_n=0.$$
\end{lem}
Therefore, for a strictly convex Calabi hypersurface with parallel Fubini-Pick form, we have to deal with $(n+1)$ cases as follows:

\noindent\textbf{Case} $\mathfrak{C}_0.$  $\mu_1=0$.

\noindent\textbf{Case} $\mathfrak{C}_1.$ $\mu_1>0;\mu_2=\mu_3=\cdots=\mu_n=0.$

\noindent\textbf{Case }$\mathfrak{C}_i. $ $\mu_2=\mu_3=\cdots=\mu_i=\frac{1}{2}\mu_1>0;\quad \mu_{i+1}=\cdots=\mu_n=0 \quad \text{for}\;\; 2\leq i\leq n-1.$

\noindent\textbf{Case} $\mathfrak{C}_n. $ $\mu_2=\mu_3=\cdots=\mu_n= \frac{1}{2}\mu_1>0.$

When working at the point $p\in M^n$, we will always assume that an orthonormal basis is chosen such that Lemma 4.1 is satisfied.
About the Cases $\mathfrak{C}_0$, $\mathfrak{C}_n$, we proved the following conclusions in \cite{XL1}.
\begin{lem}\label{lemma4.3}\cite{XL1} If the Case $\mathfrak{C}_0$  occurs, then $M^n$ is an open part of elliptic paraboloid.
\end{lem}
\begin{lem}\label{lemma4.4}\cite{XL1} The Case $\mathfrak{C}_n$ does not occur.
\end{lem}

Applying of Theorem 3.6 to the case $n_1=0$, we have
\begin{cor}
If Case $\mathfrak{C}_1$ occurs, then $M^n$ can be obtained as the Calabi product of an $(n-1)$-dimensional Calabi hypersurface in $\bbr^n$ with parallel cubic form and a point.
\end{cor}

From Theorem 3.3 above and (3.10) in \cite{XL2}, we get
\begin{cor} If Case $\mathfrak{C}_{n-1}$ occurs, then $M^n$ can be obtained as the Calabi product a locally strongly convex hyperbolic centroaffine hypersurface $\varphi: M_1^{n-1}\rightarrow \bbr^{n}$ with parallel cubic form and a point.
\end{cor}

In fact, by induction we obtain its explicit expression (see also (2) in example 3.1).
\begin{thm}\cite{XL2}
If the Case $\mathfrak{C}_{n-1}$ occurs, then $M^n$ is Calabi affine equivalent to an open part of the hypersurface $$ x_{n+1}=\frac{ (n-1)(n-2)}{2R}\ln(x_1^2-(x_2^2+\cdots+x_n^2)), $$
where $R$ is the scalar curvature of $M^n$.
\end{thm}

\section{Discussions in terms of a Typical Basis}

Now, we consider the Cases $\{\mathfrak{C}_m\}_{2\leq m\leq n-2}$. In this section, we follow closely the same procedure as in Cheng-Hu-Moruz paper \cite{CHM} and \cite{HLV}.

\subsection{Intermediary Cases $\{\mathfrak{C}_m\}_{2\leq m\leq n-2}$ and Isotropic Mapping $L$}

In these cases, we denote by $\mathcal{D}_{2}$ and $\mathcal{D}_{3}$ the two subspaces of $T_{p} M$:
$$
\mathcal{D}_{2}=\operatorname{span}\left\{e_{2}, \ldots, e_{m}\right\} \text { and } \mathcal{D}_{3}=\operatorname{span}\left\{e_{m+1}, \ldots, e_{n}\right\}.
$$
First of all, we have the following
\begin{lem}\label{lemma5.1}
Associated with the direct sum decomposition $T_pM=\mathcal{D}_1\oplus \mathcal{D}_2\oplus \mathcal{D}_3$, where $\mathcal{D}_1=\operatorname{span}\{e_1\}$, the following relations hold:
\begin{enumerate}
  \item [(i)]$A(e_1,v)=\frac{1}{2}\mu_1v,\quad A(e_1,w)=0, \quad \forall\ v\in \mathcal{D}_2,\ w\in \mathcal{D}_3$.
  \item [(ii)] $A(v_1,v_2)-\frac{1}{2}\mu_1 G(v_1,v_2) e_1\in \mathcal{D}_3, \quad \forall\ v_1,v_2\in \mathcal{D}_2$.
  \item [(iii)]$A(v,w)\in \mathcal{D}_2, \quad \forall\ v\in \mathcal{D}_2,\ w\in \mathcal{D}_3$.
\end{enumerate}
\end{lem}
\textbf{Proof.}
By definition we have (i). The claim (ii) follows from (ii) of Lemma \ref{lemma4.1}. By $\hat\nabla A=0$, we know that the curvature operator of Levi-Civita connection $\hat R$ and Fubini-Pick tensor $A$ satisfy
\be\label{5.1} \hat{R}(X,Y)A(Z,W)=A(\hat{R}(X,Y)Z,W)+A(Z,\hat{R}(X,Y)W).\ee
In order to prove the third claim, we take $X=v \in \mathcal{D}_{2}$, $Y=w \in \mathcal{D}_{3}$ and $Z=W=e_{1}$ to obtain that
$$
\lambda_{1} \hat{R}(v, w) e_{1}=2 A\left(\hat{R}(v, w) e_{1}, e_{1}\right).
$$
Thus we have $\hat{R}(v, w) e_{1} \in \mathcal{D}_{2}$.
On the other hand, a direct calculation gives
$$
\hat{R}(v, w) e_{1}=-A_{v} A_{w} e_{1}+A_{w} A_{v} e_{1}=\frac{1}{2} \mu_{1} A_{v} w.
$$
Therefore, $A_{v} w \in \mathcal{D}_{2}$.   \hfill $\Box$

With the remarkable conclusions of Lemma \ref{lemma5.1}, similar to that \cite{HLV},  we can now introduce a bilinear map $L: \mathcal{D}_{2} \times \mathcal{D}_{2} \rightarrow \mathcal{D}_{3}$, defined by
$$
L\left(v_{1}, v_{2}\right):=A_{v_{1}} v_{2}-\frac{1}{2} \mu_{1} G\left(v_{1}, v_{2}\right) e_{1}, \ \ \forall v_{1}, v_{2} \in \mathcal{D}_{2}.
$$

The following lemmas show that the operator $L$ enjoys remarkable properties and it becomes an important tool for exploring information of the difference tensor.
\begin{lem}\label{lemma5.2}
The bilinear map $L$ is isotropic in the sense that
\be\label{5.2}
G(L(v,v),L(v,v))=\fr{1}{4}\mu_1^2 G(v,v)^2,\quad \forall\ v\in \mathcal{D}_2.\ee
Moreover, linearizing (\ref{5.2}), it follows for any $v_1,v_2,v_3,v_4\in \mathcal{D}_2$ that
\begin{align}\label{5.3}
G(L&(v_1,v_2),L(v_3,v_4))+G(L(v_1,v_3),L(v_2,v_4))+G(L(v_1,v_4),L(v_2,v_3))    \notag   \\
=&\fr{1}{4}\mu_1^2 \left(G(v_1,v_2)G(v_3,v_4)+G(v_1,v_3)G(v_2,v_4)+G(v_1,v_4)G( v_2,v_3)\right).
\end{align}
\end{lem}
\textbf{Proof.}
We use (\ref{5.1}) and take $X=e_{1}$ and $Y=v_{1}, Z=v_{2}, W=v_{3}$ in $\mathcal{D}_{2} .$ By the definition of $L$, it follows immediately that
\begin{align*}
LHS&=R(e_1,v_1)A(v_2,v_3)=R(e_1,v_1)(L(v_2,v_3)+\fr{1}{2}\mu_1\langle  v_2,v_3\rangle e_1)  \\
&=A_{v_1}A_{e_1}L(v_2,v_3)-A_{e_1}A_{v_1}L(v_2,v_3)+\fr{1}{2}\mu_1\langle v_2,v_3\rangle (A_{v_1}A_{e_1}e_1-A_{e_1}A_{v_1}e_1)     \\
&=-\fr{1}{2}\mu_1A(L(v_2,v_3),v_1)+\fr{1}{8}\mu_1^3\langle v_2,v_3\rangle v_1 ;   \\
RHS&=A(R(e_1,v_1)v_2,v_3)+A(R(e_1,v_1)v_3,v_2),   \\ &=A_{v_3}(A_{v_1}A_{e_1}v_2-A_{e_1}A_{v_1}v_2)+A_{v_2}(A_{v_1}A_{e_1}v_3-A_{e_1}A_{v_1}v_3)\\
&=\fr{1}{2}\mu_1\left(A(L(v_1,v_2),v_3)+A(L(v_1,v_3),v_2)\right)-\fr{1}{8}\mu_1^3(\langle v_1,v_2\rangle v_3+\langle v_1,v_3\rangle v_2).
\end{align*}
Hence
\be\label{5.4}
A(L(v_1,v_2),v_3)+A(L(v_1,v_3),v_2)+A(L(v_2,v_3),v_1)=\fr{1}{4}\mu_1^2(\lang v_1,v_2\rang v_3+\lang v_1,v_3\rang v_2+\lang v_2,v_3\rang v_1).
\ee
Taking the product of (\ref{5.4}) with $v_4\in \mathcal{D}_2$ and by the fact
$$\lang A(L(v_1,v_2),v_3),v_4\rang=\lang L(v_1,v_2),A(v_3,v_4)\rang=\lang L(v_1,v_2),L(v_3,v_4)\rang,$$
we can obtain (\ref{5.3}). By taking $v_1=v_2=v_3=v_4=v$, we get (\ref{5.2}).
\hfill $\Box$

Since $L: \mathcal{D}_{2} \times \mathcal{D}_{2} \rightarrow \mathcal{D}_{3}$ is isotropic, we see from (\ref{5.2}) that, if $\operatorname{dim} \mathcal{D}_{2} \geq$ 1, then the image space of $L$ has positive dimension, i.e. $\operatorname{dim}(\operatorname{Im} L) \geq 1$ Moreover, the following well-known properties hold, see \cite{HLSV} or\cite{HLV}.

\begin{lem}\label{lemma5.3}
If $\operatorname{dim} \mathcal{D}_{2} \geq 1$, for orthonormal vectors $v_1, v_2, v_3$, and $ v_4\in \mathcal{D}_2$, there hold
\begin{align}
\label{5.5}  &\lang L(v_1,v_1),L(v_1,v_2)\rang=0,\\
\label{5.6}
&\lang L(v_1,v_1),L(v_2,v_2)\rang+2\lang L(v_1,v_2),L(v_1,v_2)\rang=\fr{1}{4}\mu_1^2,\\
\label{5.7}&\lang L(v_1,v_1),L(v_2,v_3)\rang+2\lang L(v_1,v_2),L(v_1,v_3)\rang=0,\\
\label{5.8}&\lang L(v_1,v_2),L(v_3,v_4)\rang+\lang L(v_1,v_3),L(v_2,v_4)\rang+\lang L(v_1,v_4),L(v_2,v_3)\rang=0.
\end{align}
\end{lem}

\begin{lem}\label{lemma5.4}
In Cases $\left\{\mathfrak{C}_{m}\right\}_{2 \leq m \leq n-2}$, if it occurs that $\operatorname{Im} L \neq \mathcal{D}_{3}$, then for any $v_{1}, v_{2} \in \mathcal{D}_{2}$ and $w \in \mathcal{D}_{3}$ with $w \perp \operatorname{Im} L$, we have
\be\label{5.9} A(L(v_1,v_2),w)=0. \ee
\end{lem}
\textbf{Proof.}
For every $v \in \mathcal{D}_{2}$ and $w \perp \operatorname{Im} L$, we apply (iii) of Lemma \ref{lemma5.1} to obtain
\begin{align*}
A(v,w)&=\sum_{i=2}^m\lang A(v,w),e_i\rang e_i=\sum_{i=2}^m\lang A(v,e_i),w\rang e_i \\
      &=\sum_{i=2}^m\lang L(v,e_i),w\rang e_i=0, \\
\hat R(e_1,v)w&=-[A_{e_1},A_{v}]w=0.
\end{align*}
Then, for $v_{1}, v_{2}$ and $w$ as in the assumptions, the following equation
$$\hat R(e_1,v_1)A(v_2,w)=A(\hat R(e_1,v_1)v_2,w)+A(v_2,\hat R(e_1,v_1)w)$$
becomes equivalent to $A(\hat R(e_1,v_1)v_2,w)=0$. On the other hand, direct calculation gives that
\begin{align*}
\hat R(e_1,v_1)v_2=A_{v_1}A_{e_1}v_2-A_{e_1}A_{v_1}v_2
=\fr{1}{2}\mu_1L(v_1,v_2)-\fr{1}{4}\mu_1^2\lang v_1,v_2\rang e_1.
\end{align*}
Then (\ref{5.9}) immediately follows.  \hfill $\Box$

\begin{lem}\label{lemma5.5}
 In Cases $\left\{\mathfrak{C}_{m}\right\}_{2 \leq m \leq n-2}$, let $v_{1}, v_{2}, v_{3}, v_{4} \in \mathcal{D}_{2}$ and $\left\{u_{1}, \ldots, u_{m-1}\right\}$ be an orthonormal basis of $\mathcal{D}_{2}$, then we have
\begin{align}\label{5.10}
A(L(v_1,v_2),L(v_3,v_4))=&\sum_{i=1}^{m-1}\lang L(v_1,u_i),L(v_3,v_4)\rang L(v_2,u_i)
\notag\\
 &+\sum_{i=1}^{m-1}\lang L(v_2,u_i),L(v_3,v_4)\rang L(v_1,u_i).
\end{align}
\end{lem}
\textbf{Proof.}
 By (\ref{5.1}), we have, for $v_{1}, v_{2}, v_{3}, v_{4} \in \mathcal{D}_{2}$, that
\begin{align}\label{5.11}
\hat R(e_1,v_1)A(v_2,L(v_3,v_4))=A(\hat R(e_1,v_1)v_2,L(v_3,v_4))+A(v_2,\hat R(e_1,v_1)L(v_3,v_4)).
\end{align}
For $v,v_1,v_2\in \mathcal{D}_2$ and $w\in \mathcal{D}_3$, a direct calculation shows that
\begin{align*}
\hat R(e_1,v)w&=-\fr{1}{2}\mu_1A(v,w),   \\
\hat R(e_1,v_1)v_2&=\fr{1}{2}\mu_1A_{v_1}v_2-A_{e_1}A_{v_1}v_2
=\fr{1}{2}\mu_1L(v_1,v_2)-\fr{1}{4}\mu_1^2\lang v_1,v_2\rang e_1.
\end{align*}
By $A(v_2, L(v_3,v_4))\in \mathcal{D}_2$, we can write
$$A(v_2, L(v_3,v_4))=\sum_{i=1}^{m-1}\lang L(u_i,v_2),L(v_3,v_4)\rang u_i.$$
Now, a direct calculation on both sides of (\ref{5.11}) shows that
\begin{align*}
LHS=&\sum_{i=1}^{m-1}\lang L(u_i,v_2),L(v_3,v_4)\rang \hat R(e_1,v_1)u_i  \\
   =&\fr{1}{2}\mu_1\sum_{i=1}^{m-1}\lang L(u_i,v_2),L(v_3,v_4)\rang L(v_1,u_i) -\fr{1}{4}\mu_1^2\lang L(v_1,v_2),L(v_3,v_4)\rang e_1,
\end{align*}
\begin{align*}
RHS=&A\left(\fr{1}{2}\mu_1L(v_1,v_2)-\fr{1}{4}\mu_1^2\lang v_1,v_2\rang e_1,L(v_3,v_4)\right)-\fr{1}{2}\mu_1 A\left(v_2,A(v_1,L(v_3,v_4))\right) \\
   =&\fr{1}{2}\mu_1 A\left(L(v_1,v_2),L(v_3,v_4)\right)-\fr{1}{2}\mu_1 \sum_{i=1}^{m-1}\lang L(u_i,v_1),L(v_3,v_4)\rang\left(L(u_i,v_2)+\fr{1}{2}\mu_1\lang u_i,v_2\rang e_1\right)   \\
   =&\fr{1}{2}\mu_1 A\left(L(v_1,v_2),L(v_3,v_4)\right)-\fr{1}{4}\mu_1^2\lang L(v_1,v_2),L(v_3,v_4)\rang e_1   \\
   &-\fr{1}{2}\mu_1 \sum_{i=1}^{m-1}\lang L(u_i,v_1),L(v_3,v_4)\rang L(u_i,v_2).
\end{align*}
One can immediately get (\ref{5.10}) from these computations.    \hfill $\Box$

We note that (\ref{5.10}) has very important consequences which will be used in sequel sections. For example, we have

\begin{lem}\label{lem4.12}
For Case $\mathfrak{C}_{m}$ with $m \geq 3$, let $\left\{u_{1}, \ldots, u_{m-1}\right\}$ be an orthonormal basis of $\mathcal{D}_{2}$, then for $p \neq j$, we have
\be\label{5.12}
0=\left(\fr{1}{2}\mu_1^2-4\lang L(u_j,u_p),L(u_j,u_p)\rang\right)L(u_j,u_p)
-\sum_{i\neq p}4\lang L(u_j,u_i),L(u_j,u_p)\rang L(u_j,u_i).\ee
In particular, if $L\left(u_{1}, u_{2}\right) \neq 0$ and $L\left(u_{1}, u_{i}\right)$ is orthogonal to $L\left(u_{1}, u_{2}\right)$ for all $i \neq 2$, then
\be\label{5.13}
\lang L(u_1,u_2),L(u_1,u_2)\rang=\fr{1}{8}\mu_1^2=:\tau.\ee
\end{lem}
\textbf{Proof.}
By (\ref{5.10}), interchanging the couples of indices $\{1,2\}$ and $\{3,4\}$, we find the following condition:
\begin{align}\label{5.14}
0=&\sum_{i=1}^{m-1}\lang L(v_1,u_i),L(v_3,v_4)\rang L(v_2,u_i)+\sum_{i=1}^{m-1}\lang L(v_2,u_i),L(v_3,v_4)\rang L(v_1,u_i)    \notag\\
  &-\sum_{i=1}^{m-1}\lang L(v_3,u_i),L(v_1,v_2)\rang L(v_4,u_i)-\sum_{i=1}^{m-1}\lang L(v_4,u_i),L(v_1,v_2)\rang L(v_3,u_i).
\end{align}
Taking $v_2=v_3=v_4=u_j\neq u_p=v_1$ in (\ref{5.14}) and using also the isotropy condition, we have
\begin{align*}
0=&\lang L(u_p,u_p),L(u_j,u_j)\rang L(u_j,u_p)+\sum_{i\neq p}\lang L(u_i,u_p),L(u_j,u_j)\rang L(u_i,u_j)   \\
&+\lang L(u_j,u_j),L(u_j,u_j)\rang L(u_j,u_p)-2\lang L(u_j,u_p),L(u_j,u_p)\rang L(u_j,u_p)  \\
&-2\sum_{i\neq p}\lang L(u_i,u_j),L(u_j,u_p)\rang L(u_i,u_j)   \\
=&\left(\fr{1}{2}\mu_1^2-4\lang L(u_j,u_p),L(u_j,u_p)\rang\right)L(u_j,u_p)-4\sum_{i\neq p}\lang L(u_i,u_j),L(u_p,u_j)\rang L(u_i,u_j).
\end{align*}
If $j=1$ and $p=2$ in (\ref{5.12}), we obtain (\ref{5.13}).    \hfill $\Box$

\subsection{ The Mapping $P_{v}: \mathcal{D}_{2} \rightarrow \mathcal{D}_{2}$ with Unit Vector $v \in \mathcal{D}_{2}$}

We now define for any given unit vector $v \in \mathcal{D}_{2}$ a linear map $P_{v}: \mathcal{D}_{2} \rightarrow \mathcal{D}_{2}$ by
$$
P_{v} \bar{v}=A_{v} L(v, \bar{v}), \quad \forall \bar{v} \in \mathcal{D}_{2}.
$$
It is easily seen that $P_{v}$ is a symmetric linear operator satisfying
$$
G\left(P_{v} \bar{v}, v^{\prime}\right)=G\left(L(v, \bar{v}), L\left(v, v^{\prime}\right)\right)=G\left(P_{v} v^{\prime}, \bar{v}\right),
$$
for any $\bar{v}, v^{\prime} \in \mathcal{D}_{2}$. Moreover, we have

\begin{lem}\label{lemma5.7}
For any unit vector $v \in \mathcal{D}_{2}$, the operator $P_{v}: \mathcal{D}_{2} \rightarrow \mathcal{D}_{2}$ has $\sigma=\frac{1}{4} \mu_1^2$ as an eigenvalue with eigenvector $v .$ In the orthogonal complement $\{v\}^{\perp}$ of $\{v\}$ in $\mathcal{D}_{2}$ the operator $P_{v}$ has at most two eigenvalues, namely 0 and $\tau$, defined as in (\ref{5.13}).
\end{lem}

\textbf{Proof.}
 By (\ref{5.2}), we have
$$ G\left(P_{v} v, v\right)=G(L(v, v), L(v, v))=\frac{1}{4} \mu_1^2.$$
Taking $v^{\prime} \perp v$, we get
$$G\left(P_{v} v, v^{\prime}\right)=G\left(L\left(v, v^{\prime}\right), L(v, v)\right)=0.$$
The above relations imply that $P_{v} v=\fr{1}{4}\mu_1^2v$.

Next, we take an orthonormal basis $\left\{u_{1}, \ldots, u_{m-1}\right\}$ of $\mathcal{D}_{2}$ consisting of eigenvectors of $P_{v}$ such that $P_{v} u_{i}=\sigma_{i} u_{i}, i=1, \ldots, m-1$, with $u_{1}=v$ and $\sigma_{1}=\sigma .$ We take the inner product of (\ref{5.12}) with $L\left(u_{1}, u_{p}\right)$ for $j=1$ and any $p \geq 2 .$ We obtain that
$$0=\left(\fr{1}{2}\mu_1^2-4\lang L(u_1,u_p),L(u_1,u_p)\rang\right)\lang L(u_1,u_p),L(u_1,u_p)\rang,$$
which implies that
$$\left(\tau-\lang P_{u_1}u_p,u_p\rang\right)\lang P_{u_1}u_p,u_p\rang=0.$$
The previous equation shows the remaining assertion.     \hfill $\Box$

In the following we denote by $V_{v}(0)$ and $V_{v}(\tau)$ the eigenspaces of $P_{v}$ (in the orthogonal complement of $\{v\}$ ) with respect to the eigenvalues 0 and $\tau$, respectively.

\begin{lem}\label{lemma5.8}
Let $v,u\in \mathcal{D}_{2}$ be two unit orthogonal vectors. Then the following statements are equivalent:
\begin{enumerate}
  \item [(i)]  $u\in V_v(0)$.
  \item [(ii)]  $L(u,v)=0$.
  \item [(iii)]  $L(u,u)=L(v,v)$.
  \item [(iv)]  $v\in V_u(0)$. \newline Moreover, any of the previous statements implies that
  \item [(v)]\ \ $P_v=P_u$ on $\{u,v\}^\perp$.
\end{enumerate}
\end{lem}

\textbf{Proof.}
As $G\left(P_{v} u, u\right)=G(L(v, u), L(v, u))=G\left(P_{u} v, v\right)$, the equivalence of (i), (ii) and (iv) follows immediately. As $u$ and $v$ are orthogonal, (\ref{5.2}) and (\ref{5.6}) imply that
$$\lang L(v,v)-L(u,u),L(v,v)-L(u,u)\rang=\fr{1}{2}\mu_1^2-2\lang L(v,v),L(u,u)\rang=4\lang L(u,v),L(u,v)\rang.$$
It follows that (ii) is equivalent to (iii).

Now we assume that (i)-(iv) are satisfied. In order to prove (v), we see that the space spanned by $\{u, v\}$ is invariant by $P_{v}$ and $P_{u}$, also its orthogonal complement is invariant. By taking $v_1,v_2\in\{u,v\}^\perp$ and using (\ref{5.3}), we find
\begin{align*}
\lang P_vv_1,v_2\rang=&\lang L(v,v_1),L(v,v_2)\rang=\fr{1}{8}\mu_1^2\lang v_1,v_2\rang-\fr{1}{2}\lang L(v,v),L(v_1,v_2)\rang   \\
                     =&\fr{1}{8}\mu_1^2\lang v_1,v_2\rang-\fr{1}{2}\lang L(u,u),L(v_1,v_2)\rang=\lang P_uv_1,v_2\rang,
\end{align*}
and
\begin{align*}
\lang P_vv_1,u\rang=&\lang P_uv_1,u\rang =0= \lang P_vv_1,v\rang=\lang P_uv_1,v\rang.
\end{align*}
This completes the proof.
 \hfill $\Box$

\begin{lem}\label{lemma5.9}
Let $v, \td{v} \in \mathcal{D}_{2}$ be two unit orthogonal vectors, then
\be\label{5.15}
G(L(v, \td{v}), L(v, \td{v}))=\tau  \ee
holds if and only if $\td{v} \in V_{v}(\tau) .$ Moreover, if we assume $u \in V_{v}(0)$ and the equality in (\ref{5.15}) holds, then $u \in V_{\td{v}}(\tau)$.
\end{lem}

\textbf{Proof.}
If $\td{v}\in V_v(\tau)$, then $\lang L(v,\td{v}),L(v,\td{v})\rang=\lang P_v\td{v},\td{v}\rang=\tau$.

Conversely, if $\lang L(v,\td{v}),L(v,\td{v})\rang=\tau$, we should consider the following three cases:

(i) $V_v(0)=\emptyset$. By $v\perp\td{v}$, it follows that $\td{v}\in V_v(\tau)$.

(ii)$V_v(\tau)=\emptyset$. In this case, Lemma \ref{lemma5.7} shows that $\td{v}\in V_v(0)$. By Lemma \ref{lemma5.8}, we have $\lang L(v,\td{v}),L(v,\td{v})\rang=0$. This is a contradiction.

(iii) $V_v(0)\neq\emptyset$ and $V_v(\tau)\neq\emptyset$. We can write
$$\td{v}=cos\theta v_0+sin\theta v_1, \quad \lang v_0,v_0\rang=\lang v_1,v_1\rang=1,$$
where $v_0\in V_v(0)$ and $v_1\in V_v(\tau)$. Then we get
$$\tau=\lang L(v,\td{v}),L(v,\td{v})\rang=sin^2\theta\tau,$$
which means that $cos\theta=0$. Hence $\td{v}\in V_v(\tau)$.

Taking unit vector $u\in V_v(0)$, we have $L(u,u)=L(v,v)$. By (\ref{5.6}), we get
\begin{align*}
\lang L(\td{v},u),L(\td{v},u)\rang=&\fr{1}{8}\mu_1^2-\fr{1}{2}\lang L(\td{v},\td{v}),L(u,u)\rang\\
=&\fr{1}{8}\mu_1^2-\fr{1}{2}\lang L(\td{v},\td{v}),L(v,v)\rang=\lang L(\td{v},v),L(\td{v},v)\rang=\tau.
\end{align*}
Applying the first assertion of Lemma \ref{lemma5.9}, one get $u\in V_{\td{v}}(\tau)$.
 \hfill $\Box$

\begin{lem}\label{lemma5.10}
 Let $v_{1}, v_{2}, v_{3} \in \mathcal{D}_{2}$ be orthonormal vectors satisfying $v_{1}, v_{2} \in$ $V_{v_{3}}(\tau)$, then for any vector $v \in \mathcal{D}_{2}$, we have $G\left(L\left(v_{1}, v_{2}\right), L\left(v, v_{3}\right)\right)=0$.
 \end{lem}

\textbf{Proof.} Using the linearity of the assertion with $v$, we may assume that $v$ is an eigenvector of $P_{v_{3}}$. Let $\left\{u_{1}, \ldots, u_{m-1}\right\}$ be an orthonormal basis of $\mathcal{D}_{2}$ consisting of eigenvectors of $P_{v_{3}}$ such that $u_{1}=v_{1}, u_{2}=v_{2}$ and $u_{3}=v_{3}$. We now use (\ref{5.14}) for $v_{3}=v_{4}$ to obtain
\begin{align}\label{5.16}
0=&\sum_{i=1}^{m-1}\lang L(v_1,u_i),L(v_3,v_3)\rang L(v_2,u_i)+\sum_{i=1}^{m-1}\lang L(v_2,u_i),L(v_3,v_3)\rang L(v_1,u_i)    \notag\\
  &-2\sum_{i=1}^{m-1}\lang L(v_3,u_i),L(v_1,v_2)\rang L(v_3,u_i).
\end{align}
On the other hand, from (\ref{5.5})-(\ref{5.7}), we have
$$
\begin{aligned}
&G\left(L\left(v_{1}, u_{i}\right), L\left(v_{3}, v_{3}\right)\right)=G\left(L\left(v_{2}, u_{j}\right), L\left(v_{3}, v_{3}\right)\right)=0, i \neq 1, j \neq 2, \\
&G\left(L\left(v_{1}, v_{1}\right), L\left(v_{3}, v_{3}\right)\right)=G\left(L\left(v_{2}, v_{2}\right), L\left(v_{3}, v_{3}\right)\right)=0.
\end{aligned}
$$
Inserting the above into (\ref{5.16}), we obtain
\be\label{5.17}
0=\sum_{i=1}^{m-1} G\left(L\left(v_{3}, u_{i}\right), L\left(v_{1}, v_{2}\right)\right) L\left(u_{i}, v_{3}\right).   \ee

Taking the inner product of (\ref{5.17}) with $ L\left(v_{1}, v_{2}\right)$, we complete the proof of Lemma \ref{lemma5.10}. \hfill $\Box$

\subsection{ Direct Sum Decomposition for $\mathcal{D}_{2}$}

For our purpose, a crucial matter is to introduce a direct sum decomposition for $\mathcal{D}_{2}$ based on the preceding Lemmas. First, pick any unit vector $v_{1} \in \mathcal{D}_{2}$ and recall that $\tau=\frac{1}{8}\mu_{1}^2$, then by Lemma \ref{lemma5.7}, we have a direct sum decomposition for $\mathcal{D}_{2}$ :
$$
\mathcal{D}_{2}=\left\{v_{1}\right\} \oplus V_{v_{1}}(0) \oplus V_{v_{1}}(\tau),
$$
where, here and later on, we denote also by $\{\cdot\}$ the vector space spanned by its elements. If $V_{v_{1}}(\tau) \neq \emptyset$, we take an arbitrary unit vector $v_{2} \in V_{v_{1}}(\tau)$. Then by Lemma \ref{lemma5.9} we have:
$$
v_{1} \in V_{v_{2}}(\tau), \quad V_{v_{1}}(0) \subset V_{v_{2}}(\tau) \text { and } V_{v_{2}}(0) \subset V_{v_{1}}(\tau).
$$
From this we deduce that
$$
\mathcal{D}_{2}=\left\{v_{1}\right\} \oplus V_{v_{1}}(0) \oplus\left\{v_{2}\right\} \oplus V_{v_{2}}(0) \oplus\left(V_{v_{1}}(\tau) \cap V_{v_{2}}(\tau)\right).
$$
If $V_{v_{1}}(\tau) \cap V_{v_{2}}(\tau) \neq \emptyset$, we further pick a unit vector $v_{3} \in V_{v_{1}}(\tau) \cap V_{v_{2}}(\tau)$. Then
$$
\mathcal{D}_{2}=\left\{v_{3}\right\} \oplus V_{v_{3}}(0) \oplus V_{v_{3}}(\tau),
$$
and by Lemma \ref{lemma5.9} we have
$$
v_{1}, v_{2} \in V_{v_{3}}(\tau) ;\quad V_{v_{1}}(0), V_{v_{2}}(0) \subset V_{v_{3}}(\tau)
$$
It follows that
$$\mathcal{D}_{2}=\{v_1\}\oplus V_{v_1}(0)\oplus\{v_2\}\oplus V_{v_2}(0)\oplus\{v_3\}\oplus V_{v_3}(0)\oplus\left(V_{v_1}(\tau)\cap V_{v_2}(\tau)\cap V_{v_3}(\tau)\right).$$
Considering that $\operatorname{dim}\left(\mathcal{D}_{2}\right)=m-1$ is finite, by induction, we get
\begin{prop}\label{prop5.11} In Cases $\left\{\mathfrak{C}_{m}\right\}_{2 \leq m \leq n-2}$, there exists an integer $k_{0}$ and unit vectors $v_{1}, \ldots, v_{k_{0}} \in \mathcal{D}_{2}$ such that
\be\label{5.18}
\mathcal{D}_{2}=\left\{v_{1}\right\} \oplus V_{v_{1}}(0) \oplus \cdots \oplus\left\{v_{k_{0}}\right\} \oplus V_{v_{k_{0}}}(0).  \ee
\end{prop}

In what follows, we will study the decomposition (\ref{5.18}) in more details.

\begin{lem}\label{lemma5.12}
\begin{enumerate}
\item [(i)] For any unit vector $u_1\in\{v_1\}\oplus V_{v_1}(0)$, we have
$$\{v_1\}\oplus V_{v_1}(0)=\{u_1\}\oplus V_{u_1}(0).$$
\item [(ii)] For any orthonormal vectors $u_1, \td{u}_1\in\{v_1\}\oplus V_{v_1}(0)$, we have
$L(u_1, \td{u}_1)=0$.
\end{enumerate}
\end{lem}

\textbf{Proof.}
(i) We first assume the special case that $u_{1} \perp v_{1}$. Then we have $u_{1} \in V_{v_{1}}(0)$ and thus $L\left(u_{1}, v_{1}\right)=0$, hence $v_{1} \in V_{u_{1}}(0) .$ Let $u \in V_{v_{1}}(0)$ and write $u=x_{1} u_{1}+u^{\prime}$ with $u^{\prime} \perp u_{1} .$ By $(\mathrm{v})$ in Lemma \ref{lemma5.8}, we have $P_{u_{1}} u^{\prime}=P_{v_{1}} u^{\prime}=$ $P_{v_{1}}\left(u-x_{1} u_{1}\right)=0 .$ Therefore, $u^{\prime} \in V_{u_{1}}(0)$ and $\left\{v_{1}\right\} \oplus V_{v_{1}}(0) \subset\left\{u_{1}\right\} \oplus V_{u_{1}}(0)$. Similarly, we obtain $\left\{u_{1}\right\} \oplus V_{u_{1}}(0) \subset\left\{v_{1}\right\} \oplus V_{v_{1}}(0)$.

Next we consider the general case in three subcases. (a) If $V_{v_{1}}(0)=\emptyset$, there is nothing to prove. (b) If $\operatorname{dim}\left(V_{v_{1}}(0)\right) \geq 2$, we can take a vector $\tilde{u} \in$ $V_{v_{1}}(0)$ which is orthogonal to both $u_{1}$ and $v_{1}$. (1) By $\td u \perp v_{1}$, $\td u \in V_{v_{1}}(0)$ and the previous result, we get $\left\{v_{1}\right\} \oplus V_{v_{1}}(0)=\left\{\td u\right\} \oplus V_{\td u}(0)$.
(2) $u_1 \in \left\{v_{1}\right\} \oplus V_{v_{1}}(0)=\left\{\td u\right\} \oplus V_{\td u}(0)$ and $u_1 \perp \td u$ show that $u_1 \in V_{\td u}(0)$. Applying the special case then completes the proof. (c) If $\operatorname{dim}\left(V_{v_{1}}(0)\right)=1$, there exists a unit vector $v_{0} \in V_{v_{1}}(0)$ such that $V_{v_{1}}(0)=\left\{v_{0}\right\}$. Denote $u_{1}=\cos \theta v_{1}+\sin \theta v_{0}$. By Lemma \ref{lemma5.8}, we see that
$$
L\left(\cos \theta v_{1}+\sin \theta v_{0}, \cos \theta v_{0}-\sin \theta v_{1}\right)=0,
$$
thus $\cos \theta v_{0}-\sin \theta v_{1} \in V_{u_{1}}(0) .$ Therefore, $\left\{v_{1}\right\} \oplus V_{v_{1}}(0) \subset\left\{u_{1}\right\} \oplus V_{u_{1}}(0)$. If $\left\{v_{1}\right\} \oplus V_{v_{1}}(0) \subsetneq\left\{u_{1}\right\} \oplus V_{u_{1}}(0)$, we have a unit vector $x \in\left\{u_{1}\right\} \oplus V_{u_{1}}(0)$ which is orthogonal to both $u_{1}$ and $v_{1}$. As
$x \perp u_1$ and $x \in \left\{u_{1}\right\} \oplus V_{u_{1}}(0)$, we deduce that $\left\{u_{1}\right\} \oplus V_{u_{1}}(0)=\left\{x\right\} \oplus V_{x}(0)$, hence $\left\{v_{1}\right\} \oplus V_{v_{1}}(0) \subsetneq\left\{x\right\} \oplus V_{x}(0)$. By $x \perp v_1$, we have $v_1 \in V_{x}(0)$. It follows that $x \in V_{v_1}(0)$, which is contradict to $x \notin \left\{v_{1}\right\} \oplus V_{v_1}(0)$.

(ii) From (i) we have that $\left\{v_{1}\right\} \oplus V_{v_{1}}(0)=\left\{u_{1}\right\} \oplus V_{u_{1}}(0)$. As $u_{1}$ and $\td{u}_{1}$ are orthogonal, this implies that $\td{u}_{1} \in V_{u_{1}}(0)$. Consequently, we have $L\left(u_{1}, \td{u}_{1}\right)=0$.\hfill $\Box$

\begin{lem}\label{lemma5.13}
 In the decomposition (\ref{5.18}), if we pick a unit vector $u_{2} \in V_{v_{2}}(0)$, then there exists a unique unit vector $u_{1} \in V_{v_{1}}(0)$ such that $L\left(u_{1}, v_{2}\right)=$ $L\left(v_{1}, u_{2}\right)$ and $L\left(v_{1}, v_{2}\right)=-L\left(u_{1}, u_{2}\right)$.
\end{lem}

 \textbf{Proof.}
Let $u_{1}^{l}, \ldots, u_{p_l}^{l}$ be an orthonormal basis of $V_{v_{l}}(0), 1 \leq l \leq k_{0}$, such that $u_{1}^{2}=u_{2}$. Then
$$
\left\{v_{1}, \ldots, v_{k_{0}}, u_{1}^{1}, \ldots, u_{p_1}^{1}, \ldots, u_{1}^{k_{0}}, \ldots, u_{p_{k_0}}^{k_{0}}\right\}=:\left\{\td{u}_{i}\right\}_{1 \leq i \leq m-1}
$$
forms an orthonormal basis of $\mathcal{D}_{2}$. Now we use (\ref{5.10}) with the vectors $v_{2}, u_{2}, v_{1}, v_{2} .$ As by Lemma \ref{lemma5.8} $L\left(v_{2}, u_{2}\right)=0$ and by our decomposition $v_{1} \in V_{v_{2}}(\tau)$, we obtain
\begin{align*}
0=&A(L(v_2,u_2),L(v_1,v_2))\\
 =&\sum_{i=1}^{m-1}\lang L(v_2,\td{u}_i),L(v_1,v_2)\rang L(u_2,\td{u}_i)+\sum_{i=1}^{m-1}\lang L(u_2,\td{u}_i),L(v_1,v_2)\rang L(v_2,\td{u}_i)\\
 =&\tau L(u_2,v_1)+\sum_{i=1}^{m-1}\lang L(u_2,\td{u}_i),L(v_1,v_2)\rang L(v_2,\td{u}_i).
\end{align*}
Let us take
$$u_1=-\fr{1}{\tau}\sum_{i=1}^{m-1}\lang L(u_2,\td{u}_i),L(v_1,v_2)\rang \td{u}_i,$$
which satisfies $L(u_1,v_2)=L(u_2,v_1)$.

Firstly, we prove $u_1\in V_{v_1}(0)$.\\
(a) $\td{u}_{i} \notin\left\{v_{1}\right\} \oplus V_{v_{1}}(0) \oplus\left\{v_{2}\right\} \oplus V_{v_{2}}(0)$. By $u_2\in V_{v_2}(0)\subset V_{v_1}(\tau)$ and Lemma \ref{lemma5.10}, we have
$$
G\left(L\left(u_{2}, \td{u}_{i}\right), L\left(v_{1}, v_{2}\right)\right)=0. \quad
$$
(b) $\td{u}_{i} \in\left\{v_{2}\right\} \oplus V_{v_{2}}(0)$. Note the fact $u_2\in\left\{v_{2}\right\} \oplus V_{v_2}(0)$, we have

(1) $\td u_i\neq u_2$, which shows that $\td u_i\perp u_2$. From (ii) of Lemma \ref{lemma5.12}, we know $L(u_2,\td u_i)=0$.

(2) $\td u_i= u_2\in V_{v_2}(0)$. By Lemma \ref{lemma5.8}, we get $L(u_2,u_2)=L(v_2,v_2)$.\\
Applying the above relations, we get
$$
G\left(L\left(u_{2}, \tilde {u}_{i}\right), L\left(v_{1}, v_{2}\right)\right)=0. \quad
$$
(c) $\td u_i= v_1$. The fact $v_2\perp u_2$ shows that
$$
G\left(L\left(u_{2}, v_{1}\right), L\left(v_{1}, v_{2}\right)\right)=0.
$$
It follows from $(a-c)$ that $u_{1} \in V_{v_{1}}(0)$.

In order to prove the uniqueness of $u_{1} \in V_{v_{1}}(0)$, suppose that $\bar{u}_{1} \in V_{v_{1}}(0)$ such that $L\left(\bar{u}_{1}, v_{2}\right)=L\left(v_{1}, u_{2}\right)$, then we have $L\left(u_{1}-\bar{u}_{1}, v_{2}\right)=0$. It follows from Lemma \ref{lemma5.8} that $u_{1}-\bar{u}_{1} \in V_{v_{2}}(0).$ On the other hand, we also have $u_{1}-\bar{u}_{1} \in V_{v_{1}}(0)$; so we must have $u_{1}=\bar{u}_{1}$.

Next, we prove $u_1$ is a unit vector. By the following fact
$$
u_1\in V_{v_{1}}(0) \subset V_{v_{2}}(\tau), \quad u_2\in V_{v_{2}}(0) \subset V_{v_{1}}(\tau), \quad G(u_2,u_2)=1,
$$
we have
$$G\left(u_{1}, u_{1}\right) \tau=G(P_{v_2}u_1,u_1)=G\left(L\left(u_{1}, v_{2}\right), L\left(u_{1}, v_{2}\right)\right)=G\left(L\left(v_{1}, u_{2}\right), L\left(v_{1}, u_{2}\right)\right)=\tau.$$
Hence, $u_{1}$ is a unit vector.

In order to prove the fact that $L\left(u_{1}, v_{2}\right)=L\left(v_{1}, u_{2}\right)$ and $L\left(v_{1}, v_{2}\right)=$ $-L\left(u_{1}, u_{2}\right)$ are equivalent, we use (\ref{5.3}) and the Cauchy-Schwarz inequality. In fact, if we first suppose $L\left(u_{1}, v_{2}\right)=L\left(v_{1}, u_{2}\right)$, then applying (\ref{5.3}) we get
$$
G\left(L\left(v_{1}, v_{2}\right),-L\left(u_{1}, u_{2}\right)\right)=G\left(L\left(v_{1}, u_{2}\right), L\left(v_{2}, u_{1}\right)\right)=G\left(L\left(v_{2}, u_{1}\right), L\left(v_{2}, u_{1}\right)\right)=\tau.
$$
On the other hand, as $u_2\in V_{v_2}(0)$ and Lemma \ref{lemma5.12}, we have $\{v_2\}\oplus V_{v_2}(0)=\{u_2\}\oplus V_{u_2}(0)$. So $V_{v_2}(\tau)=V_{u_2}(\tau)$. By $v_1,u_1\in V_{v_2}(\tau)=V_{u_2}(\tau)$, we get
$$\lang L(v_1,v_2),L(v_1,v_2)\rang=\lang L(u_1,u_2),L(u_1,u_2)\rang=\tau.$$
Then, by Cauchy-Schwarz inequality, we have $L(u_1,u_2)=-L(v_1,v_2)$.

The converse can be proved in a similar way.
\hfill $\Box$

To state the next lemma, we denote $V_{l}=\left\{v_{l}\right\} \oplus V_{v_{l}}(0)$ in the decomposition (\ref{5.18}) for each $1 \leq l \leq k_{0} .$ Then we have

\begin{lem}\label{lemma5.14}
 With respect to the decomposition (\ref{5.18}), the following hold.
\begin{enumerate}
  \item [(i)] For any unit vector $a\in V_j$,
\be\label{5.19}
A(L(a,a),L(a,a))=\fr{1}{2}\mu_1^2L(a,a).\ee
  \item [(ii)] For $j\neq l$ and any unit vectors $a\in V_j,\ b\in V_l$,
\begin{align}
\label{5.20}A(L(a,a),L(a,b))=&2\tau L(a,b), \\
\label{5.21}A(L(a,a),L(b,b))=&0, \\
\label{5.22}A(L(a,b),L(a,b))=&\tau (L(a,a)+L(b,b)).
\end{align}
\item [(iii)] For distinct $j, l, q, s$ and any unit vectors $a\in V_j,\ b, b'\in V_l,\ c\in V_q,\ d\in V _s$, where $b$ and $b'$ are orthogonal, the following relations hold
\begin{align}
\label{5.23}A(L(a,b),L(a,c))=&\tau L(b,c), \\
\label{5.24}A(L(a,a),L(b,c))=&0, \\
\label{5.25}A(L(a,b),L(a,b'))=&0,\\
\label{5.26}A(L(a,b),L(c,d))=&0.
\end{align}
\item [(iv)] For distinct $j, l, q$ and orthogonal unit vectors $a_1, a_2\in V_j$, and unit vectors $\ b\in V_l,\ c\in V_q$, it holds
\be\label{5.27} A(L(a_1,b),L(a_2,c))=\tau L(b,c'),\ee
where $c'\in V_q$ is the unique unit vector satisfying $L(a_1,c')=L(a_2,c)$.
\end{enumerate}
\end{lem}

\textbf{Proof.}
Take an orthonormal basis of $\mathcal{D}_{2}$ such that it consists of the orthonormal basis of all $V_{l}, 1 \leq l \leq k_{0}$, the assertions are direct consequences of Lemma \ref{lemma5.5}.

Taking $v_1=v_2=v_3=v_4=a$ in (\ref{5.10}) and $v_1=v_2=v_3=a$, $v_4=u_i$ in (\ref{5.3}), we get
\begin{align*}
A(L(a,a),L(a,a))&=2\sum_{i=1}^{m-1}\lang L(a,u_i),L(a,a)\rang L(a,u_i),\\
\lang L(a,u_i),L(a,a)\rang&=\fr{1}{4}\mu_1^2\lang a,u_i\rang.
\end{align*}
Thus (\ref{5.19}) holds. Similarly, (\ref{5.20})-(\ref{5.24}) and (\ref{5.26}) can be proved.

From the fact $b\perp b'\in V_l$, we get $b'\in V_b(0)$. From Lemma \ref{lemma5.13} we see that there exists a unique unit vector  $a'\in V_{a}(0)$ such that $L(a,b')=L(a',b)$. By taking $v_1=v_3=a,\ v_2=b$ and $v_4=b'$ in (\ref{5.10}), we have
\begin{align*}
A(L(a,b),L(a,b'))=&\sum_{i=1}^{m-1}\lang L(a,u_i),L(a,b')\rang L(b,u_i)
                  +\sum_{i=1}^{m-1}\lang L(b,u_i),L(a,b')\rang L(a,u_i)  \\
                 =&\sum_{i=1}^{m-1}\lang P_ab',u_i\rang L(b,u_i)
                  +\sum_{i=1}^{m-1}\lang P_ba',u_i\rang L(a,u_i)  \\
                 =&\tau L(b,b')+\tau L(a,a')=0.
\end{align*}
This complete the proof of (\ref{5.25}).

From Lemma \ref{lemma5.12}, the fact $c\in V_q$ implies that $V_c=V_q$. By $a_1\perp a_2\in V_j$, we know $a_2\in V_{a_1}(0)$. Lemma \ref{lemma5.13} shows that there exists a unique unit vector  $c'\in V_c(0)$ such that $L(a_1,c')=L(a_2,c)$. By $a_2,c\in V_b(\tau)$ and taking $v_1=a_1,\ v_2=b,\ v_3=a_2$ and $v_4=c$ in (\ref{5.10}), we have
\begin{align*}
A(L(a_1,b),L(a_2,c))=&\sum_{i=1}^{m-1}\lang L(a_1,u_i),L(a_2,c)\rang L(b,u_i)
                  +\sum_{i=1}^{m-1}\lang L(b,u_i),L(a_2,c)\rang L(a_1,u_i)  \\
                 =\sum_{i=1}^{m-1}&\lang L(a_1,u_i),L(a_1,c')\rang L(b,u_i)
                 =\sum_{i=1}^{m-1}\lang P_{a_1}c',u_i\rang L(b,u_i)=\tau L(b,c').
\end{align*}
This complete the proof of (\ref{5.27}).
 \hfill $\Box$

\begin{prop}\label{prop5.15}
In the decomposition (\ref{5.18}), if $k_{0}=1$, then $\operatorname{dim}(\operatorname{Im} L)=1$. If $k_{0} \geq 2$, then $\operatorname{dim} V_{v_{1}}(0)=\cdots=\operatorname{dim} V_{v_{k_{0}}}(0)$ and the dimension which we denote by $\mathfrak{p}$ can only be equal to $0,1,3$ or $7 .$
\end{prop}

\textbf{Proof.} If $k_0=1$, then $\mathcal{D}_{2}=\{v_1\}\oplus V_{v_1}(0)$. Let $\{u_1,\cdots,u_{m-1}\}$ be an orthonormal basis of $\mathcal{D}_{2}$, and $u_1=v_1$. By (ii) of Lemma \ref{lemma5.12} we know that $L(u_i,u_j)=0$, for any $i\neq j$. Lemma \ref{lemma5.8} shows that $L(v_1,v_1)=L(u_1,u_1)=L(u_i,u_i)$ for any $u_i\in V_{v_1}(0),\ i>1$. Thus, for any $a,b\in \mathcal{D}_{2}$, we have
$$L(a,b)=\sum_{i=1}^{m-1}a_ib_iL(u_i,u_i)=\sum_{i=1}^{m-1}a_ib_iL(v_1,v_1).$$
It follows that $L(v_1,v_1)$ is a basis of the image $\operatorname{Im} L$, which implies that
$\operatorname{dim}(\operatorname{Im} L)=1$.

If $k_0\geq2$, as a direct consequence of Lemma \ref{lemma5.13}, we can define a one-to-one linear mapping $\varphi: V_{v_j}(0)\rightarrow V_{v_l}(0)$ for any $j\neq l$, which preserves the length of vectors.

In fact, $\varphi$ is an injective mapping, which means that if $u_2\neq\td{u}_2$, then $\varphi(u_2)\neq\varphi(\td{u}_2)$. Lemma \ref{lemma5.13} implies that for any unit vector $u_2\in V_{v_j}(0)$, then there exists a unique unit vector  $\varphi(u_2)\in V_{v_l}(0)$ such that $L(\varphi(u_2),v_j)=L(u_2,v_l)$. In order to prove the assertion, we assume that $\varphi(u_2)=\varphi(\td{u}_2)$, then $L(\varphi(u_2),v_j)=L(\varphi(\td{u}_2),v_j)$. So we have $L(u_2,v_l)=L(\td{u}_2,v_l)$. It follows from Lemma \ref{lemma5.8} that $u_2-\td{u}_2\in V_{v_l}(0)$. On the other hand, $u_2-\td{u}_2\in V_{v_j}(0)$. Thus, we get $u_2=\td{u}_2$.

We can also define an injective mapping $\phi:V_{v_l}(0)\rightarrow V_{v_j}(0)$ by Lemma \ref{lemma5.13}. Hence $V_{v_j}(0)$ and $V_{v_l}(0)$ are isomorphic and have the same dimension which we denote by $\mathfrak{p}$. To make the following discussion meaningful, we now assume $\mathfrak{p}\geq1$.

Let $\{v_l,u_1^l,\cdots,u_{\mathfrak{p}}^l\}$ be an orthonormal basis of $V_l=\{v_l\}\oplus V_{v_l}(0)$. Note that $u_j^2\in V_{v_2}(0)$ for each $j=1,\cdots,\mathfrak{p}$. Lemma \ref{lemma5.13} shows that there exists a unique unit vector  $u_1\in V_{v_1}(0)$ such that $L(v_1,u_j^2)=L(v_2,u_1)$.

It follows that we can define a linear map $\mathcal{T}_j:V_1\rightarrow V_1$ such that, for any unit vector $v\in V_1\ (V_1=\{v_1\}\oplus V_{v_1}(0)=\{v\}\oplus V_{v}(0))$, the image $\mathcal{T}_j(v)$ satisfies
\be\label{5.28} L(v,u_j^2)=L(v_2,\mathcal{T}_j(v)).\ee
The linear map $\mathcal{T}_j:V_1\rightarrow V_1$ has the following properties:\\
(a) For any $v\in V_1,\ \lang \mathcal{T}_j(v),\mathcal{T}_j(v)\rang=\lang v,v\rang$, i.e., $\mathcal{T}_j$ preserves the length of vectors.\\
(b) For all $v\in V_1$, we have $\mathcal{T}_j(v)\perp v$.\\
(c) $\mathcal{T}_j^2=-id$.\\
(d) For any $j\neq l$, we have $\lang \mathcal{T}_j(v),\mathcal{T}_l(v)\rang=0$ for all $v\in V_1$.

(a) and (b) can be easily seen from Lemma \ref{lemma5.13} and the definition of $\mathcal{T}_j(v)$. We now verify (c) and (d) separately. For any unit vector $v\in V_1$, we have
\be\label{5.29} L(\mathcal{T}_j(v),u_j^2)=L(v_2,\mathcal{T}_j^2(v)).\ee
Using the fact $V_1=\{\mathcal{T}_j(v)\}\oplus V_{\mathcal{T}_j(v)}(0)$ and $u_j^2\in V_{v_2}(0)\subset V_{\mathcal{T}_j(v)}(\tau)$, we have
\begin{align*}
\lang L(u_j^2,\mathcal{T}_j(v)),L(u_j^2,\mathcal{T}_j(v))\rang=\lang L(v_2,\mathcal{T}_j(v)),L(v_2,\mathcal{T}_j(v))\rang=\lang L(v_2,v),L(v_2,v)\rang=\tau.
\end{align*}
Since $v, \mathcal{T}_j(v), v_2, u_j^2$ are orthonormal vectors, by (\ref{5.8}), (\ref{5.28}) and $L(v_2,u_j^2)=0$, we see that
\begin{align*}
0=&\lang L(v,v_2),L(u_j^2,\mathcal{T}_j(v))\rang+\lang L(v,u_j^2),L(v_2,\mathcal{T}_j(v))\rang+\lang L(v,\mathcal{T}_j(v)),L(u_j^2,v_2)\rang\\
 =&\lang L(v,v_2),L(u_j^2,\mathcal{T}_j(v))\rang+\lang L(v_2,\mathcal{T}_j(v)),L(v_2,\mathcal{T}_j(v))\rang.
\end{align*}
Thus $\lang L(v,v_2),L(u_j^2,\mathcal{T}_j(v))\rang=-\tau$. Applying the Cauchy-Schwarz inequality we deduce that
\be\label{5.30} L(\mathcal{T}_j(v),u_j^2)=-L(v,v_2).\ee
Combining (\ref{5.29}) and (\ref{5.30}), we get $L(v_2,\mathcal{T}_j^2(v)+v)=0$, which implies that $\mathcal{T}_j^2(v)+v\in V_{v_2}(0)$. As $\mathcal{T}_j^2(v)+v\in V_1\subset V_{v_2}(\tau)$, it follows that $\mathcal{T}_j^2(v)=-v$ for a unit vector $v$ and then by linearity for all $v\in V_1$, as claimed by (c).

To verify (d), we note that, for $j\neq l$, $L(v,u_j^2)\perp L(v,u_l^2)$. It follows from (\ref{5.28}) that
\begin{align*}
0&=\lang L(v,u_j^2), L(v,u_l^2)\rang=\lang L(v_2,\mathcal{T}_j(v)), L(v_2,\mathcal{T}_l(v))\rang \\
&=\lang P_{v_2}\mathcal{T}_j(v),\mathcal{T}_l(v)\rang=\tau\lang \mathcal{T}_j(v),\mathcal{T}_l(v)\rang.
\end{align*}
Hence, $\mathcal{T}_j(v)\perp \mathcal{T}_l(v)$.

We look at the unit hypersphere $S^\mathfrak{p}(1)\subset V_1$, the above properties (a)-(d) show that at $v\in S^\mathfrak{p}(1)$ one has
$$T_vS^\mathfrak{p}(1)=\operatorname{span}\{\mathcal{T}_1(v),\cdots,\mathcal{T}_\mathfrak{p}(v)\}.$$

Hence, by the properties (a)-(d), the $\mathfrak{p}$-dimensional sphere $S^\mathfrak{p}(1)$ is parallelizable. Then, according to Bott and Milnor \cite{BM} and Kervaire \cite{K}, the dimension $\mathfrak{p}$ can only be equal to 1,3 or 7.   \hfill $\Box$
\section{ Case $\left\{\mathfrak{C}_{m}\right\}_{2 \leq m \leq n-2}$ with $k_0=1$}

In this section, we consider Case  $\mathfrak{C}_{m}$ with $2 \leq m \leq$ $n-2$ with the condition that in the decomposition (\ref{5.18}), $k_0 = 1$. We will prove the following theorem.

\begin{thm}\label{thm6.1}
Let $M^{n}$ be a Calabi hypersurface in $\mathbb{R}^{n+1}$ which has parallel and non-vanishing cubic form. If $\mathfrak{C}_{m}$ with $2 \leq m \leq$ $n-2$ accurs and the integer $k_{0}$, as defined in section 5.3, satisfies $k_{0}=1$, then $M^{n}$ can be decomposed as the Calabi product of a locally strongly convex hyperbolic centroaffine hypersurfaces with parallel cubic form  and a Calabi hypersurface with parallel cubic form, or the Calabi product of a locally strongly convex hyperbolic centroaffine hypersurface with parallel cubic form and a point.
\end{thm}

To prove Theorem \ref{thm6.1}, we first note that if $k_{0}=1$ then by Proposition \ref{prop5.15} we have $\operatorname{dim}(\operatorname{Im} L)=1$. Moreover, we can prove the following result.

\begin{lem}\label{lemma6.2}
If $\operatorname{dim}(\operatorname{Im} L)=1$, then there is a unit vector $w_{1} \in \operatorname{Im} L \subset \mathcal{D}_{3}$ such that $L$ has the expression
\be\label{6.1}
L\left(v_{1}, v_{2}\right)=\frac{1}{2} \mu_{1} G\left(v_{1}, v_{2}\right) w_{1}, \quad \forall v_{1}, v_{2} \in \mathcal{D}_{2}.\ee
\end{lem}

\textbf{Proof.}
The fact $\operatorname{dim}(\operatorname{Im} L)=1$ implies that we have a unit vector $\bar{w} \in \operatorname{Im} L \subset$ $\mathcal{D}_{3}$ and a symmetric bilinear form $\alpha$ over $\mathcal{D}_{2}$ such that
\be\label{6.2}
L\left(v_{1}, v_{2}\right)=\alpha\left(v_{1}, v_{2}\right) \bar{w}, \quad \forall v_{1}, v_{2} \in \mathcal{D}_{2}.\ee

We define $Q: \mathcal{D}_{2} \rightarrow \mathcal{D}_{2}$ by $G\left(Q v_{1}, v_{2}\right):=\alpha\left(v_{1}, v_{2}\right)$. From Lemma \ref{lemma5.12} we have
\be\label{6.3}
L\left(v_{1}, v_{2}\right)=0, \quad \text { if } G\left(v_{1}, v_{2}\right)=0.\ee

Now we see that $G\left(Q v_{1}, v_{2}\right)=0$ if $G\left(v_{1}, v_{2}\right)=0 .$
Hence, $Qv=bv$ for any $v\in \mathcal{D}_{2}$.

By (\ref{5.2}), for any $v\in \mathcal{D}_{2}$, we get
$$\fr{1}{4}\mu_1^2\lang v,v\rang^2=\lang L(v,v),L(v,v)\rang=\alpha^2(v,v)=b^2\lang v,v\rang^2.$$
Thus $b=\pm\fr{1}{2}\mu_1=:\varepsilon(v)\fr{1}{2}\mu_1$. It follows that
\be\label{6.4}
L(v_1,v_2)=\alpha(v_1,v_2)\bar w=\varepsilon(v_1)\fr{1}{2}\mu_1\lang v_1,v_2\rang\bar w. \ee
On the other hand, we get
$$L(v_2,v_1)=\varepsilon(v_2)\fr{1}{2}\mu_1\lang v_2,v_1\rang\bar w.$$
The above two equations show that, for any $v_1,v_2\in \mathcal{D}_{2}$, $\varepsilon(v_1)=\varepsilon(v_2)$ holds, which implies that $\varepsilon (v)$ is independent of $v$. We finally get the assertion by putting $w_1=\varepsilon(v_1)\bar w$ in (\ref{6.4}).                    \hfill $\Box$

In sequel of this section, we will fix the unit vector $w_{1} \in \mathcal{D}_{3}$ as in Lemma \ref{lemma6.2}. Then, besides $A_{e_{1}} w_{1}=0$, the next three lemmas give all informations about the Fubini-Pick tensor $A$.

\begin{lem}\label{lemma6.3}
There exists an orthonormal basis $\left\{v_{1}, \ldots, v_{m-1}\right\}$ of $\mathcal{D}_{2}$ such that
\begin{align*}
A(e_1,v_i)&=\fr{1}{2}\mu_1v_i,\quad A(w_1,v_i)=\fr{1}{2}\mu_1v_i,\quad 1\leq i\leq m-1,\\
A(v_i,v_j)&=\fr{1}{2}\mu_1(e_1+w_1)\delta_{ij}, \quad 1\leq i,j\leq m-1.
\end{align*}
\end{lem}

\textbf{Proof.}
From Lemma \ref{lemma5.1}, we see that $A_{w_{1}}$ maps $\mathcal{D}_{2}$ to $\mathcal{D}_{2}$. Note that $A_{w_{1}}$ is self-adjoint, then there exists an orthonormal basis $\left\{v_{1}, \ldots, v_{m-1}\right\}$ of $\mathcal{D}_{2}$ such that $A_{w_{1}} v_{i}=b_{i} v_{i}$ with eigenvalues $b_{i}$. As $v_{i} \in \mathcal{D}_{2}$, we have $A_{e_1} v_{i}=\frac{1}{2} \mu_{1} v_{i}$. By Lemma \ref{lemma6.2} we get
$$b_i=\lang A_{w_1}v_i,v_i\rang=\lang A(v_i,v_i),w_1\rang=\lang L(v_i,v_i),w_1\rang=\fr{1}{2}\mu_1.$$
By the definition of $L$ and Lemma \ref{lemma6.2}, we have
$$A(v_i,v_j)=L(v_i,v_j)+\fr{1}{2}\mu_1\delta_{ij}e_1=\fr{1}{2}\mu_1(e_1+w_1)\delta_{ij}.$$
\hfill $\Box$

\begin{lem}\label{lemma6.4}    $A(w_1,w_1)=\mu_1w_1$.
\end{lem}

\textbf{Proof.}
By Lemma \ref{lemma6.2}, for any unit vector $a\in V_j$, we compute both sides of (\ref{5.19}). The assertion immediately follows.       \hfill $\Box$

Finally, in case $\mathcal{D}_{3} \neq \mathbb{R} w_{1}$ and let $\left\{w_{2}, \ldots, w_{n-m}\right\}$ be an orthonormal basis of $\mathcal{D}_{3} \backslash \mathbb{R} w_{1}$, by Lemmas \ref{lemma5.4} and \ref{lemma6.2}, we immediately have:

\begin{lem}\label{lemma6.5}    $A(w_1,w_i)=0, \quad 2\leq i\leq n-m$.
\end{lem}

Now, we are ready to complete the proof of Theorem \ref{thm6.1}.

\textbf{Proof of Theorem \ref{thm6.1}.}
Based on Lemmas \ref{lemma6.2}, \ref{lemma6.3}, \ref{lemma6.4}, and \ref{lemma6.5}, by putting
$$t=\fr{\sqrt2}{2}e_1+\fr{\sqrt2}{2}w_1, \quad v=-\fr{\sqrt2}{2}e_1+\fr{\sqrt2}{2}w_1,$$
we see that if $\mathcal{D}_{3}=\mathbb{R} w_{1}$, then $\left\{t, v, v_{1}, \ldots, v_{m-1}\right\}$ (or, resp. if $\mathcal{D}_{3} \neq \mathbb{R} w_{1}$, then $\left.\left\{t, v, v_{1}, \ldots, v_{m-1}, w_{2}, \ldots, w_{n-m}\right\}\right)$ forms an orthonormal basis of $T_{p} M^{n}$, with respect to which, the Fubini-Pick tensor takes the following form
\be\label{6.5}\left\{
\begin{array}{ll}
A(t,t)=\fr{\sqrt2}{2}\mu_1t;\; A(t,v)=\fr{\sqrt2}{2}\mu_1v;\\
A(t,v_i)=\fr{\sqrt2}{2}\mu_1v_i, \quad (1\leq i\leq m-1);   \\
\text { if } \mathcal{D}_{3} \neq \mathbb{R} w_{1}, \ A(t, w_i)=0, \quad (2\leq i\leq n-m).
\end{array}
\right.
\ee
By parallel translation along geodesics (with respect to $\hat{\nabla}$ ) through $p$, we can extend $\left\{t, v, v_{1}, \ldots, v_{m-1}\right\}$ (if $\left.\mathcal{D}_{3}=\mathbb{R} w_{1}\right)$, or, resp. $\left\{t, v, v_{1}, \ldots, v_{m-1}, w_{2}\right.$ ...,$\left.w_{n-m}\right\}$ (if $\mathcal{D}_{3} \neq \mathbb{R} w_{1}$ ) to obtain a local orthonormal basis $\left\{T, V, V_{1}, \ldots,\right.$, $\left.V_{m-1}\right\}$, or, resp. $\left\{T, V, V_{1}, \ldots, V_{m-1}, W_{2}, \ldots, W_{n-m}\right\}$ such that
$$
\left\{\begin{array}{l}
A(T, T)=\fr{\sqrt2}{2}\mu_1 T ; \quad A(T, V)=\fr{\sqrt2}{2}\mu_1 V ; \quad A\left(T, V_{i}\right)=\fr{\sqrt2}{2}\mu_1 V_{i}, \quad 1 \leq i \leq m-1; \\
\text { if } \mathcal{D}_{3} \neq \mathbb{R} w_{1}, \quad A\left(T, W_{i}\right)=0, \quad 2 \leq i \leq n-m.
\end{array}\right.
$$
Now, the above fact implies that, if $\mathcal{D}_{3} \neq \mathbb{R} w_{1}$ we can apply Theorem 3.6 to conclude that $M^{n}$ is decomposed as the Calabi product of a locally strongly convex hyperbolic centroaffine hypersurfaces with parallel cubic form and a Calabi hypersurface with parallel cubic form. If $\mathcal{D}_{3}= \mathbb{R} w_{1}$, then we can apply Theorem 3.3 to conclude that $M$ can be decomposed as the Calabi product of a locally strongly convex hyperbolic centroaffine hypersurfaces with parallel cubic form and a point. \hfill $\Box$

\section{ Case $\left\{\mathfrak{C}_{m}\right\}_{2 \leq m \leq n-2}$ with $k_0\geq2$ and $\mathfrak{p}=0$}

In this section, we will prove the following theorem.

\begin{thm}\label{thm7.1} Let $M^{n}$ be a Calabi hypersurface in $\mathbb{R}^{n+1}$ which has parallel and non-vanishing cubic form. If $\mathfrak{C}_{m}$ with $2 \leq m \leq n-2$ accurs and the integer $k_{0}$ and $\mathfrak{p}$, satisfy  $k_{0} \geq 2$ and $\mathfrak{p}=0$, then $n \geq \frac{1}{2} m(m+1)$. Moreover, we have either
\begin{enumerate}
\item [(i)] $n=\frac{1}{2} m(m+1)$, $M^{n}$ can be decomposed as the Calabi product of a locally strongly convex hyperbolic centroaffine hypersurface with parallel cubic form and a point, or
\item [(ii)] $n>\frac{1}{2} m(m+1)$, $M^{n}$ can be decomposed as the Calabi product of a locally strongly convex hyperbolic centroaffine hypersurfaces with parallel cubic form  and a Calabi hypersurface with parallel cubic form.
      \end{enumerate}
\end{thm}

In the present situation, the decomposition (5.18) reduces to $\mathcal{D}_{2}=\left\{v_{1}\right\} \oplus$ $\cdots \oplus\left\{v_{k_{0}}\right\}$, and $\left\{v_{1}, \ldots, v_{k_{0}}\right\}$ forms an orthonormal basis of $\mathcal{D}_{2}$. Then $\operatorname{dim} \mathcal{D}_{2}=k_{0}=m-1, m \geq 3$, which means that $\dim(\operatorname{Im} L)\leq k_0+\cdots+1=\fr{m(m-1)}{2}$.

According to (\ref{5.3}), Lemmas \ref{lemma5.3}, \ref{lemma5.10} and the fact that for $j \neq l, v_{j} \in V_{v_{l}}(\tau)$, we have
\begin{align}
\label{7.1}\lang L(v_j,v_l),L(v_j,v_l)\rang&=\tau, \quad j\neq l, \\
\label{7.2}\lang L(v_j,v_{l_1}),L(v_j,v_{l_2})\rang&=0, \quad j, l_1, l_2\, \text{distinct}, \\
\label{7.3}\lang L(v_{j_1},v_{j_2}),L(v_{j_3},v_{j_4})\rang&=0, \quad j_1,  j_2, j_3, j_4\,\text{distinct}, \\
\label{7.4}\lang L(v_j,v_j),L(v_j,v_j)\rang&=\fr{1}{4}\mu_1^2,  \\
\label{7.5}\lang L(v_j,v_j),L(v_l,v_l)\rang&=0, \quad j\neq l, \\
\label{7.6}\lang L(v_j,v_j),L(v_j,v_{l})\rang&=0, \quad j\neq l, \\
\label{7.7}\lang L(v_j,v_j),L(v_{l_1},v_{l_2})\rang&=0, \quad j, l_1, l_2\, \text{distinct}.
\end{align}
Denote $L_{j}:=L\left(v_{1}, v_{1}\right)+\cdots+L\left(v_{j}, v_{j}\right)-j L\left(v_{j+1}, v_{j+1}\right), 1 \leq j \leq k_{0}-1$. Then it is easy to check
$$G\left(L_{j}, L_{j}\right)=\fr{1}{4}\mu_1^2\times j+j^2\times\fr{1}{4}\mu_1^2=2j(j+1)\tau\neq0,$$
and that
\be\label{7.8}
\begin{cases}w_{j}=\frac{1}{\sqrt{2 j(j+1) \tau}} L_{j}, & 1 \leq j \leq k_{0}-1 \\ w_{k l}=\frac{1}{\sqrt{\tau}} L\left(v_{k}, v_{l}\right), & 1 \leq k<l \leq k_{0}\end{cases}
\ee
give $\frac{1}{2}(m+1)(m-2)$ orthonormal vectors in $\operatorname{Im} L \subset \mathcal{D}_{3}$. Direct computations show that $\operatorname{Tr} L=L\left(v_{1}, v_{1}\right)+\cdots+L\left(v_{k_{0}}, v_{k_{0}}\right)\in \operatorname{Im} L\subset \mathcal{D}_{3}$ is orthogonal to all vectors in (\ref{7.8}), and by using (\ref{7.4}), (\ref{7.5}) and the fact that $v_{i} \in V_{v_{j}}(\tau), i \neq j$, we get
$$
G(\operatorname{Tr} L, \operatorname{Tr} L) =\fr{1}{4}k_0\mu_1^2=: \rho^{2}>0
$$
Hence, we get $\fr{1}{2}(m+1)(m-2)+1=\fr{m(m-1)}{2}$ orthogonal vectors in $\operatorname{Im} L\subset \mathcal{D}_{3}$. It follows that $\dim(\operatorname{Im} L)=\fr{m(m-1)}{2}$.

Thus, we have the estimate of the dimension
$$
\begin{aligned}
n  =1+\operatorname{dim}\left(\mathcal{D}_{2}\right)+
\operatorname{dim}\left(\mathcal{D}_{3}\right)
  \geq 1+m-1+\frac{1}{2}m(m-1)=\frac{1}{2} m(m+1).
\end{aligned}
$$

Now, we are ready to complete the proof of Theorem \ref{thm7.1}.\\
\textbf{Proof of Theorem \ref{thm7.1}.} We need to consider two cases:

Case (1).\ $n=\fr{m(m+1)}{2}$,\ (if $\operatorname{Im} L= \mathcal{D}_{3}$).
Case (2).\ $n>\fr{m(m+1)}{2}$,\ ( if $\operatorname{Im} L \varsubsetneqq \mathcal{D}_{3}$). \\
For both cases (1) and (2), as $\operatorname{Tr} L\neq0$, we can define a unit vector $t:=\fr{1}{\rho}\operatorname{Tr} L$.

In Case (1), from the previous discussions, we see that
$$\{t, w_j|_{1\leq j\leq k_0-1}, w_{kl}|_{1\leq k<l\leq k_0}\}$$
forms an orthonormal basis of $\operatorname{Im} L=\mathcal{D}_{3}$.

In Case (2), we choose $\{\td w_i|_{1\leq i\leq \td n}\}$ as an orthonormal basis of $\mathcal{D}_{3}\backslash \operatorname{Im} L$ such that
$$\{t, w_j|_{1\leq j\leq k_0-1}, w_{kl}|_{1\leq k<l\leq k_0}, \td w_i|_{1\leq i\leq \td n}\}$$
is an orthonormal basis of $\mathcal{D}_{3}$,  where $\td n=n-\fr{1}{2}m(m+1)$.

\begin{lem}\label{lemma7.2}
the Fubini-Pick tensor satisfies
\be\label{7.9}
\left\{\begin{array}{lll}
A(e_1,t)=0, \quad A(t,v_i)=\fr{\rho}{k_0}v_i, \ 1\leq i\leq k_0, \\
A(t,t)=\fr{2\rho}{k_0}t, \quad A(t,w_j)=\fr{2\rho}{k_0}w_j, \quad 1\leq j\leq k_0-1,\\
A(t,w_{kl})=\fr{2\rho}{k_0}w_{kl}, \quad 1\leq k<l\leq k_0, \\
\text { if } \operatorname{Im} L \varsubsetneqq \mathcal{D}_{3}, \quad A(t,\td w_i)=0, \quad 1\leq i\leq \td n.
\end{array}
\right.
\ee
\end{lem}

\textbf{Proof.}
By (iii) of Lemma \ref{lemma5.1}, we know $A(t,v_i)\in \mathcal{D}_2$. Hence
$$\lang A(t,v_i),v_j\rang=\fr{1}{\rho}\sum_{l=1}^{k_0}\lang A(v_i,v_j),L(v_l,v_l)\rang=\fr{\rho}{k_0}\delta_{ij},\quad 1\leq i,j\leq k_0.$$
From Lemma \ref{lemma5.14} we get
\begin{align*}
A(t,t)&=\fr{1}{\rho^2} \sum_{i=1}^{k_0}A(L(v_i,v_i),L(v_i,v_i))=\fr{2\rho}{k_0}t, \\
A(t,w_j)&=\fr{\mu_1^2}{2\rho\sqrt{2j(j+1)\tau}}L_j=\fr{2\rho}{k_0}w_j, \quad 1\leq j\leq k_0-1,\\
A(t,w_{kl})&=\fr{1}{\rho\sqrt \tau}(A(L(v_k,v_k),L(v_k,v_l))+A(L(v_l,v_l),L(v_l,v_k)))\\ &=\fr{2\rho}{k_0}w_{kl}, \quad 1\leq k<l\leq k_0.
\end{align*}
If $\operatorname{Im} L \varsubsetneqq \mathcal{D}_{3}$, by Lemma \ref{lemma5.4}, one get
$$A(t,\td w_i)=\fr{1}{\rho}\sum_{j=1}^{k_0}A(L(v_j,v_j),\td w_i)=0, \quad 1\leq i\leq \td n.$$

This completes the proof of Lemma \ref{lemma7.2}.    \hfill $\Box$

Put
\be\label{7.10}
T=\fr{1}{\sqrt{k_0+1}}e_1+\sqrt{\fr{k_0}{k_0+1}}t,\quad T^*=-\sqrt{\fr{k_0}{k_0+1}}e_1+\fr{1}{\sqrt{k_0+1}}t.
\ee
It is easily to see that if $\operatorname{Im} L=\mathcal{D}_{3}$, then $\{T, T^*, v_j|_{1\leq j\leq k_0}, w_j|_{1\leq j\leq k_0-1}, w_{kl}|_{1\leq k<l\leq k_0}\}$ (or, resp. if $\operatorname{Im} L \varsubsetneqq \mathcal{D}_{3}$, then $\{T, T^*, v_j|_{1\leq j\leq k_0}, w_j|_{1\leq j\leq k_0-1}, w_{kl}|_{1\leq k<l\leq k_0}, \td w_i|_{1\leq i\leq \td n}\}$) is an orthonormal basis of $T_pM$. By Lemmas \ref{lemma5.1} and \ref{lemma7.2} we have the following
\begin{lem}\label{lemma7.3}
With respect to the above notations, the Fubini-Pick tensor takes the following form
\be\label{7.11}\left\{
\begin{array}{ll}
A(T,T)=\sigma T, \quad A(T,T^*)=\sigma T^*,\quad A(T,v_i)=\sigma v_i, \  1\leq i\leq k_0; \\
A(T,w_j)=\sigma w_j, \  1\leq j\leq k_0-1; \quad
A(T,w_{kl})=\sigma w_{kl}, \  1\leq k<l\leq k_0; \\
\text { if } \operatorname{Im} L \varsubsetneqq \mathcal{D}_{3}, \quad A(T,\td w_i)=0, \quad 1\leq i\leq \td n.
\end{array}
\right.
\ee
where $\sigma=\fr{1}{\sqrt{k_0+1}}\mu_1$.
\end{lem}

Given the parallelism of the Fubini-Pick form $A$,  Lemma \ref{lemma7.3} and Theorem \ref{thm3.3}, we conclude that in case (1), $M$ can be decomposed as the Calabi product of a locally strongly convex hyperbolic centroaffine hypersurfaces with parallel cubic form and a point.  Similarly, using Theorem \ref{thm3.6}, we conclude that in case (2), $M$ can be decomposed as the Calabi product of a locally strongly convex hyperbolic centroaffine hypersurfaces with parallel cubic form and a Calabi hypersurface with parallel cubic form.

The combination of the preceding two cases' discussion then completes the proof of Theorem \ref{thm7.1}. \hfill $\Box$
\section{ Case $\left\{\mathfrak{C}_{m}\right\}_{2 \leq m \leq n-2}$ with $k_0\geq2$ and $\mathfrak{p}=1$}

In this section, we will prove the following theorem.
\begin{thm}\label{thm8.1}  Let $M^{n}$ be a Calabi hypersurface in $\mathbb{R}^{n+1}$ which has parallel and non-vanishing cubic form. If $\mathfrak{C}_{m}$ with $2 \leq m \leq n-2$ accurs and the integer $k_{0}$ and $\mathfrak{p}$, satisfy $k_{0} \geq 2$ and $\mathfrak{p}=1$, then $n \geq \frac{1}{4}(m+1)^{2}.$ Moreover, we have either
\begin{enumerate}
\item [(i)] $n= \frac{1}{4}(m+1)^2$, $ M^{n}$ can be decomposed as the Calabi product of a locally strongly convex hyperbolic centroaffine hypersurface with parallel cubic form and a point, or
\item [(ii)] $n>\frac{1}{4}(m+1)^2$, $ M^{n}$ can be decomposed as the Calabi product of a locally strongly convex hyperbolic centroaffine hypersurfaces with parallel cubic form  and a Calabi hypersurface with parallel cubic form.
      \end{enumerate}
\end{thm}

Now we have $\operatorname{dim} \mathcal{D}_{2}=m-1=2 k_{0}$ and $m \geq 5$. Firstly, we will prove the following

\begin{lem}\label{lemma8.2} In the decomposition (\ref{5.18}), if we have $k_{0} \geq 2$ and $\mathfrak{p}=1$, then there exist unit vectors $u_{j} \in V_{v_{j}}(0)\left(1 \leq j \leq k_{0}\right)$ such that the orthonormal basis $\left\{v_{1}, u_{1}, \ldots, v_{k_{0}}, u_{k_{0}}\right\}$ of $\mathcal{D}_{2}$ satisfies the relations
\be\label{8.1}
L\left(u_{l}, v_{j}\right)=-L\left(v_{l}, u_{j}\right), \quad L\left(v_{l}, v_{j}\right)=L\left(u_{l}, u_{j}\right), \quad 1 \leq j, l \leq k_{0}.
\ee
\end{lem}

\textbf{Proof.}
As for each $1 \leq j \leq k_{0}$ it holds $\operatorname{dim}\left(V_{v_{j}}(0)\right)=1$, we assume $V_{v_{2}}(0)=$ $\left\{u_{2}\right\}$ for a unit vector $u_{2}$. Then, for each $j \neq 2$, by Lemma \ref{lemma5.13}, we have a unique unit vector $u_{j} \in V_{v_{j}}(0)$ satisfying
\be\label{8.2}
L\left(u_{2}, v_{j}\right)=-L\left(v_{2}, u_{j}\right), \quad L\left(v_{2}, v_{j}\right)=L\left(u_{2}, u_{j}\right), \quad 1 \leq j \leq k_{0}, j \neq 2.
\ee

Moreover, Lemma \ref{lemma5.8} implies that (\ref{8.2}) also holds for $j=2$. Next, we state\\
\textbf{Claim 1.} $L(u_l,v_j)=-L(u_j,v_l)$, \ $L(v_l,v_j)=L(u_l,u_j)$, \ $1\leq j,l\leq k_0, \ j,l \neq 2$.

To verify the claim, as $u_{j} \in V_{v_{j}}(0)$, we first see by Lemma \ref{lemma5.8} that $L\left(u_{j}, v_{j}\right)=0$ and $L\left(v_{j}, v_{j}\right)=L\left(u_{j}, u_{j}\right) .$ Hence the claim is true for $j=l$.

Now we fix $j \neq l$ and $j, l \neq 2$. By the fact $u_l\in V_{v_l}(0)$ and Lemma \ref{lemma5.13}, there exists a unique unit vector $u_{j}^{(l)} \in V_{v_{j}}(0)$, such that
\be\label{8.3}
L\left(u_{l}, v_{j}\right)=-L\left(v_l, u_{j}^{(l)}\right), \quad L\left(v_l, v_{j}\right)=L\left(u_{l}, u_{j}^{(l)}\right).
\ee

In the following, we prove $u_j=u_j^{(l)}$ in (\ref{8.2}) and (\ref{8.3}).
Noting that $\operatorname{dim}\left(V_{v_{j}}(0)\right)=1$ and $u_{j}^{(l)}, u_{j} \in V_{v_{j}}(0)$ are unit vectors, we have $u_{j}^{(l)}=\epsilon u_{j}$ with $\epsilon=\pm 1$. Hence from (\ref{8.3}) we have
\be\label{8.4}
L\left(u_{l}, v_{j}\right)=-\epsilon L\left(v_{l}, u_{j}\right), \quad L\left(v_{l}, v_{j}\right)=\epsilon L\left(u_{l}, u_{j}\right).
\ee
On the other hand, by using (\ref{5.23}), (\ref{8.2}) and (\ref{8.4}), we get
$$
\begin{aligned}
&A\left(L\left(v_{j}, v_{l}\right), L\left(v_{2}, u_{j}\right)\right)=A\left(L\left(v_{j}, v_{l}\right),-L\left(v_{j}, u_{2}\right)\right)=-\tau L\left(v_{l}, u_{2}\right), \\
&A\left(L\left(v_{j}, v_{l}\right), L\left(v_{2}, u_{j}\right)\right)=A\left(\epsilon L\left(u_{j}, u_{l}\right), L\left(v_{2}, u_{j}\right)\right)=-\epsilon \tau L\left(v_{l}, u_{2}\right).
\end{aligned}
$$
From the comparison of the above two equations we get $\epsilon=1$.
From (\ref{8.4}) we have verified Claim 1 and the proof of Lemma \ref{lemma8.2} is fullled.
\hfill $\Box$

By the symmetry of $L$ and Lemma \ref{lemma8.2}, we get
\be\label{8.5}
\operatorname{dim}(\operatorname{Im} L)\leq k_0+\cdots+1+(k_0-1)+\cdots+1=\fr{1}{4}(m-1)^2.
\ee

To continue the proof of Theorem \ref{thm8.1}, we now assume that $k_{0} \geq 2$ and let $\left\{v_{1}, u_{1}, \ldots, v_{k_{0}}, u_{k_{0}}\right\}$ be the orthonormal basis of $\mathcal{D}_{2}$ as constructed in Lemma \ref{lemma8.2}.

Given Lemmas \ref{lemma5.2}, \ref{lemma5.3}, \ref{lemma5.8}, \ref{lemma5.10} and that for $j \neq l, v_{j}, u_{j} \in V_{v_{l}}(\tau)=V_{u_{l}}(\tau)$, we have the following calculations:
\begin{align}
\label{8.6}
&G\left(L\left(v_{j}, u_{l}\right), L\left(v_{j}, u_{l}\right)\right)=G\left(L\left(v_{j}, v_{l}\right), L\left(v_{j}, v_{l}\right)\right)=\tau, \quad j \neq l,\\
\label{8.7}
&G\left(L\left(u_{j}, v_{l_{1}}\right), L\left(u_{j}, v_{l_{2}}\right)\right)=G\left(L\left(v_{j}, u_{l_{1}}\right), L\left(v_{j}, u_{l_{2}}\right)\right)       \notag \\
&\quad=G\left(L\left(v_{j}, v_{l_{1}}\right), L\left(v_{j}, v_{l_{2}}\right)\right)=0, \quad j, l_{1}, l_{2} \text { distinct}, \\
\label{8.8}
&G\left(L\left(v_{j_{1}}, v_{j_{2}}\right), L\left(v_{j_{3}}, v_{j_{4}}\right)\right)=0, \quad j_{1}, j_{2}, j_{3}, j_{4} \text { distinct}, \\
\label{8.9}
&G\left(L\left(v_{j}, v_{l}\right), L\left(v_{j_{1}}, u_{l_{1}}\right)\right)=0, \quad j \neq l \text { and } j_{1} \neq l_{1}, \\
\label{8.10}
&G\left(L\left(v_{j}, v_{j}\right), L\left(v_{j}, v_{j}\right)\right)=\frac{1}{4} \mu_{1}^2, \quad 1 \leq j \leq k_{0}, \\
\label{8.11}
&G\left(L\left(v_{j}, v_{j}\right), L\left(v_{l}, v_{l}\right)\right)=0, \quad 1 \leq j \neq l \leq k_{0}, \\
\label{8.12}
&G\left(L\left(v_{j}, v_{j}\right), L\left(v_{j}, v_{l}\right)\right)=G\left(L\left(v_{j}, v_{j}\right), L\left(v_{j}, u_{l}\right)\right)    \notag  \\
&\quad=G\left(L\left(v_{j}, v_{j}\right), L\left(v_{l}, u_{j}\right)\right)=0, \quad 1 \leq j \neq l \leq k_{0}, \\
\label{8.13}
&G\left(L\left(v_{j}, v_{j}\right), L\left(v_{l_{1}}, v_{l_{2}}\right)\right)=G\left(L\left(v_{j}, v_{j}\right), L\left(v_{l_{1}}, u_{l_{2}}\right)\right)=0,   \notag   \\
&\quad 1 \leq j, l_{1}, l_{2} \text { distinct } \leq k_{0}.
\end{align}

Similar as in the proof of Theorem \ref{thm7.1}, we denote
$$
L_{j}:=L\left(v_{1}, v_{1}\right)+\cdots+L\left(v_{j}, v_{j}\right)-j L\left(v_{j+1}, v_{j+1}\right), \quad 1 \leq j \leq k_{0}-1.
$$
(\ref{8.10}) and (\ref{8.11}) show that $G\left(L_{j}, L_{j}\right)=2 j(j+1) \tau \neq 0$ for each $j$. Moreover (\ref{8.6})-(\ref{8.13}) and Lemma \ref{lemma5.10} show that
\be\label{8.14}
\left\{\begin{array}{cl}
w_{j}=\frac{1}{\sqrt{2 j(j+1) \tau}} L_{j}, & 1 \leq j \leq k_{0}-1, \\
w_{k l}=\frac{1}{\sqrt{\tau}} L\left(v_{k}, v_{l}\right), & 1 \leq k<l \leq k_{0}, \\
w_{k l}^{\prime}=\frac{1}{\sqrt{\tau}} L\left(v_{k}, u_{l}\right), & 1 \leq k<l \leq k_{0}
\end{array}\right.
\ee
give $\frac{1}{4}(m+1)(m-3)$ mutually orthogonal unit vectors in $\operatorname{Im} L \subset \mathcal{D}_{3}$.

Direct computations and Lemma \ref{lemma8.2} show that $\operatorname{Tr} L=L(v_1,v_1)+L(u_1,u_1)+\cdots+L(v_{k_0},v_{k_0})+L(u_{k_0},u_{k_0})=2\left[L\left(v_{1}, v_{1}\right)+\cdots+L\left(v_{k_{0}}, v_{k_{0}}\right)\right]\in \operatorname{Im} L$ is orthogonal to all vectors in (\ref{8.14}), and by using (\ref{8.10}), (\ref{8.11}) and the fact that $v_{i} \in V_{v_{j}}(\tau)$ for $i \neq j$, we get
\be\label{8.15}
\frac{1}{4} G(\operatorname{Tr} L, \operatorname{Tr} L)=\frac{1}{4} k_{0}\mu_1^2=: \rho^{2}>0.\ee
(\ref{8.14}) and (\ref{8.15}) give $\frac{1}{4}(m+1)(m-3)+1=\fr{1}{4}(m-1)^2$ mutually orthogonal vectors in $\operatorname{Im} L \subset \mathcal{D}_{3}$. Combining (\ref{8.5}) one can get $\operatorname{dim}(\operatorname{Im} L)=\fr{1}{4}(m-1)^2$.

Thus, we have the estimate of the dimension
$$
\begin{aligned}
n  =1+\operatorname{dim}\left(\mathcal{D}_{2}\right)+
\operatorname{dim}\left(\mathcal{D}_{3}\right)
  \geq 1+m-1+\fr{1}{4}(m-1)^2=\frac{1}{4}(m+1)^{2}.
\end{aligned}
$$

Now, we are ready to complete the proof of Theorem \ref{thm8.1}.\\
\textbf{Proof of Theorem \ref{thm8.1}.} We need to consider two cases:

Case (1).\ $n=\fr{(m+1)^2}{4}$,\ ( if $\operatorname{Im} L= \mathcal{D}_{3}$).

Case (2).\ $n>\fr{(m+1)^2}{4}$,\ (if $\operatorname{Im} L \varsubsetneqq \mathcal{D}_{3}$).\\
For both cases (1) and (2), as $\operatorname{Tr} L\neq0$, we can define a unit vector $t:=\fr{1}{2\rho}\operatorname{Tr} L$.

In Case (1), from the previous discussions, we see that
$$\{t, w_j|_{1\leq j\leq k_0-1}, w_{kl}|_{1\leq k<l\leq k_0}, w_{k l}^{\prime}|_{1\leq k<l\leq k_0}\} $$
forms an orthonormal basis of $\operatorname{Im} L=\mathcal{D}_{3}$.

In Case (2), we choose $\{\td w_i|_{1\leq i\leq \td n}\}$
as an orthonormal basis of $\mathcal{D}_{3}\backslash \operatorname{Im} L$ such that
$$\{t, w_j|_{1\leq j\leq k_0-1}, w_{kl}|_{1\leq k<l\leq k_0}, w_{k l}^{\prime}|_{1\leq k<l\leq k_0}, \td w_i|_{1\leq i\leq \td n}\}$$
forms an orthonormal basis of $\mathcal{D}_{3}$,  where $\td n=n-\fr{(m+1)^2}{4}$.

\begin{lem}\label{lemma8.3}
the Fubini-Pick tensor satisfies
\be\label{8.16}
\left\{\begin{array}{lll}
A(e_1,t)=0, \quad A(t,v_i)=\fr{\rho}{k_0}v_i, \quad A(t,u_i)=\fr{\rho}{k_0}u_i, \ 1\leq i\leq k_0, \\
A(t,t)=\fr{2\rho}{k_0}t, \quad A(t,w_j)=\fr{2\rho}{k_0}w_j, \quad 1\leq j\leq k_0-1,\\
A(t,w_{kl})=\fr{2\rho}{k_0}w_{kl}, \quad A(t,w_{kl}')=\fr{2\rho}{k_0}w_{kl}' ,\quad 1\leq k<l\leq k_0, \\
\text { if } \operatorname{Im} L \varsubsetneqq \mathcal{D}_{3}, \quad A(t,\td w_i)=0, \quad 1\leq i\leq \td n.
\end{array}
\right.
\ee
\end{lem}

\textbf{Proof.}
By (iii) of Lemma \ref{lemma5.1}, we know $A(t,v_i)\in \mathcal{D}_2$. (\ref{8.10}), (\ref{8.11}) and (\ref{8.13}) show that
\begin{align*}
\lang A(t,v_i),v_j\rang&=\fr{1}{\rho}\sum_{l=1}^{k_0}\lang A(v_i,v_j),L(v_l,v_l)\rang=\fr{\rho}{k_0}\delta_{ij},\quad 1\leq i,j\leq k_0,\\
\lang A(t,v_i),u_j\rang&=\fr{1}{\rho}\sum_{l=1}^{k_0}\lang A(v_i,u_j),L(v_l,v_l)\rang=0,\quad 1\leq i,j\leq k_0.
\end{align*}
Hence, we have $A(t,v_i)=\fr{\rho}{k_0}v_i, \ 1\leq i\leq k_0$. Similarly, $A(t,u_i)=\fr{\rho}{k_0}u_i, \ 1\leq i\leq k_0$ can be proved.

From Lemmas \ref{lemma5.14} and \ref{lemma8.2} we get
\begin{align*}
A(t,t)&=\fr{1}{\rho^2} \sum_{i=1}^{k_0}A(L(v_i,v_i),L(v_i,v_i))=\fr{2\rho}{k_0}t, \\
A(t,w_j)&=\fr{\mu_1^2}{2\rho\sqrt{2j(j+1)\tau}}L_j=\fr{2\rho}{k_0}w_j, \quad 1\leq j\leq k_0-1,\\
A(t,w_{kl})&=\fr{1}{\rho\sqrt \tau}(A(L(v_k,v_k),L(v_k,v_l))+A(L(v_l,v_l),L(v_l,v_k)))=\fr{2\rho}{k_0}w_{kl}, \quad 1\leq k<l\leq k_0,\\
A(t,w_{kl}')&=\fr{1}{\rho\sqrt \tau}(A(L(v_k,v_k),L(v_k,u_l))+A(L(v_l,v_l),L(u_l,v_k)))   \\
&=\fr{1}{\rho\sqrt \tau}(A(L(v_k,v_k),L(v_k,u_l))+A(L(u_l,u_l),L(u_l,v_k)))=\fr{2\rho}{k_0}w_{kl}', \quad 1\leq k<l\leq k_0.
\end{align*}

If $\operatorname{Im} L \varsubsetneqq \mathcal{D}_{3}$, by Lemma \ref{lemma5.4}, one get
$$A(t,\td w_i)=\fr{1}{\rho}\sum_{j=1}^{k_0}A(L(v_j,v_j),\td w_i)=0, \quad 1\leq i\leq \td n.$$

This completes the proof of Lemma \ref{lemma8.3}.    \hfill $\Box$

Put
\be\label{8.17}
T=\fr{1}{\sqrt{k_0+1}}e_1+\sqrt{\fr{k_0}{k_0+1}}t,\quad T^*=-\sqrt{\fr{k_0}{k_0+1}}e_1+\fr{1}{\sqrt{k_0+1}}t.
\ee
It is easily to see that if $\operatorname{Im} L=\mathcal{D}_{3}$, then $\{T, T^*, v_j|_{1\leq j\leq k_0} , u_j|_{1\leq j\leq k_0}, w_j|_{1\leq j\leq k_0-1},\\ w_{kl}|_{1\leq k<l\leq k_0}, w_{kl}'|_{1\leq k<l\leq k_0}\}$ (or, resp. if $\operatorname{Im} L \varsubsetneqq \mathcal{D}_{3}$, then $\{T, T^*, v_j|_{1\leq j\leq k_0} , u_j|_{1\leq j\leq k_0}, \\ w_j|_{1\leq j\leq k_0-1}, w_{kl}|_{1\leq k<l\leq k_0}, w_{kl}'|_{1\leq k<l\leq k_0}, \td w_i|_{1\leq i\leq \td n}\}$) is an orthonormal basis of $T_pM$. By Lemmas \ref{lemma5.1} and \ref{lemma8.3} we have the following
\begin{lem}\label{lemma8.4}
With respect to the above notations, the Fubini-Pick tensor takes the following form
$$\label{8.18}\left\{
\begin{array}{ll}
A(T,T)=\sigma T, \quad A(T,T^*)=\sigma T^*,\quad A(T,v_i)=\sigma v_i, \quad A(T,u_i)=\sigma u_i, \  1\leq i\leq k_0, \\
A(T,w_j)=\sigma w_j, \  1\leq j\leq k_0-1, \quad
A(T,w_{kl})=\sigma w_{kl}, \  1\leq k<l\leq k_0, \\
A(T,w_{kl}')=\sigma w_{kl}', \ 1\leq k<l\leq k_0,\\   \text { if } \operatorname{Im} L \varsubsetneqq \mathcal{D}_{3}, \quad A(T,\td w_i)=0, \quad 1\leq i\leq \td n,
\end{array}
\right.
$$
where $\sigma=\fr{1}{\sqrt{k_0+1}}\mu_1$.
\end{lem}

Now, the above fact implies that, if $n> \frac{1}{4}(m+1)^2$ we can apply Theorem 3.6 to conclude that $M^{n}$ is decomposed as the Calabi product of a locally strongly convex hyperbolic centroaffine hypersurfaces with parallel cubic form and a Calabi hypersurface with parallel cubic form. If $n= \frac{1}{4}(m+1)^2$, then we can apply Theorem 3.3 to conclude that $M$ can be decomposed as the Calabi product of a locally strongly convex hyperbolic centroaffine hypersurfaces with parallel cubic form and a point. It completes the proof of Theorem \ref{thm8.1}. \hfill $\Box$
\section{ Case $\left\{\mathfrak{C}_{m}\right\}_{2 \leq m \leq n-2}$ with $k_0\geq2$ and $\mathfrak{p}=3$}

In this section, we will prove the following theorem.

\begin{thm}\label{thm9.1}
Let $M^{n}$ be a Calabi hypersurface in $\mathbb{R}^{n+1}$ which has parallel and non-vanishing cubic form. If $\mathfrak{C}_{m}$ with $2 \leq m \leq n-2$ accurs and the integer $k_{0}$ and $\mathfrak{p}$, satisfy $k_{0} \geq 2$ and $\mathfrak{p}=3$, then $n \geq \frac{1}{8}(m+1)(m+3).$ Moreover, we have either
\begin{enumerate}
\item [(i)] $n= \frac{1}{8}(m+1)(m+3)$, $ M^{n}$ can be decomposed as the Calabi product of a locally strongly convex hyperbolic centroaffine hypersurface with parallel cubic form and a point, or
\item [(ii)] $n>\frac{1}{8}(m+1)(m+3)$, $ M^{n}$ can be decomposed as the Calabi product of a locally strongly convex hyperbolic centroaffine hypersurfaces with parallel cubic form  and a Calabi hypersurface with parallel cubic form.
      \end{enumerate}
\end{thm}

Now we have $\operatorname{dim} \mathcal{D}_{2}=m-1=4 k_{0}$ and $m \geq 9 .$ Similarly, we will prove the following lemma.

\begin{lem}\label{lemma9.2} In the decomposition (\ref{5.18}), if we have $k_{0} \geq 2$ and $\mathfrak{p}=3$, then there exist unit vectors $x_{j}, y_{j}, z_{j} \in V_{v_{j}}(0)$  $\left(1 \leq j \leq k_{0}\right)$ such that the orthonormal basis $\left\{v_{1}, x_{1}, y_{1}, z_{1} ; \ldots ; v_{k_{0}}, x_{k_{0}}, y_{k_{0}}, z_{k_{0}}\right\}$ of $\mathcal{D}_{2}$ satisfies the relations,
\be\label{9.1}
\left\{\begin{array}{l}
L\left(v_{j}, v_{l}\right)=L\left(x_{j}, x_{l}\right)=L\left(y_{j}, y_{l}\right)=L\left(z_{j}, z_{l}\right), \\
L\left(v_{j}, x_{l}\right)=-L\left(x_{j}, v_{l}\right)=-L\left(y_{j}, z_{l}\right)=L\left(z_{j}, y_{l}\right), \\
L\left(v_{j}, y_{l}\right)=-L\left(y_{j}, v_{l}\right)=-L\left(z_{j}, x_{l}\right)=L\left(x_{j}, z_{l}\right), \\
L\left(v_{j}, z_{l}\right)=-L\left(z_{j}, v_{l}\right)=-L\left(x_{j}, y_{l}\right)=L\left(y_{j}, x_{l}\right),
\end{array}\right.
\ee where  $1\leq j,l\leq k_0$.
\end{lem}

\textbf{Proof.} As doing before, we denote $V_{j}=\left\{v_{j}\right\} \oplus V_{v_{j}}(0), 1 \leq j \leq k_{0} .$ Let us fix two orthogonal unit vectors $x_{1}, y_{1} \in V_{v_{1}}(0) .$ By using Lemma \ref{lemma5.13}, for each $j \neq 1$, we have two unit vectors $x_{j}, y_{j} \in V_{v_{j}}(0)$ such that
\be\label{9.2}
\left\{\begin{array}{l}
L\left(v_{j}, v_{1}\right)=L\left(x_{j}, x_{1}\right)=L\left(y_{j}, y_{1}\right), \\
L\left(v_{j}, x_{1}\right)=-L\left(x_{j}, v_{1}\right), \quad L\left(v_{j}, y_{1}\right)=-L\left(y_{j}, v_{1}\right).
\end{array}\right.
\ee
As $v_j\in V _{y_j}(0) \quad (j\neq 1)$ and $v_1\in V _{y_1}(0)$, according to Lemma \ref{lemma5.13}, we further have unit vectors $z_{1}^{j} \in V_{x_{1}}(0)$ and $z_{j} \in V_{x_{j}}(0)$ such that
\be\label{9.3}
\left\{\begin{array}{l}
L\left(v_{j}, z_{1}^{j}\right)=L\left(y_{j}, x_{1}\right), \quad L\left(v_{j}, x_{1}\right)=-L\left(y_{j}, z_{1}^{j}\right), \\
L\left(z_{j}, v_{1}\right)=L\left(x_{j}, y_{1}\right), \quad L\left(z_{j}, y_{1}\right)=-L\left(x_{j}, v_{1}\right).
\end{array}\right.
\ee
The important is that we have the following

\textbf{Claim 1.} For each $j \neq 1,\left\{x_{1}, y_{1}, z_{1}^{j}\right\}$ is an orthonormal basis of $V_{v_{1}}(0)$ and $\left\{x_{j}, y_{j}, z_{j}\right\}$ is an orthonormal basis of $V_{v_{j}}(0)$.

To verify this claim, it suffices to show that
$$
z_{1}^{j} \perp v_{1}, \quad z_{1}^{j} \perp y_{1}, \quad x_{j} \perp y_{j}, \quad z_{j} \perp y_{j}, \quad z_{j} \perp v_{j}.
$$
In fact, by using (\ref{9.2}) and (\ref{9.3}), we obtain that
$$
\begin{aligned}
&\tau G\left(z_{1}^{j}, v_{1}\right)=G\left(L\left(z_{1}^{j}, v_{j}\right), L\left(v_{j}, v_{1}\right)\right)=G\left(L\left(y_{j}, x_{1}\right), L\left(y_{j}, y_{1}\right)\right)=0, \\
&\tau G\left(z_{1}^{j}, y_{1}\right)=G\left(L\left(z_{1}^{j}, v_{j}\right), L\left(v_{j}, y_{1}\right)\right)=G\left(L\left(y_{j}, x_{1}\right),-L\left(y_{j}, v_{1}\right)\right)=0, \\
&\tau G\left(x_{j}, y_{j}\right)=G\left(L\left(x_{j}, v_{1}\right), L\left(y_{j}, v_{1}\right)\right)=G\left(L\left(v_{j}, x_{1}\right), L\left(v_{j}, y_{1}\right)\right)=0, \\
&\tau G\left(z_{j}, y_{j}\right)=G\left(L\left(z_{j}, v_{1}\right), L\left(y_{j}, v_{1}\right)\right)=G\left(L\left(x_{j}, y_{1}\right),-L\left(v_{j}, y_{1}\right)\right)=0, \\
&\tau G\left(z_{j}, v_{j}\right)=G\left(L\left(z_{j}, v_{1}\right), L\left(v_{j}, v_{1}\right)\right)=G\left(L\left(x_{j}, y_{1}\right), L\left(x_{j}, x_{1}\right)\right)=0.
\end{aligned}
$$
From these relations, we immediately get the claim.

From (\ref{9.2}) and (\ref{9.3}), $j\neq 1$, one get
$$
\left\{\begin{array}{l}
L\left(v_{j}, v_{1}\right)=L\left(x_{j}, x_{1}\right)=L\left(y_{j}, y_{1}\right), \\
L\left(v_{j}, x_{1}\right)=-L\left(x_{j}, v_{1}\right)=-L\left(y_{j}, z_{1}^{j}\right)=L\left(z_{j}, y_{1}\right), \\
L\left(v_{j}, y_{1}\right)=-L\left(y_{j}, v_{1}\right), \\
L\left(v_{j}, z_{1}^{j}\right)=L\left(y_{j}, x_{1}\right),\\
L\left(x_{j}, y_{1}\right)=L\left(z_{j}, v_{1}\right).
\end{array}\right.
$$
Next, the above relations and (\ref{5.8}) show that
\begin{align*}
G(L(y_1, y_j),L(z_1^j, z_j))&=-G(L(y_1, z_j),L(z_1^j, y_j))=G(L(y_1, z_j),L(y_1, z_j))=\tau, \\
G(L(y_1, y_j),L(y_1, y_j))&=G(L(z_1^j, z_j),L(z_1^j, z_j))=\tau.
\end{align*}
Then, by Cauchy-Schwarz inequality, we have
\be\label{9.4}
L(z_1^j, z_j)=L(y_1, y_j)=L(x_1, x_j).
\ee
Similarly
\begin{align*}
G(L(y_1, v_j),L(z_j, x_1))&=-G(L(y_1, z_j),L(v_j, x_1))=-G(L(v_j, x_1),L(v_j, x_1))=-\tau, \\
G(L(y_1, v_j),L(y_1, v_j))&=G(L(z_j, x_1),L(z_j, x_1))=\tau.
\end{align*}
Then, by Cauchy-Schwarz inequality, we have
\be\label{9.5}
L(y_1, v_j)=-L(z_j, x_1).
\ee
From (\ref{5.8}) and (\ref{9.4})
\begin{align*}
G(L(z_j, x_1),L(x_j, z_1^j))&=-G(L(z_j, z_1^j),L(x_1, x_j))=-G(L(x_1, x_j),L(x_1, x_j))=-\tau, \\
G(L(z_j, x_1),L(z_j, x_1))&=G(L(x_j, z_1^j),L(x_j, z_1^j))=\tau.
\end{align*}
Then, by Cauchy-Schwarz inequality, we have
\be\label{9.6}
L(z_j, x_1)=-L(x_j, z_1^j).
\ee
By (\ref{5.8}) and (\ref{9.5})
\begin{align*}
G(L(y_j, x_1),L(z_j, v_1))&=-G(L(y_j, v_1),L(x_1, z_j))=G(L(y_j, v_1),-L(y_j, v_1))=-\tau, \\
G(L(y_j, x_1),L(y_j, x_1))&=G(L(z_j, v_1),L(z_j, v_1))=\tau.
\end{align*}
Then, by Cauchy-Schwarz inequality, we have
\be\label{9.7}
L(y_j, x_1)=-L(z_j, v_1).
\ee
Hence, for $j>1$, we have
\be\label{9.8}
\left\{\begin{array}{l}
L\left(v_{j}, v_{1}\right)=L\left(x_{j}, x_{1}\right)=L\left(y_{j}, y_{1}\right)=L\left(z_{j}, z_{1}^{j}\right), \\
L\left(v_{j}, x_{1}\right)=-L\left(x_{j}, v_{1}\right)=-L\left(y_{j}, z_{1}^{j}\right)=L\left(z_{j}, y_{1}\right), \\
L\left(v_{j}, y_{1}\right)=-L\left(y_{j}, v_{1}\right)=-L\left(z_{j}, x_{1}\right)=L\left(x_{j}, z_{1}^{j}\right), \\
L\left(v_{j}, z_{1}^{j}\right)=-L\left(z_{j}, v_{1}\right)=-L\left(x_{j}, y_{1}\right)=L\left(y_{j}, x_{1}\right).
\end{array}\right.
\ee
From these relations we can prove the following assertion:

\textbf{Claim 2.} $z_{1}^{2}=\cdots=z_{1}^{k_{0}}=: z_{1}$.

In fact, by Claim 1 , we know that for $j \neq l$ and $j, l \geq 2$ we have $z_{1}^{j}=\varepsilon_{j l} z_{1}^{l}$ with $\varepsilon_{j l}=\pm 1$. From (\ref{5.23}) and (\ref{9.8}) we get
$$
\begin{aligned}
\varepsilon_{j l} \tau L\left(v_{j}, v_{l}\right) &=\varepsilon_{j l}A\left(L\left(z_{1}^{j}, v_{j}\right), L\left(z_{1}^{j}, v_{l}\right)\right)=A\left(L\left(z_{1}^{j}, v_{j}\right), L\left(z_{1}^{l}, v_{l}\right)\right) \\
&=A\left(L\left(y_{j}, x_{1}\right), L\left(y_{l}, x_{1}\right)\right)=\tau L\left(y_{j}, y_{l}\right).
\end{aligned}
$$
Similarly, we get
$$
\varepsilon_{j l}\tau L\left(x_{j}, x_{l}\right)=A\left(L\left(z_{1}^{j}, x_{j}\right), L\left(z_{1}^{l}, x_{l}\right)\right)=
\left\{\begin{array}{l}
A\left(L\left(z_{j}, x_{1}\right), L\left(z_{l}, x_{1}\right)\right)=\tau L\left(z_{j}, z_{l}\right), \\
A\left(L\left(y_{j}, v_{1}\right), L\left(y_{l}, v_{1}\right)\right)=\tau L\left(y_{j}, y_{l}\right), \\
A\left(L\left(v_{j}, y_{1}\right), L\left(v_{l}, y_{1}\right)\right)=\tau L\left(v_{j}, v_{l}\right).
\end{array}\right.
$$
The above two relations show that $\varepsilon_{j l}=1$. Thus Claim 2 is verified.

Moreover, the following relations hold
\be\label{9.9}
L\left(v_{j}, v_{l}\right)=L\left(x_{j}, x_{l}\right)=L\left(y_{j}, y_{l}\right)=L\left(z_{j}, z_{l}\right), \quad j \neq l \text{ and } j, l \geq 2.
\ee
From (\ref{5.23}) and (\ref{9.8}), we get
\begin{align}
\label{9.10}
\tau L\left(v_{j}, x_{l}\right)=A\left(L\left(x_{1}, v_{j}\right), L\left(x_{1}, x_{l}\right)\right)=
\left\{\begin{array}{l}
A\left(-L\left(x_{j}, v_{1}\right), L\left(v_{l}, v_{1}\right)\right)=-\tau L\left(x_{j}, v_{l}\right), \\
A\left(-L\left(y_{j}, z_{1}\right), L\left(z_{l}, z_{1}\right)\right)=-\tau L\left(y_{j}, z_{l}\right), \\
A\left(L\left(z_{j}, y_{1}\right), L\left(y_{l}, y_{1}\right)\right)=\tau L\left(z_{j}, y_{l}\right).
\end{array}\right.
\end{align}
Similarly, we have the following relations:
\begin{align}
\label{9.11}
\tau L\left(v_{j}, y_{l}\right)=A\left(L\left(y_{1}, v_{j}\right), L\left(y_{1}, y_{l}\right)\right)=
\left\{\begin{array}{l}
A\left(-L\left(y_{j}, v_{1}\right), L\left(v_{l}, v_{1}\right)\right)=-\tau L\left(y_{j}, v_{l}\right), \\
A\left(-L\left(z_{j}, x_{1}\right), L\left(x_{l}, x_{1}\right)\right)=-\tau L\left(z_{j}, x_{l}\right), \\
A\left(L\left(x_{j}, z_{1}\right), L\left(z_{l}, z_{1}\right)\right)=\tau L\left(x_{j}, z_{l}\right).
\end{array}\right.    \\
\label{9.12}
\tau L\left(v_{j}, z_{l}\right)=A\left(L\left(z_{1}, v_{j}\right), L\left(z_{1}, z_{l}\right)\right)=
\left\{\begin{array}{l}
A\left(-L\left(z_{j}, v_{1}\right), L\left(v_{l}, v_{1}\right)\right)=-\tau L\left(z_{j}, v_{l}\right), \\
A\left(-L\left(x_{j}, y_{1}\right), L\left(y_{l}, y_{1}\right)\right)=-\tau L\left(x_{j}, y_{l}\right), \\
A\left(L\left(y_{j}, x_{1}\right), L\left(x_{l}, x_{1}\right)\right)=\tau L\left(y_{j}, x_{l}\right).
\end{array}\right.
\end{align}

Combination of Claim 2 and (\ref{9.8})-(\ref{9.12}), we get (\ref{9.1}) immediately. \hfill $\Box$

By the symmetry of $L$ and Lemma \ref{lemma9.2}, we get
\be\label{9.13}
\operatorname{dim}(\operatorname{Im} L)\leq k_0+\cdots+1+3((k_0-1)+\cdots+1)=\fr{1}{8}(m-1)(m-3).
\ee

To continue the proof of Theorem \ref{thm9.1}, we now assume that $k_{0} \geq 2$ and let $\left\{v_{1}, x_{1}, y_{1}, z_{1} ; \ldots ; v_{k_{0}}, x_{k_{0}}, y_{k_{0}}, z_{k_{0}}\right\}$ be the orthonormal basis of $\mathcal{D}_{2}$ as constructed in Lemma \ref{lemma9.2}. According to Lemmas \ref{lemma5.3}, \ref{lemma5.10} and the fact that for $j \neq l, v_{j}, x_{j}, y_{j}, z_{j} \in V_{v_{l}}(\tau)=V_{x_{l}}(\tau)=V_{y_{l}}(\tau)=V_{z_{l}}(\tau)$, we have
\begin{align}
\label{9.14}
&G\left(L(v_{j}, x_{l}), L(v_{j}, x_{l})\right)=G\left(L(v_{j}, y_{l}), L(v_{j}, y_{l})\right)=G\left(L(v_{j}, z_{l}), L(v_{j}, z_{l})\right)   \notag \\
&=G\left(L(v_{j}, v_{l}), L(v_{j}, v_{l})\right)=\tau, \quad \quad j \neq l,       \\
\label{9.15}
&G\left(L(v_{j}, x_{l}), L(v_{j}, x_{k})\right)=G\left(L(v_{j}, y_{l}), L(v_{j}, y_{k})\right)=G\left(L(v_{j}, z_{l}), L(v_{j}, z_{k})\right)   \notag \\
&=G\left(L(v_{j}, v_{l}), L(v_{j}, v_{k})\right)=G\left(L(x_{j}, v_{l}), L(x_{j}, v_{k})\right)=G\left(L(y_{j}, v_{l}), L(y_{j}, v_{k})\right) \notag \\
&=G\left(L(z_{j}, v_{l}), L(z_{j}, v_{k})\right)=0, \quad \quad j, k, l \, \text{distinct},   \\
\label{9.16}
&G\left(L(v_{i}, x_{j}), L(v_{k}, x_{l})\right)=G\left(L(v_{i}, y_{j}), L(v_{k}, y_{l})\right)=G\left(L(v_{i}, z_{j}), L(v_{k}, z_{l})\right)   \notag \\
&=G\left(L(v_{i}, v_{j}), L(v_{k}, v_{l})\right)=0, \quad \quad i, j, k, l\;\text{distinct},\\
\label{9.17}
&G\left(L\left(v_{j}, v_{l}\right), L\left(v_{j_{1}}, x_{l_{1}}\right)\right)
=G\left(L\left(v_{j}, v_{l}\right), L\left(v_{j_{1}}, y_{l_{1}}\right)\right) \notag \\
&=G\left(L\left(v_{j}, v_{l}\right), L\left(v_{j_{1}}, z_{l_{1}}\right)\right)=0, \quad j \neq l \text { and } j_{1} \neq l_{1}, \\
\label{9.18}
&G\left(L\left(v_{j}, v_{j}\right), L\left(v_{j}, v_{j}\right)\right)=\frac{1}{4} \mu_{1}^2, \quad 1 \leq j \leq k_{0}, \\
\label{9.19}
&G\left(L\left(v_{j}, v_{j}\right), L\left(v_{l}, v_{l}\right)\right)=0, \quad j \neq l, \\
\label{9.20}
&G\left(L\left(v_{j}, v_{j}\right), L\left(v_{j}, v_{l}\right)\right)=G\left(L\left(v_{j}, v_{j}\right), L\left(v_{j}, x_{l}\right)\right)=G\left(L\left(v_{j}, v_{j}\right), L\left(v_{j}, y_{l}\right)\right)  \notag    \\
&=G\left(L\left(v_{j}, v_{j}\right), L\left(v_{j}, z_{l}\right)\right)=G\left(L\left(v_{j}, v_{j}\right), L\left(v_{l}, x_{j}\right)\right)=G\left(L\left(v_{j}, v_{j}\right), L\left(v_{l}, y_{j}\right)\right)  \notag  \\
&=G\left(L\left(v_{j}, v_{j}\right), L\left(v_{l}, z_{j}\right)\right)=0, \quad j \neq l, \\
\label{9.21}
&G\left(L\left(v_{j}, v_{j}\right), L\left(v_{l_{1}}, v_{l_{2}}\right)\right)
=G\left(L\left(v_{j}, v_{j}\right), L\left(v_{l_{1}}, x_{l_{2}}\right)\right)
=G\left(L\left(v_{j}, v_{j}\right), L\left(v_{l_{1}}, y_{l_{2}}\right)\right) \notag  \\
&=G\left(L\left(v_{j}, v_{j}\right), L\left(v_{l_{1}}, z_{l_{2}}\right)\right)=0, \quad j, l_{1}, l_{2}\,\text{distinct}.
\end{align}

As in preceding sections we denote
$$
L_{j}:=L\left(v_{1}, v_{1}\right)+\cdots+L\left(v_{j}, v_{j}\right)-j L\left(v_{j+1}, v_{j+1}\right), \quad 1 \leq j \leq k_{0}-1.
$$
(\ref{9.18}) and (\ref{9.19}) show that $G\left(L_{j}, L_{j}\right)=2 j(j+1) \tau \neq 0$ for each $j$. Moreover (\ref{9.14})-(\ref{9.21}) and Lemmas \ref{lemma5.10}, \ref{lemma9.2} show that
\be\label{9.22}
\left\{\begin{array}{cl}
w_{j}=\frac{1}{\sqrt{2 j(j+1) \tau}} L_{j}, & 1 \leq j \leq k_{0}-1,\\
w_{k l}=\frac{1}{\sqrt{\tau}} L\left(v_{k}, v_{l}\right), & 1 \leq k<l \leq k_{0}, \\
w_{k l}^{\prime}=\frac{1}{\sqrt{\tau}} L\left(v_{k}, x_{l}\right), & 1 \leq k<l \leq k_{0}, \\ w_{k l}^{''}=\frac{1}{\sqrt{\tau}} L\left(v_{k}, y_{l}\right), & 1 \leq k<l \leq k_{0}, \\
w_{k l}^{'''}=\frac{1}{\sqrt{\tau}} L\left(v_{k}, z_{l}\right), & 1 \leq k<l \leq k_{0}
\end{array}\right.
\ee
give $\frac{1}{8}(m+1)(m-5)$ mutually orthogonal unit vectors in $\operatorname{Im} L \subset \mathcal{D}_{3}$.

Direct computations and Lemma \ref{lemma9.2} show that $\operatorname{Tr} L=L(v_1,v_1)+L(x_1,x_1)+L(y_1,y_1)+L(z_1,z_1)+\cdots+L(v_{k_0},v_{k_0})+L(x_{k_0},x_{k_0})
+L(y_{k_0},y_{k_0})+L(z_{k_0},z_{k_0})=4\left[L\left(v_{1}, v_{1}\right)+\cdots+L\left(v_{k_{0}}, v_{k_{0}}\right)\right]\in \operatorname{Im} L$ is orthogonal to all vectors in (\ref{9.22}), and by using (\ref{9.18}) and (\ref{9.19}), we get
\be\label{9.23}
\frac{1}{16} G(\operatorname{Tr} L, \operatorname{Tr} L)=\frac{1}{4} k_{0}\mu_1^2=: \rho^{2}>0.   \ee

(\ref{9.22}) and (\ref{9.23}) give $\frac{1}{8}(m+1)(m-5)+1=\fr{1}{8}(m-1)(m-3)$ mutually orthogonal vectors in $\operatorname{Im} L \subset \mathcal{D}_{3}$. Combining (\ref{9.13}) one can get $\operatorname{dim}(\operatorname{Im} L)=\fr{1}{8}(m-1)(m-3)$.
Thus, we have the estimate of the dimension
$$
\begin{aligned}
n &=1+\operatorname{dim}\left(\mathcal{D}_{2}\right)+
\operatorname{dim}\left(\mathcal{D}_{3}\right) \\
& \geq 1+m-1+\fr{1}{8}(m-1)(m-3)=\fr{1}{8}(m+1)(m+3).
\end{aligned}
$$

Now, we are ready to complete the proof of Theorem \ref{thm9.1}.\\
\textbf{Proof of Theorem \ref{thm9.1}.} We need to consider two cases:

Case (1).\ $n=\fr{1}{8}(m+1)(m+3)$,\ ($\operatorname{Im} L= \mathcal{D}_{3}$).

Case (2).\ $n>\fr{1}{8}(m+1)(m+3)$,\ ($\operatorname{Im} L \varsubsetneqq \mathcal{D}_{3}$).\\
For both cases (1) and (2), as $\operatorname{Tr} L\neq0$, we can define a unit vector $t:=\fr{1}{4\rho}\operatorname{Tr} L$.

In Case (1), from the previous discussions, we see that
$$\{t, w_j|_{1\leq j\leq k_0-1}, w_{kl}|_{1\leq k<l\leq k_0}, w_{k l}'|_{1\leq k<l\leq k_0}, w_{k l}''|_{1\leq k<l\leq k_0}, w_{k l}'''|_{1\leq k<l\leq k_0}\}$$
forms an orthonormal basis of $\operatorname{Im} L=\mathcal{D}_{3}$.

In Case (2), we choose that $\{\td w_i|_{1\leq i\leq \td n}\}$
forms an orthonormal basis of $\mathcal{D}_{3}\backslash \operatorname{Im} L$, and such that
$$\{t, w_j|_{1\leq j\leq k_0-1}, w_{kl}|_{1\leq k<l\leq k_0}, w_{k l}'|_{1\leq k<l\leq k_0}, w_{k l}''|_{1\leq k<l\leq k_0}, w_{k l}'''|_{1\leq k<l\leq k_0}, \td w_i|_{1\leq i\leq \td n}\}$$
forms an orthonormal basis of $\mathcal{D}_{3}$,  where $\td n=n-\fr{1}{8}(m+1)(m+3)$.

\begin{lem}\label{lemma9.3}
The Fubini-Pick tensor satisfies
$$\label{9.24}
\left\{\begin{array}{lll}
A(t, e_1)=0,\; A(t,t)=\fr{2\rho}{k_0}t,
\\ A(t,v_i)=\fr{\rho}{k_0}v_i, \; A(t,x_i)=\fr{\rho}{k_0}x_i, \; A(t,y_i)=\fr{\rho}{k_0}y_i,\;  A(t,z_i)=\fr{\rho}{k_0}z_i, \; \ 1\leq i\leq k_0, \\ A(t,w_j)=\fr{2\rho}{k_0}w_j, \; \ 1\leq j\leq k_0-1,\\
A(t,w_{kl})=\fr{2\rho}{k_0}w_{kl}, \; A(t,w_{kl}')=\fr{2\rho}{k_0}w_{kl}', \; A(t,w_{kl}'')=\fr{2\rho}{k_0}w_{kl}'', \;\\ A(t,w_{kl}''')=\fr{2\rho}{k_0}w_{kl}''',\; 1\leq k<l\leq k_0, \\
\text { if } \operatorname{Im} L \varsubsetneqq \mathcal{D}_{3}, \quad A(t,\td w_i)=0, \quad 1\leq i\leq \td n.
\end{array}
\right.
$$
\end{lem}

\textbf{Proof.}
By (iii) of Lemma \ref{lemma5.1}, we know $A(t,v_i)\in \mathcal{D}_2$. (\ref{9.14})-(\ref{9.21}) show that
\begin{align*}
\lang A(t,v_i),v_j\rang&=\fr{1}{\rho}\sum_{l=1}^{k_0}\lang A(v_i,v_j),L(v_l,v_l)\rang=\fr{\rho}{k_0}\delta_{ij},\quad 1\leq i,j\leq k_0,\\
\lang A(t,v_i),x_j\rang&=\fr{1}{\rho}\sum_{l=1}^{k_0}\lang A(v_i,x_j),L(v_l,v_l)\rang=0,\quad 1\leq i,j\leq k_0, \\
\lang A(t,v_i),y_j\rang&=\fr{1}{\rho}\sum_{l=1}^{k_0}\lang A(v_i,y_j),L(v_l,v_l)\rang=0,\quad 1\leq i,j\leq k_0, \\
\lang A(t,v_i),z_j\rang&=\fr{1}{\rho}\sum_{l=1}^{k_0}\lang A(v_i,z_j),L(v_l,v_l)\rang=0,\quad 1\leq i,j\leq k_0.
\end{align*}
Hence, we have $A(t,v_i)=\fr{\rho}{k_0}v_i, \ 1\leq i\leq k_0$. Similarly, $A(t,x_i)=\fr{\rho}{k_0}x_i,\ A(t,y_i)=\fr{\rho}{k_0}y_i,\ A(t,z_i)=\fr{\rho}{k_0}z_i,\ 1\leq i\leq k_0$ can be proved.

From Lemmas \ref{lemma5.14} and \ref{lemma9.2} we get
\begin{align*}
A(t,t)&=\fr{1}{\rho^2} \sum_{i=1}^{k_0}A(L(v_i,v_i),L(v_i,v_i))=\fr{2\rho}{k_0}t, \\
A(t,w_j)&=\fr{\mu_1^2}{2\rho\sqrt{2j(j+1)\tau}}L_j=\fr{2\rho}{k_0}w_j, \quad 1\leq j\leq k_0-1,\\
A(t,w_{kl})&=\fr{1}{\rho\sqrt \tau}(A(L(v_k,v_k),L(v_k,v_l))+A(L(v_l,v_l),L(v_l,v_k)))=\fr{2\rho}{k_0}w_{kl}, \\
A(t,w_{kl}')&=\fr{1}{\rho\sqrt \tau}(A(L(v_k,v_k),L(v_k,x_l))+A(L(v_l,v_l),L(x_l,v_k)))   \\
&=\fr{1}{\rho\sqrt \tau}(A(L(v_k,v_k),L(v_k,x_l))+A(L(x_l,x_l),L(x_l,v_k)))=\fr{2\rho}{k_0}w_{kl}', \end{align*}
and
$$A(t,w_{kl}'')=\fr{2\rho}{k_0}w_{kl}'', \quad A(t,w_{kl}''')=\fr{2\rho}{k_0}w_{kl}''',$$ where $1\leq k<l\leq k_0$.
If $\operatorname{Im} L \varsubsetneqq \mathcal{D}_{3}$, by Lemma \ref{lemma5.4}, one get
$$A(t,\td w_i)=\fr{1}{\rho}\sum_{j=1}^{k_0}A(L(v_j,v_j),\td w_i)=0, \quad 1\leq i\leq \td n.$$
This completes the proof of Lemma \ref{lemma9.3}.    \hfill $\Box$

Put
\be\label{9.25}
T=\fr{1}{\sqrt{k_0+1}}e_1+\sqrt{\fr{k_0}{k_0+1}}t,\quad T^*=-\sqrt{\fr{k_0}{k_0+1}}e_1+\fr{1}{\sqrt{k_0+1}}t.
\ee
It is easily to see that if $\operatorname{Im} L=\mathcal{D}_{3}$, then $\{T, T^*, v_j|_{1\leq j\leq k_0} , x_j|_{1\leq j\leq k_0}, y_j|_{1\leq j\leq k_0}, z_j|_{1\leq j\leq k_0} , \\ w_j|_{1\leq j\leq k_0-1}, w_{kl}|_{1\leq k<l\leq k_0}, w_{kl}'|_{1\leq k<l\leq k_0}, w_{kl}''|_{1\leq k<l\leq k_0}, w_{kl}'''|_{1\leq k<l\leq k_0}\}$ (or, resp. if $\operatorname{Im} L \varsubsetneqq \mathcal{D}_{3}$, then $\{T, T^*, v_j|_{1\leq j\leq k_0} , x_j|_{1\leq j\leq k_0}, y_j|_{1\leq j\leq k_0}, z_j|_{1\leq j\leq k_0} , w_j|_{1\leq j\leq k_0-1}, w_{kl}|_{1\leq k<l\leq k_0}, \\
w_{kl}'|_{1\leq k<l\leq k_0}, w_{kl}''|_{1\leq k<l\leq k_0}, w_{kl}'''|_{1\leq k<l\leq k_0}, \td w_i|_{1\leq i\leq \td n}\}$) is an orthonormal basis of $T_pM$. By Lemmas \ref{lemma5.1} and \ref{lemma9.3} we have the following
\begin{lem}\label{lemma9.4}
With respect to the above notations, the Fubini-Pick tensor takes the following form
$$\label{9.26}\left\{
\begin{array}{ll}
A(T,T)=\sigma T, \, A(T,T^*)=\sigma T^*, \\ A(T,v_i)=\sigma v_i, \, A(T,x_i)=\sigma x_i,
A(T,y_i)=\sigma y_i, \, A(T,z_i)=\sigma z_i, \, \ 1\leq i\leq k_0, \\
A(T,w_j)=\sigma w_j, \ 1\leq j\leq k_0-1,      \\
A(T,w_{kl})=\sigma w_{kl}, \quad A(T,w_{kl}')=\sigma w_{kl}', \quad A(T,w_{kl}'')=\sigma w_{kl}'', \quad 1\leq k<l\leq k_0, \\
A(T,w_{kl}''')=\sigma w_{kl}''',\quad 1\leq k<l\leq k_0, \\
\text { if } \operatorname{Im} L \varsubsetneqq \mathcal{D}_{3}, \quad A(T,\td w_i)=0, \quad 1\leq i\leq \td n,
\end{array}
\right.
$$
where $\sigma=\fr{1}{\sqrt{k_0+1}}\mu_1$.
\end{lem}

Using the parallelism of Fubini-Pick form, Lemma 9.4, Theorems 4.3 and 4.6, we complete the proof of Theorem \ref{thm9.1}.
\section{ Case $\left\{\mathfrak{C}_{m}\right\}_{2 \leq m \leq n-2}$ with $k_0\geq2$ and $\mathfrak{p}=7$}

In this section, we will prove the following theorem.

\begin{thm}\label{thm10.1}
 Let $M^{n}$ be a Calabi hypersurface in $\mathbb{R}^{n+1}$ which has parallel and non-vanishing cubic form. If $\mathfrak{C}_{m}$ with $2 \leq m \leq n-2$ accurs and the integer $k_{0}$ and $\mathfrak{p}$, satisfy $k_{0} \geq 2$ and $\mathfrak{p}=7$, then $k_{0}=2, m=17$ and $n \geq 27$. Moreover, we have either
\begin{enumerate}
\item [(i)] $n= 27$, $ M^{n}$ can be decomposed as the Calabi product of a locally strongly convex hyperbolic centroaffine hypersurface with parallel cubic form and a point, or
\item [(ii)] $n>27$, $ M^{n}$ can be decomposed as the Calabi product of a locally strongly convex hyperbolic centroaffine hypersurfaces with parallel cubic form  and a Calabi hypersurface with parallel cubic form.
      \end{enumerate}
\end{thm}

To prove Theorem \ref{thm10.1}, a key ingredient is the following lemma whose proof is similar to that of Lemma 8.1 in \cite{HLV}.

\begin{lem}\label{lemma10.2}
If in the decomposition (\ref{5.18}), $k_{0} \geq 2$ and $\mathfrak{p}=7$, then we can choose an orthonormal basis $\left\{x_{j}\right\}_{1 \leq j \leq 7}$ for $V_{v_{1}}(0)$ and an orthonormal basis $\left\{y_{j}\right\}_{1 \leq j \leq 7}$ for $V_{v_{2}}(0)$ so that by identifying $e_{j}\left(v_{1}\right)=x_{j}$ and $e_{j}\left(v_{2}\right)=y_{j}$, we have the relations
\be\label{10.1}
L\left(e_{j}\left(v_{1}\right), e_{l}\left(v_{2}\right)\right)=-L\left(v_{1}, e_{j} e_{l}\left(v_{2}\right)\right)=-L\left(e_{l} e_{j}\left(v_{1}\right), v_{2}\right)
\ee
for $1 \leq j, l \leq 7$, where $e_{j} e_{l}$ denotes a product defined by the following multiplication table.\end{lem}
$$
\begin{tabular}{|c|c|c|c|c|c|c|c|}
  \hline     & $e_1$ & $e_2$ & $e_3$ & $e_4$ & $e_5$ & $e_6$ & $e_7$ \\
\hline $e_1$ & -id & $e_3$ & $-e_2$ & $e_5$ & $-e_4$ & $-e_7$ & $e_6$ \\
\hline $ e_2$ & $-e_3$ & -id  & $e_1$ & $e_6$ & $e_7$ & $-e_4$ & $-e_5$ \\
\hline  $e_3$ & $e_2$ & $-e_1$ & -id & $e_7$ & $-e_6$ &  $e_5$ & $-e_4 $\\
 \hline $e_4$ & $-e_5$ &  $-e_6$ & $-e_7$ & -id & $e_1$ & $e_2$ & $e_3$ \\
\hline  $e_5$ & $e_4$ & $-e_7$ & $ e_6$ & $-e_1$ & -id & $-e_3$ & $e_2$ \\
\hline  $e_6$ & $e_7$ & $e_4$ & $-e_5$ & $-e_2$ & $e_3$ & -id & $-e_1$ \\
\hline  $e_7$ &  $-e_6$ & $e_5$ & $e_4$ & $-e_3$ & $-e_2$ & $e_1$  & -id \\
\hline
\end{tabular}
$$

\textbf{Proof.} As before we denote $V_{j}=\left\{v_{j}\right\} \oplus V_{v_{j}}(0), 1 \leq j \leq k_{0} .$ First we fix any two orthogonal unit vectors $x_{1}, x_{2} \in V_{v_{1}}(0)$. Then, by Lemmas \ref{lemma5.12} and \ref{lemma5.13}, we can consecutively find unit vectors $y_{1}, y_{2} \in V_{v_{2}}(0)$, such that
\begin{align}
\label{10.2}
L\left(y_{1}, v_{1}\right)&=-L\left(x_{1}, v_{2}\right), \quad  L\left(y_{1}, x_{1}\right)=L\left(v_{1}, v_{2}\right), \\
\label{10.3}
L\left(y_{2}, v_{1}\right)&=-L\left(x_{2}, v_{2}\right), \quad  L\left(y_{2}, x_{2}\right)=L\left(v_{1}, v_{2}\right).
\end{align}
As $y_{1}\in V_{v_{2}}(0)$, we find a unit vector $x_{3} \in V_{x_{2}}(0)$ such that
\be\label{10.4}
L\left(y_{1}, x_{2}\right)=-L\left(x_{3}, v_{2}\right), \quad   L\left(y_{1}, x_{3}\right)=L\left(x_{2}, v_{2}\right).   \ee
From (\ref{10.2}) and (\ref{10.4}), we have
$$
\tau G\left(x_{3}, v_{1}\right)=G\left(L\left(x_{3}, v_{2}\right), L\left(v_{1}, v_{2}\right)\right)=G\left(-L\left(y_{1}, x_{2}\right), L\left(y_{1}, x_{1}\right)\right)=0,
$$
which means that $x_{3} \perp v_{1}$. By $x_{3} \in \{x_2\} \oplus V_{x_2}(0)=\{v_1\} \oplus V_{v_1}(0)$, we get $x_{3} \in V_{v_{1}}(0)$. Thus, we can further take a unit vector $y_{3} \in V_{v_{2}}(0)$ such that
\be\label{10.5}
L\left(y_{3}, v_{1}\right)=-L\left(x_{3}, v_{2}\right), \quad L\left(y_{3}, x_{3}\right)=L\left(v_{1}, v_{2}\right).
\ee
\textbf{Claim 1.} $\left\{x_{1}, x_{2}, x_{3}, v_{1}\right\}$ are orthonormal vectors. Similarly, $\left\{y_{1}, y_{2}, y_{3}, v_{2}\right\}$ are orthonormal vectors.

To verify this claim, it suffices to show that
$$x_{1} \perp x_{3}, \quad y_{1} \perp y_{2}, \quad y_{1} \perp y_{3}, \quad y_{2} \perp y_{3}.$$
In fact, by using (\ref{10.2})-(\ref{10.5}), we have
$$
\begin{aligned}
\tau G\left(x_{3}, x_{1}\right)&=G\left(L\left(x_{3}, v_{2}\right), L\left(x_{1}, v_{2}\right)\right)=
\left\{\begin{array}{l}
G\left(L\left(y_{3}, v_{1}\right), L\left(y_{1}, v_{1}\right)\right)=\tau G\left(y_{3}, y_{1}\right),    \\
G\left(L\left(y_{1}, x_{2}\right), L\left(y_{1}, v_{1}\right)\right)=0;
\end{array}\right.   \\
\tau G\left(y_{1}, y_{2}\right)&=G\left(L\left(y_{1}, v_{1}\right), L\left(y_{2}, v_{1}\right)\right)=G\left(-L\left(x_{1}, v_{2}\right), -L\left(x_{2}, v_{2}\right)\right)=0; \\
\tau G\left(y_{2}, y_{3}\right)&=G\left(L\left(y_{2}, v_{1}\right), L\left(y_{3}, v_{1}\right)\right)=G\left(-L\left(x_{2}, v_{2}\right), -L\left(x_{3}, v_{2}\right)\right)=0.
\end{aligned}
$$
Hence we have the Claim 1.

By (\ref{10.2}),(\ref{10.3}) and (\ref{10.5}), we get the relation
\be\label{10.6}
L\left(y_{1}, x_{1}\right)=L\left(y_{2}, x_{2}\right)=L\left(y_{3}, x_{3}\right)=L\left(v_{1}, v_{2}\right).
\ee
Claim 1 shows that $x_1\in V_{x_2}(0)$ and $y_1\in V_{y_2}(0)$, which together with $L\left(y_{1}, x_{1}\right)=L\left(y_{2}, x_{2}\right)$ in (\ref{10.6}) and Lemma \ref{lemma5.13}, imply that
\be\label{10.7}
L\left(x_{1}, y_{2}\right)=-L\left(y_{1}, x_{2}\right)\overset{(\ref{10.4})}=L\left(x_{3}, v_{2}\right)\overset{(\ref{10.5})}=-L\left(y_{3}, v_{1}\right). \ee
Similarly, $L\left(y_{1}, x_{1}\right)=L\left(y_{3}, x_{3}\right)$ and $L\left(y_{2}, x_{2}\right)=L\left(y_{3}, x_{3}\right)$ in (\ref{10.6}) show that
\begin{align}
\label{10.8}
&L\left(x_{1}, y_{3}\right)=-L\left(y_{1}, x_{3}\right)\overset{(\ref{10.4})}=-L\left(x_{2}, v_{2}\right)\overset{(\ref{10.3})}=L\left(y_{2}, v_{1}\right), \\
\label{10.9}
&L\left(y_{3}, x_{2}\right)=-L\left(x_{3}, y_{2}\right)= -L\left(y_{1}, v_{1}\right)
\overset{(\ref{10.2})}=L\left(x_{1}, v_{2}\right).
\end{align}
The second "="  in (\ref{10.9}) holds because the following fact
$$G(L(x_{3}, y_{2}),L(y_{1}, v_{1}))\overset{(\ref{5.8})}=-G(L(x_{3}, y_{1}),L(y_{2}, v_{1}))\overset{(\ref{10.8})}=G(L(x_{2}, v_{2}),L(x_{2}, v_{2}))=\tau$$
and
$$G(L(x_{3}, y_{2}),L(x_{3}, y_{2}))=G(L(y_{1}, v_{1}),L(y_{1}, v_{1}))=\tau. $$

Now we pick an arbitrary unit vector $x_{4} \in V_{v_{1}}(0)$ such that it is orthogonal to all $x_{1}, x_{2}$ and $x_{3}$. As $x_4\in \{v_1\}\oplus V_{v_1}(0)=\{x_1\}\oplus V_{x_1}(0)$, we have $x_4\in V_{x_1}(0)$. It follows Lemma \ref{lemma5.13} that there exists a unit vector
$y_4\in V_{y_1}(0)$ such that
\be\label{10.10}
L\left(x_{4}, y_{1}\right)=-L\left(x_{1}, y_{4}\right), \quad
L\left(x_{4}, y_{4}\right)=L\left(x_{1}, y_{1}\right)\overset{(\ref{10.6})}=L\left(v_{1}, v_{2}\right).
\ee
By
$$\tau G(y_4,v_2)=G(L(y_4,v_1),L(v_2,v_1))\overset{(\ref{10.10})}
=G(L(y_4,v_1),L(y_4,x_4))=0$$
and $y_4\in \{y_1\}\oplus V_{y_1}(0)=\{v_2\}\oplus V_{v_2}(0)$, we have $y_4\in V_{v_2}(0)$.

The fact $y_1\in V_{v_2}(0)$ shows that there exists a unit vector
$x_5\in V_{x_4}(0)$ such that
\be\label{10.11}
L\left(x_{4}, y_{1}\right)=-L\left(x_{5}, v_{2}\right), \quad
L\left(x_{4}, v_{2}\right)=L\left(x_{5}, y_{1}\right).
\ee
As
$$\tau G(x_5,v_1)=G(L(x_5,v_2),L(v_1,v_2))\overset{(\ref{10.11})(\ref{10.6})}
=G(-L(y_1,x_4),L(y_1,x_1))=0$$
and $x_5\in \{x_4\}\oplus V_{x_4}(0)=\{v_1\}\oplus V_{v_1}(0)$, we have $x_5\in V_{v_1}(0)$. It follows Lemma \ref{lemma5.13} that there exists a unit vector
$y_5\in V_{v_2}(0)$ such that
\be\label{10.12}
L\left(x_{5}, v_{2}\right)=-L\left(y_{5}, v_{1}\right), \quad
L\left(x_{5}, y_{5}\right)=L\left(v_{1}, v_{2}\right).
\ee

The fact $y_2\in V_{v_2}(0)$ shows that there exists a unit vector
$x_6\in V_{x_4}(0)$ such that
\be\label{10.13}
L\left(x_{4}, y_{2}\right)=-L\left(x_{6}, v_{2}\right), \quad
L\left(x_{6}, y_{2}\right)=L\left(x_{4}, v_{2}\right).
\ee
As
$$\tau G(x_6,v_1)=G(L(x_6,v_2),L(v_1,v_2))\overset{(\ref{10.13})(\ref{10.6})}
=G(-L(y_2,x_4),L(y_2,x_2))=0$$
and $x_6\in \{x_4\}\oplus V_{x_4}(0)=\{v_1\}\oplus V_{v_1}(0)$, we have $x_6\in V_{v_1}(0)$. It follows Lemma \ref{lemma5.13} that there exists a unit vector
$y_6\in V_{v_2}(0)$ such that
\be\label{10.14}
L\left(x_{6}, v_{2}\right)=-L\left(y_{6}, v_{1}\right), \quad
L\left(x_{6}, y_{6}\right)=L\left(v_{1}, v_{2}\right).
\ee

The fact $y_3\in V_{v_2}(0)$ shows that there exists a unit vector
$x_7\in V_{x_4}(0)$ such that
\be\label{10.15}
L\left(x_{4}, y_{3}\right)=-L\left(x_{7}, v_{2}\right), \quad
L\left(x_{7}, y_{3}\right)=L\left(x_{4}, v_{2}\right).
\ee
As
$$\tau G(x_7,v_1)=G(L(x_7,v_2),L(v_1,v_2))\overset{(\ref{10.15})(\ref{10.6})}
=G(-L(y_3,x_4),L(y_3,x_3))=0$$
and $x_7\in \{x_4\}\oplus V_{x_4}(0)=\{v_1\}\oplus V_{v_1}(0)$, we have $x_7\in V_{v_1}(0)$. It follows Lemma \ref{lemma5.13} that there exists a unit vector
$y_7\in V_{v_2}(0)$ such that
\be\label{10.16}
L\left(x_{7}, v_{2}\right)=-L\left(y_{7}, v_{1}\right), \quad
L\left(x_{7}, y_{7}\right)=L\left(v_{1}, v_{2}\right).
\ee

Then, by (\ref{10.10})-(\ref{10.16}) , we can find unit vectors $x_{5}, x_{6}, x_{7} \in V_{v_{1}}(0)$ and $y_{4}, y_{5}, y_{6}$, $y_{7} \in V_{v_{2}}(0)$ such that the following relations hold:
\begin{align}
\label{10.17}
&\left\{\begin{array}{l}
L\left(x_{4}, y_{1}\right)=-L\left(x_{1}, y_{4}\right)=-L\left(x_{5}, v_{2}\right)=L\left(y_{5}, v_{1}\right), \\
L\left(x_{4}, y_{4}\right)=L\left(x_{1}, y_{1}\right)=L\left(x_{5}, y_{5}\right)=L\left(v_{1}, v_{2}\right),\; L\left(x_{4}, v_{2}\right)=L\left(x_{5}, y_{1}\right),
\end{array}\right. \\
\label{10.18}
&\left\{\begin{array}{l}
L\left(x_{4}, y_{2}\right)=-L\left(x_{6}, v_{2}\right)=L\left(y_{6}, v_{1}\right), \\
L\left(x_{4}, v_{2}\right)=L\left(x_{6}, y_{2}\right),\; L\left(x_{6}, y_{6}\right)=L\left(v_{1}, v_{2}\right),
\end{array}\right. \\
\label{10.19}
&\left\{\begin{array}{l}
L\left(x_{4}, y_{3}\right)=-L\left(x_{7}, v_{2}\right)=L\left(y_{7}, v_{1}\right), \\
L\left(x_{4}, v_{2}\right)=L\left(x_{7}, y_{3}\right),\; L\left(x_{7}, y_{7}\right)=L\left(v_{1}, v_{2}\right).
\end{array}\right.
\end{align}

\textbf{Claim 2.} $\left\{x_{1}, \ldots, x_{7}, v_{1}\right\}$ are orthonormal vectors. Similarly, $\left\{y_{1}, \ldots, y_{7}, v_{2}\right\}$ are orthonormal vectors.

From \textbf{Claim 1} and the proof of above, we know that $\left\{x_{1}, x_{2}, x_{3}, x_{4}, v_{1}\right\}$, $\left\{y_{1}, y_{2}, y_{3}, v_{2}\right\}$ are orthonormal vectors and  $x_5\perp v_1$, $x_5\perp x_4$, $x_6\perp v_1$, $x_6\perp x_4$, $x_7\perp v_1$, $x_7\perp x_4$, $y_4\perp y_1$, $y_4\perp v_2$, $y_5\perp v_2$, $y_6\perp v_2$, $y_7\perp v_2$. It suffices to show that
\begin{align*}
\tau G(x_1,x_5)=G(L(x_5,v_2),L(x_1,v_2))
\left\{\begin{array}{l}
\overset{(\ref{10.17})(\ref{10.2})}=G(L(y_5,v_1),L(y_1,v_1))=\tau G(y_1,y_5),  \\
\overset{(\ref{10.17})(\ref{10.2})}=G(L(x_4,y_1),L(y_1,v_1))=0,
\end{array}\right.
\end{align*}
\begin{align*}
\tau G(x_1,x_6)=G(L(x_6,v_2),L(x_1,v_2))
\left\{\begin{array}{l}
\overset{(\ref{10.18})}=G(-L(y_6,v_1),-L(y_1,v_1))=\tau G(y_6,y_1), \\
\overset{(\ref{10.18})}=G(-L(x_4,y_2),-L(y_2,x_3))=0,
\end{array}\right.
\end{align*}
(As $L(x_1,y_2)=L(x_3,v_2)$ in (\ref{10.7}) and Lemma \ref{lemma5.13}, we get $L(x_1,v_2)=-L(y_2,x_3)$.)
\begin{align*}
\tau G(x_1,x_7)=G(L(x_7,v_2),L(x_1,v_2))
\left\{\begin{array}{l}
\overset{(\ref{10.19})}=G(-L(y_7,v_1),-L(y_1,v_1))=\tau G(y_7,y_1),  \\
\overset{(\ref{10.19})}=G(-L(x_4,y_3),L(y_3,x_2))=0,
\end{array}\right.
\end{align*}
(As $L(x_1,y_3)=-L(x_2,v_2)$ in (\ref{10.8}) and Lemma \ref{lemma5.13}, we get $L(x_1,v_2)=L(x_2,y_3)$.)
\begin{align*}
\tau G(x_2,x_5)=G(L(x_2,v_2),L(x_5,v_2))
\left\{\begin{array}{l}
\overset{(\ref{10.3})(\ref{10.17})}=G(-L(y_2,v_1),-L(y_5,v_1))=\tau G(y_2,y_5),  \\
\overset{(\ref{10.4})(\ref{10.17})}=G(L(y_1,x_3),-L(x_4,y_1))=0,
\end{array}\right.
\end{align*}
\begin{align*}
\tau G(x_2,x_6)=G(L(x_2,v_2),L(x_6,v_2))
\left\{\begin{array}{l}
\overset{(\ref{10.8})(\ref{10.18})}=G(-L(y_2,v_1),-L(y_6,v_1))=\tau G(y_2,y_6),  \\
\overset{(\ref{10.8})(\ref{10.18})}=G(-L(y_2,v_1),-L(x_4,y_2))=0,
\end{array}\right.
\end{align*}
\begin{align*}
\tau G(x_2,x_7)=G(L(x_2,v_2),L(x_7,v_2))
\left\{\begin{array}{l}
\overset{(\ref{10.8})(\ref{10.19})}=G(-L(y_2,v_1),-L(y_7,v_1))=\tau G(y_2,y_7),  \\
\overset{(\ref{10.8})(\ref{10.19})}=G(-L(x_1,y_3),-L(x_4,y_3))=0,
\end{array}\right.
\end{align*}
\begin{align*}
\tau G(x_3,x_5)=G(L(x_3,v_2),L(x_5,v_2))
\overset{(\ref{10.4})(\ref{10.17})}=G(-L(y_1,x_2),-L(x_4,y_1))=0,
\end{align*}
\begin{align*}
\tau G(x_3,x_6)=G(L(x_3,v_2),L(x_6,v_2))
\left\{\begin{array}{l}
\overset{(\ref{10.7})(\ref{10.18})}=G(-L(y_3,v_1),-L(y_6,v_1))=\tau G(y_3,y_6),  \\
\overset{(\ref{10.7})(\ref{10.18})}=G(L(x_1,y_2),-L(x_4,y_2))=0,
\end{array}\right.
\end{align*}
\begin{align*}
\tau G(x_3,x_7)=G(L(x_3,v_2),L(x_7,v_2))
\left\{\begin{array}{l}
\overset{(\ref{10.7})(\ref{10.19})}=G(-L(y_3,v_1),-L(y_7,v_1))=\tau G(y_3,y_7),  \\
\overset{(\ref{10.7})(\ref{10.19})}=G(-L(y_3,v_1),-L(x_4,y_3))=0,
\end{array}\right.
\end{align*}
\begin{align*}
\tau G(x_5,x_6)=G(L(x_5,v_2),L(x_6,v_2))
\left\{\begin{array}{l}
\overset{(\ref{10.17})(\ref{10.18})}=G(-L(y_5,v_1),-L(y_6,v_1))=\tau G(y_5,y_6),  \\
\overset{(\ref{10.17})(\ref{10.18})}=G(-L(x_4,y_1),-L(x_4,y_2))=0,
\end{array}\right.
\end{align*}
\begin{align*}
\tau G(x_5,x_7)=G(L(x_5,v_2),L(x_7,v_2))
\left\{\begin{array}{l}
\overset{(\ref{10.17})(\ref{10.19})}=G(-L(y_5,v_1),-L(y_7,v_1))=\tau G(y_5,y_7),  \\
\overset{(\ref{10.17})(\ref{10.19})}=G(-L(x_4,y_1),-L(x_4,y_3))=0,
\end{array}\right.
\end{align*}
\begin{align*}
\tau G(x_6,x_7)=G(L(x_6,v_2),L(x_7,v_2))
\left\{\begin{array}{l}
\overset{(\ref{10.18})(\ref{10.19})}=G(-L(y_6,v_1),-L(y_7,v_1))=\tau G(y_6,y_7),  \\
\overset{(\ref{10.18})(\ref{10.19})}=G(-L(x_4,y_2),-L(x_4,y_3))=0,
\end{array}\right.
\end{align*}
\begin{align*}
\tau G(y_2,y_4)=G(L(y_2,v_1),L(y_4,v_1))
\left\{\begin{array}{l}
\overset{(\ref{10.8})}=G(L(y_3,x_1),L(y_5,x_1))=\tau G(y_3,y_5),  \\
\overset{(\ref{10.8})}=G(-L(x_2,v_2),-L(x_4,v_2))=0,
\end{array}\right.
\end{align*}
(As $L(x_1,y_4)=-L(y_5,v_1)$, $L(x_4,y_4)=L(v_1,v_2)$ in (\ref{10.17}) and Lemma \ref{lemma5.13}, we get $L(y_4,v_1)=L(y_5,x_1)$, $L(y_4,v_1)=-L(x_4,v_2)$.)
\begin{align*}
\tau G(y_3,y_4)=G(L(y_3,v_1),L(y_4,v_1))
\overset{(\ref{10.7})}=G(-L(x_3,v_2),-L(x_4,v_2))=0,
\end{align*}
(As $L(x_4,y_4)=L(v_1,v_2)$ in (\ref{10.17}) and Lemma \ref{lemma5.13}, we get $L(y_4,v_1)=-L(x_4,v_2)$.)
\begin{align*}
\tau G(y_4,y_5)=G(L(y_4,v_1),L(y_5,v_1))
\overset{(\ref{10.17})}=G(L(y_4,v_1),-L(y_4,x_1))=0,
\end{align*}
\begin{align*}
\tau G(y_4,y_6)=G(L(y_4,v_1),L(y_6,v_1))
\overset{(\ref{10.18})}=G(-L(x_4,v_2),L(x_4,y_2))=0,
\end{align*}
(As $L(x_4,y_4)=L(v_1,v_2)$ in (\ref{10.17}) and Lemma \ref{lemma5.13}, we get $L(y_4,v_1)=-L(x_4,v_2)$.)
\begin{align*}
\tau G(y_4,y_7)=G(L(y_4,v_1),L(y_7,v_1))
\overset{(\ref{10.19})}=G(-L(x_4,v_2),-L(x_7,v_2))=0.
\end{align*}

Hence, we get claim 2.

From (\ref{10.6}),(\ref{10.17})-(\ref{10.19}), it follows immediately that
\be\label{10.20}
L\left(x_{i}, y_{i}\right)=L\left(v_{1}, v_{2}\right), \quad i=1, \ldots, 7
\ee
and therefore, by Lemma \ref{lemma5.13}, we obtain
\be\label{10.21}
L\left(x_{i}, v_{2}\right)=-L\left(y_{i}, v_{1}\right).
\ee
For $i\neq j$, we have
$$G(L\left(x_{i}, y_{j}\right),-L\left(y_{i}, x_{j}\right))\overset{(\ref{5.8})}=G(L\left(x_{i}, y_{i}\right),L\left(y_{j}, x_{j}\right))=G(L\left(v_{1}, v_{2}\right),L\left(v_{1}, v_{2}\right))=\tau$$
and
$$G(L\left(x_{i}, y_{j}\right),L\left(x_{i}, y_{j}\right))=G(L\left(y_{i}, x_{j}\right),L\left(y_{i}, x_{j}\right))=\tau.$$
Thus
\be\label{10.22}
L\left(x_{i}, y_{j}\right)=-L\left(y_{i}, x_{j}\right),\; 1 \leq i \neq j \leq 7.  \ee

Finally, based on the relations (\ref{10.2})-(\ref{10.9}) and (\ref{10.17})-(\ref{10.22}), the following relations can be established:
\begin{align}
\label{10.23}
&\left\{\begin{array}{l}
L\left(x_{4}, y_{5}\right)=-L\left(v_{1}, y_{1}\right),(\Leftrightarrow
L\left(x_{4}, y_{1}\right)=L\left(v_{1}, y_{5}\right)\; in\; (\ref{10.17})) \\
L\left(x_{4}, y_{6}\right)=-L\left(v_{1}, y_{2}\right),(\Leftrightarrow
L\left(x_{4}, y_{2}\right)=L\left(v_{1}, y_{6}\right)\; in\; (\ref{10.18})) \\
L\left(x_{4}, y_{7}\right)=-L\left(v_{1}, y_{3}\right),(\Leftrightarrow
L\left(x_{4}, y_{3}\right)=L\left(v_{1}, y_{7}\right)\; in\; (\ref{10.19}))
\end{array}\right.
\end{align}
\begin{align}
\label{10.24}
&\left\{\begin{array}{l}
L\left(x_{5}, y_{1}\right)=-L\left(v_{1}, y_{4}\right), (by\; (\ref{10.17})\; and\; (\ref{10.21})) \\
L\left(x_{5}, y_{2}\right)=L\left(v_{1}, y_{7}\right) \Leftrightarrow
L\left(x_{5}, y_{7}\right)=-L\left(v_{1}, y_{2}\right),\\
L\left(x_{5}, y_{3}\right)=-L\left(v_{1}, y_{6}\right) \Leftrightarrow
L\left(x_{5}, y_{6}\right)=L\left(v_{1}, y_{3}\right).
\end{array}\right.
\end{align}
From
\begin{align*}
&G(L\left(v_{1}, y_{7}\right),L\left(x_{5}, y_{2}\right))\overset{(\ref{10.19})(\ref{10.22})}=
G(L\left(x_{4}, y_{3}\right),-L\left(y_{5}, x_{2}\right)) \\
\overset{(\ref{5.8})}=&G(L\left(x_{4}, y_{5}\right),L\left(y_{3}, x_{2}\right))\overset{(\ref{10.9})(\ref{10.23})}=
G(-L\left(v_{1}, y_{1}\right),-L\left(v_{1}, y_{1}\right))=\tau,  \\
&G(L\left(x_{5}, y_{3}\right),-L\left(v_{1}, y_{6}\right))\overset{(\ref{10.18})}=
G(L\left(x_{5}, y_{3}\right),-L\left(x_{4}, y_{2}\right)) \\
\overset{(\ref{5.8})}=&G(L\left(x_{5}, y_{2}\right),L\left(y_{3}, x_{4}\right))\overset{(\ref{10.22})}=
G(-L\left(y_{5}, x_{2}\right),L\left(y_{3}, x_{4}\right)))\overset{(\ref{10.23})(\ref{10.24})}=\tau,
\end{align*}
and Cauchy-Schwarz inequality, we have $L\left(v_{1}, y_{7}\right)=L\left(x_{5}, y_{2}\right)$ and \\ $L\left(x_{5}, y_{3}\right)=-L\left(v_{1}, y_{6}\right)$.
\begin{align}
\label{10.25}
&\left\{\begin{array}{l}
L\left(x_{6}, y_{1}\right)=-L\left(v_{1}, y_{7}\right) \Leftrightarrow L\left(x_{6}, y_{7}\right)=L\left(v_{1}, y_{1}\right),  \\
L\left(x_{6}, y_{3}\right)=L\left(v_{1}, y_{5}\right), \quad L\left(x_{6}, y_{2}\right)=-L\left(v_{1}, y_{4}\right).
\end{array}\right.
\end{align}
From
\begin{align*}
&G(-L\left(v_{1}, y_{7}\right),L\left(x_{6}, y_{1}\right))\overset{(\ref{10.19})(\ref{10.22})}=
G(-L\left(x_{4}, y_{3}\right),-L\left(y_{6}, x_{1}\right)) \\
\overset{(\ref{5.8})}=&G(-L\left(x_{4}, y_{6}\right),L\left(y_{3}, x_{1}\right))\overset{(\ref{10.8})(\ref{10.23})}=
G(L\left(v_{1}, y_{2}\right),L\left(v_{1}, y_{2}\right))=\tau,  \\
&G(L\left(x_{6}, y_{3}\right),L\left(v_{1}, y_{5}\right))\overset{(\ref{10.17})}=
G(L\left(x_{6}, y_{3}\right),-L\left(x_{1}, y_{4}\right)) \\
\overset{(\ref{5.8})}=&G(L\left(x_{6}, y_{4}\right),L\left(y_{3}, x_{1}\right))\overset{(\ref{10.22})(\ref{10.8})}=G(-L\left(x_{4}, y_{6}\right),L\left(y_{2}, v_{1}\right))\\ &\overset{(\ref{10.23})}=
G(L\left(y_{2}, v_{1}\right),L\left(y_{2}, v_{1}\right))=\tau, \\
&G(L\left(x_{6}, y_{2}\right),-L\left(v_{1}, y_{4}\right))\overset{(\ref{10.24})}=
G(L\left(x_{6}, y_{2}\right),L\left(x_{5}, y_{1}\right)) \\
\overset{(\ref{5.8})}=&G(-L\left(x_{6}, y_{1}\right),L\left(y_{2}, x_{5}\right))\overset{(\ref{10.24})(\ref{10.25})}=
G(L\left(y_{7}, v_{1}\right),L\left(y_{7}, v_{1}\right))=\tau,
\end{align*}
and Cauchy-Schwarz inequality, we obtain (\ref{10.25}).
\begin{align}
\label{10.26}
&\left\{\begin{array}{l}
L\left(x_{7}, y_{1}\right)=L\left(v_{1}, y_{6}\right), \quad L\left(x_{7}, y_{2}\right)=-L\left(v_{1}, y_{5}\right), \\
L\left(x_{7}, y_{3}\right)=-L\left(v_{1}, y_{4}\right) (by\; (\ref{10.19})\; and\; (\ref{10.21})).
\end{array}\right.
\end{align}
From
\begin{align*}
&G(L\left(x_{7}, y_{1}\right),L\left(v_{1}, y_{6}\right))\overset{(\ref{10.24})}=
G(L\left(x_{7}, y_{1}\right),-L\left(x_{5}, y_{3}\right)) \\
\overset{(\ref{5.8})}=&G(L\left(x_{7}, y_{3}\right),L\left(y_{1}, x_{5}\right))\overset{(\ref{10.19})(\ref{10.24})}=
G(L\left(x_{4}, v_{2}\right),-L\left(v_{1}, y_{4}\right))\\ & \overset{(\ref{10.21})}=
G(L\left(x_{4}, v_{2}\right),G(L\left(x_{4}, v_{2}\right))=\tau,  \\
&G(L\left(x_{7}, y_{2}\right),-L\left(v_{1}, y_{5}\right))\overset{(\ref{10.25})}=
G(L\left(x_{7}, y_{2}\right),-L\left(x_{6}, y_{3}\right)) \\
\overset{(\ref{5.8})}=&G(L\left(x_{7}, y_{3}\right),L\left(y_{2}, x_{6}\right))\overset{(\ref{10.19})(\ref{10.25})}=G(L\left(x_{4}, v_{2}\right),-L\left(y_{4}, v_{1}\right))\\ &\overset{(\ref{10.21})}=
G(L\left(x_{4}, v_{2}\right),L\left(x_{4}, v_{2}\right))=\tau,
\end{align*}
and Cauchy-Schwarz inequality, we obtain (\ref{10.26}).

All relations in (\ref{10.1}) verified, and thus we complete the proof of Lemma \ref{lemma10.2}.   \hfill $\Box$

Now, we can present the following crucial and remarkable lemma as in \cite{CHM}.

\begin{lem}\label{lemma10.3}
Suppose that in the decomposition (\ref{5.18}) we have $k_{0} \geq 2$ and $\mathfrak{p}=7$. Then it must be the case that $k_{0}=2$.
\end{lem}

\textbf{Proof.}
Suppose on the contrary that $k_{0} \geq 3$. Following the same argument as in the proof of Lemma \ref{lemma10.2} for $V_{v_{1}}(0)$ and $V_{v_{2}}(0)$, we choose a basis $\left\{x_{1}, x_{2}, \td{x}_{3}, x_{4}, \td{x}_{5}, \td{x}_{6}, \td{x}_{7}\right\}$ of $V_{v_{1}}(0)$ and a basis $\left\{z_{1}, z_{2}, z_{3}, z_{4}, z_{5}, z_{6}, z_{7}\right\}$ of $V_{v_{3}}(0)$ such that all the following relations hold:
\be\label{10.27}
L\left(e_{j}\left(v_{1}\right), e_{l}\left(v_{3}\right)\right)=-L\left(v_{1}, e_{j} e_{l}\left(v_{3}\right)\right)=-L\left(e_{l} e_{j}\left(v_{1}\right), v_{3}\right), \;1 \leq j, l \leq 7.
\ee
\textbf{Step 1.}
Now, we have two orthonormal bases of $V_{v_{1}}(0)$, i.e. $\left\{x_{1}, x_{2}, \td{x}_{3}, x_{4}, \td{x}_{5}, \td{x}_{6}, \td{x}_{7}\right\}$ and $\left\{x_{1}, x_{2}, x_{3}, x_{4}, x_{5}, x_{6}, x_{7}\right\}$. We first show that $\td{x}_{i}=x_{i}$ for $i=3,5,6,7$.

By (\ref{5.23}), (\ref{10.1}) and (\ref{10.27}), we get
\begin{align*}
&\tau L\left(y_{1}, z_{1}\right)=A\left(L\left(y_{1}, x_{1}\right), L\left(x_{1}, z_{1}\right)\right)=A\left(L\left(v_{1}, v_{2}\right),L\left(v_{1}, v_{3}\right)\right)=\tau L\left(v_{2}, v_{3}\right),  \\
&\tau L\left(y_{2}, z_{2}\right)=A\left(L\left(y_{2}, x_{2}\right), L\left(x_{2}, z_{2}\right)\right)=A\left(L\left(v_{1}, v_{2}\right),L\left(v_{1}, v_{3}\right)\right)=\tau L\left(v_{2}, v_{3}\right),   \\
&\tau L\left(y_{4}, z_{4}\right)=A\left(L\left(y_{4}, x_{4}\right), L\left(x_{4}, z_{4}\right)\right)=A\left(L\left(v_{1}, v_{2}\right),L\left(v_{1}, v_{3}\right)\right)=\tau L\left(v_{2}, v_{3}\right).
\end{align*}

Since $\left\{x_{1}, x_{2}, \td{x}_{3}, x_{4}, \td{x}_{5}, \td{x}_{6}, \td{x}_{7}\right\}$ and $\left\{x_{1}, x_{2}, x_{3}, x_{4}, x_{5}, x_{6}, x_{7}\right\}$ are two orthonormal bases for $V_{v_{1}}(0)$, we may assume that $x_{3}=b_{3} \td{x}_{3}+b_{5} \td{x}_{5}+b_{6} \td{x}_{6}+b_{7} \td{x}_{7}$. Then we have the following calculation
\begin{align}\label{10.28}
\tau L\left(y_{2}, z_{2}\right)
\overset{(\ref{5.23})}=& A\left(L\left(v_{1}, y_{2}\right), L\left(v_{1}, z_{2}\right)\right)\overset{(\ref{10.1})}=-A\left(L\left(x_{3}, y_{1}\right), L\left(v_{1}, z_{2}\right)\right)  \notag    \\
\overset{(\ref{10.27})}=& b_{3} A\left(L\left(\td{x}_{3}, y_{1}\right), L\left(\td{x}_{3}, z_{1}\right)\right)+b_{5} A\left(L\left(\td{x}_{5}, y_{1}\right), L\left(\td{x}_{5}, z_{7}\right)\right)  \notag  \\
&-b_{6} A\left(L\left(\td{x}_{6}, y_{1}\right), L\left(\td{x}_{6}, z_{4}\right)\right)-b_{7} A\left(L\left(\td{x}_{7}, y_{1}\right), L\left(\td{x}_{7}, z_{5}\right)\right)   \notag \\
=\;& b_{3} \tau L\left(y_{1}, z_{1}\right)+b_{5} \tau L\left(y_{1}, z_{7}\right)-b_{6} \tau L\left(y_{1}, z_{4}\right)-b_{7} \tau L\left(y_{1}, z_{5}\right).
\end{align}

On the other hand, since $L\left(y_{1}, z_{1}\right)(=L\left(y_{2}, z_{2}\right)), L\left(y_{1}, z_{4}\right), L\left(y_{1}, z_{5}\right)$ and $L\left(y_{1}, z_{7}\right)$ are mutually orthogonal, (\ref{10.28}) implies that $b_{3}=1, b_{5}=b_{6}=b_{7}=0$ and hence $x_{3}=\td{x}_{3}$.

Let $x_{5}=b_{5} \td{x}_{5}+b_{6} \td{x}_{6}+b_{7} \td{x}_{7}$. Then we have the following calculation
\begin{align}\label{10.29}
\tau L\left(y_{4}, z_{4}\right)
\overset{(\ref{5.23})}=& A\left(L\left(v_{1}, y_{4}\right), L\left(v_{1}, z_{4}\right)\right)\overset{(\ref{10.1})}=-A\left(L\left(x_{5}, y_{1}\right), L\left(v_{1}, z_{4}\right)\right)  \notag    \\
\overset{(\ref{10.27})}=& b_{5} A\left(L\left(\td{x}_{5}, y_{1}\right), L\left(\td{x}_{5}, z_{1}\right)\right)  \notag  \\
&+b_{6} A\left(L\left(\td{x}_{6}, y_{1}\right), L\left(\td{x}_{6}, z_{2}\right)\right)+b_{7} A\left(L\left(\td{x}_{7}, y_{1}\right), L\left(\td{x}_{7}, z_{3}\right)\right)   \notag \\
=\; & b_{5} \tau L\left(y_{1}, z_{1}\right)+b_{6} \tau L\left(y_{1}, z_{2}\right)+b_{7} \tau L\left(y_{1}, z_{3}\right).
\end{align}

On the other hand, since $L\left(y_{1}, z_{1}\right)(=L\left(y_{4}, z_{4}\right)), L\left(y_{1}, z_{2}\right)$ and $ L\left(y_{1},z_{3}\right)$ are mutually orthogonal, (\ref{10.29}) implies that $b_{5}=1, b_{6}=b_{7}=0$ and hence $x_{5}=\td{x}_{5}$.

By $x_{3}=\td{x}_{3}$ and $x_{5}=\td{x}_{5}$, we get
\begin{align*}
&\tau L\left(y_{3}, z_{3}\right)\overset{(\ref{5.23})}=A\left(L\left(y_{3}, x_{3}\right), L\left(x_{3}, z_{3}\right)\right)\overset{(\ref{10.1})(\ref{10.27})}=A\left(L\left(v_{1}, v_{2}\right),L\left(v_{1}, v_{3}\right)\right)=\tau L\left(v_{2}, v_{3}\right),  \\
&\tau L\left(y_{5}, z_{5}\right)\overset{(\ref{5.23})}=A\left(L\left(y_{5}, x_{5}\right), L\left(x_{5}, z_{5}\right)\right)\overset{(\ref{10.1})(\ref{10.27})}=A\left(L\left(v_{1}, v_{2}\right),L\left(v_{1}, v_{3}\right)\right)=\tau L\left(v_{2}, v_{3}\right).
\end{align*}

Let $x_{6}=b_{6} \td{x}_{6}+b_{7} \td{x}_{7}$. Then we have the following calculation
\begin{align}\label{10.30}
\tau L\left(y_{3}, z_{3}\right)
\overset{(\ref{5.23})}=& A\left(L\left(v_{1}, y_{3}\right), L\left(v_{1}, z_{3}\right)\right)\overset{(\ref{10.1})}=-A\left(L\left(x_{6}, y_{5}\right), L\left(v_{1}, z_{3}\right)\right)  \notag    \\
\overset{(\ref{10.27})}=& b_{6} A\left(L\left(\td{x}_{6}, y_{5}\right), L\left(\td{x}_{6}, z_{5}\right)\right)-b_{7} A\left(L\left(\td{x}_{7}, y_{5}\right), L\left(\td{x}_{7}, z_{4}\right)\right)   \notag \\
=& b_{6} \tau L\left(y_{5}, z_{5}\right)-b_{7} \tau L\left(y_{5}, z_{4}\right).
\end{align}

On the other hand, since $L\left(y_{5}, z_{5}\right)(=L\left(y_{3}, z_{3}\right))$ and $ L\left(y_{5},z_{4}\right)$ are mutually orthogonal, (\ref{10.30}) implies that $b_{6}=1, b_{7}=0$, and hence $x_{6}=\td{x}_{6}$.

Finally, let $x_7=\epsilon\td x_7$, by $L\left(y_{1}, z_{1}\right)=L\left(v_{2}, v_{3}\right)$,
\begin{align*}
\tau L\left(y_{6}, z_{6}\right)
\overset{(\ref{5.23})}=A\left(L\left(x_{6}, y_{6}\right), L\left(x_{6}, z_{6}\right)\right)\overset{(\ref{10.1})(\ref{10.27})}=A\left(L\left(v_{1}, v_{2}\right), L\left(v_{1}, v_{3}\right)\right)=\tau L\left(v_{2}, v_{3}\right)
\end{align*}
and
\begin{align*}
\tau L\left(y_{6}, z_{6}\right)\overset{(\ref{5.23})}=
&A\left(L\left(v_{1}, y_{6}\right), L\left(v_{1}, z_{6}\right)\right)\overset{(\ref{10.1})}
=A\left(L\left(x_{7}, y_{1}\right), L\left(v_{1}, z_{6}\right)\right) \\
\overset{(\ref{10.27})}=&\epsilon A\left(L\left(\td x_{7}, y_{1}\right), L\left(\td x_{7}, z_{1}\right)\right)=\epsilon \tau L\left(y_{1}, z_{1}\right),
\end{align*}
we get $\epsilon=1$ and hence $x_7=\td x_7$.

\textbf{Step 2.}
We will use (\ref{10.1}) and (\ref{10.27}) to show that we have also similar relations between $V_{v_{2}}(0)$ and $V_{v_{3}}(0)$, i.e.,
\be\label{10.31}
L\left(e_{j}\left(v_{2}\right), e_{l}\left(v_{3}\right)\right)=-L\left(v_{2}, e_{j} e_{l}\left(v_{3}\right)\right)=-L\left(e_{l} e_{j}\left(v_{2}\right), v_{3}\right), \quad 1 \leq j, l \leq 7.
\ee
In fact, for $j=l$, by (\ref{5.23}), (\ref{10.1}) and (\ref{10.27}), we have
$$
\begin{aligned}
\tau L\left(e_{j}\left(v_{2}\right), e_{j}\left(v_{3}\right)\right) &=A\left(L\left(e_{k}\left(v_{1}\right), e_{j}\left(v_{2}\right)\right), L\left(e_{k}\left(v_{1}\right), e_{j}\left(v_{3}\right)\right)\right) \\
&=A\left(L\left(v_{2}, e_{j} e_{k}\left(v_{1}\right)\right), L\left(e_{j} e_{k}\left(v_{1}\right), v_{3}\right)\right)=\tau L\left(v_{2}, v_{3}\right).
\end{aligned}
$$
For $j \neq l$, according to the multiplication table in Lemma \ref{lemma10.2}, there exists a unique $k$ and $\epsilon=\pm 1$ such that $e_{l} e_{j}=\epsilon e_{k}, e_{j} e_{k}=\epsilon e_{l}, e_{k} e_{l}=\epsilon e_{j} .$ It follows, by applying (\ref{5.23}), (\ref{10.1}) and (\ref{10.27}), that
$$
\begin{aligned}
\tau L\left(e_{j}\left(v_{2}\right), e_{l}\left(v_{3}\right)\right) &=A\left(L\left(v_{1}, e_{j}\left(v_{2}\right)\right), L\left(v_{1}, e_{l}\left(v_{3}\right)\right)\right) \\
&=A\left(-L\left(v_{1}, \epsilon e_{l} e_{k}\left(v_{2}\right)\right), L\left(v_{1}, e_{l}\left(v_{3}\right)\right)\right) \\
&=\epsilon A\left(L\left(e_{k}\left(v_{2}\right), e_{l}\left(v_{1}\right)\right),-L\left(v_{3}, e_{l}\left(v_{1}\right)\right)\right) \\
&=-\epsilon \tau L\left(e_{k}\left(v_{2}\right), v_{3}\right)=-\tau L\left(e_{l} e_{j}\left(v_{2}\right), v_{3}\right),
\end{aligned}
$$
and that
$$
\begin{aligned}
\tau L\left(v_{2}, e_{j} e_{l}\left(v_{3}\right)\right) &=A\left(L\left(e_{k}\left(v_{1}\right), v_{2}\right), L\left(e_{j} e_{l}\left(v_{3}\right), e_{k}\left(v_{1}\right)\right)\right) \\
&=A\left(L\left(v_{2}, \epsilon e_{l} e_{j}\left(v_{1}\right)\right), L\left(-\epsilon e_{k}\left(v_{3}\right), e_{k}\left(v_{1}\right)\right)\right) \\
&=A\left(L\left(v_{1},-\epsilon e_{l} e_{j}\left(v_{2}\right)\right), L\left(-\epsilon v_{3}, v_{1}\right)\right)=\tau L\left(e_{l} e_{j}\left(v_{2}\right), v_{3}\right).
\end{aligned}
$$
Thus, (\ref{10.31}) holds indeed.

\textbf{Step 3.}
From (\ref{10.1}), we have
$$L\left(v_{1}, y_{6}\right)=-L\left(x_{1}, y_{7}\right).$$
Hence
\be\label{10.32}
A\left(L\left(v_{1}, y_{6}\right)+L\left(x_{1}, y_{7}\right), L\left(x_{2}, v_{3}\right)\right)=0.
\ee
On the other hand, we have
$$
\begin{aligned}
&A\left(L\left(v_{1}, y_{6}\right), L\left(x_{2}, v_{3}\right)\right)\overset{(\ref{10.27})}=A\left(L\left(v_{1}, y_{6}\right),-L\left(v_{1}, z_{2}\right)\right)=-\tau L\left(z_{2}, y_{6}\right), \\
&A\left(L\left(x_{1}, y_{7}\right), L\left(x_{2}, v_{3}\right)\right)\overset{(\ref{10.27})}=A\left(L\left(x_{1}, y_{7}\right),-L\left(x_{1}, z_{3}\right)\right)=-\tau L\left(z_{3}, y_{7}\right).
\end{aligned}
$$
These, together with (\ref{10.32}), give that
\be\label{10.33}
L\left(z_{2}, y_{6}\right)+L\left(z_{3}, y_{7}\right)=0.
\ee
(\ref{10.31}) implies that $L\left(z_{2}, y_{6}\right)=-L\left(v_{2}, z_{4}\right)=L\left(z_{3}, y_{7}\right)$, and by $(\ref{10.33})$ we get $L\left(z_{2}, y_{6}\right)=0$.
However, we also have the relation $G\left(L\left(z_{2}, y_{6}\right), L\left(z_{2}, y_{6}\right)\right)=\tau$, which gives the contradiction.
This completes the proof of Lemma \ref{lemma10.3}.   \hfill $\Box$

Now, we are ready to complete the proof of Theorem \ref{thm10.1}.

\textbf{Proof of Theorem \ref{thm10.1}.}
First, Lemma \ref{lemma10.3} implies that $k_{0}=2$ and $\operatorname{dim}\left(\mathcal{D}_{2}\right)=16$.

Let $\left\{v_{1}, v_{2}, x_{j}, y_{j}, 1 \leq j \leq 7\right\}$ be the orthonormal basis of $\mathcal{D}_{2}$ as constructed in Lemma \ref{lemma10.2} such that all relations in (\ref{10.1}) hold. Then we easily see that the image of $L$ is spanned by
$$
\left\{L\left(v_{1}, v_{1}\right), L\left(v_{1}, v_{2}\right), L\left(v_{2}, v_{2}\right) ,\left.L\left(v_{1}, y_{j}\right)\right|_{1 \leq j \leq 7}\right\}
$$
and $\operatorname{dim}\left(\operatorname{Im} L\right)\leq 10$.

Define $L_{1}=L\left(v_{1}, v_{1}\right)-L\left(v_{2}, v_{2}\right)$ and $\operatorname{Tr} L=8\left(L\left(v_{1}, v_{1}\right)+L\left(v_{2}, v_{2}\right)\right)$, then we have
\begin{align}
\label{10.34}
G\left(L_{1}, L_{1}\right)\overset{(\ref{5.2})(\ref{5.6})}=&4 \tau \neq 0,  \\
\label{10.35}
\fr{1}{64}G\left(\operatorname{Tr} L, \operatorname{Tr} L\right)\overset{(\ref{5.2})(\ref{5.6})}=&4\tau =:\rho^2 \neq 0.
\end{align}
We now easily see that there exist ten orthonormal vectors in $\operatorname{Im} L \subset \mathcal{D}_{3}$ :
$$
t=\frac{1}{8\rho}\operatorname{Tr} L, w_{0}=\frac{1}{\sqrt{4 \tau}} L_{1}, w_{1}=\frac{1}{\sqrt{\tau}} L\left(v_{1}, v_{2}\right), w_{j+1}=\frac{1}{\sqrt{\tau}} L\left(v_{1}, y_{j}\right), 1 \leq j \leq 7.
$$
Then we have the estimate of the dimension
$$
n=1+\operatorname{dim}\left(\mathcal{D}_{2}\right)+\operatorname{dim}\left(\mathcal{D}_{3}\right) \geq 27.
$$

Now, we are ready to complete the proof of Theorem \ref{thm10.1}.

\textbf{Proof of Theorem \ref{thm10.1}.} We need to consider two cases:

Case (1): \ $n=27$\ ($\operatorname{Im} L= \mathcal{D}_{3}$), and  Case (2):\ $n>27$ \ ($\operatorname{Im} L \varsubsetneqq \mathcal{D}_{3}$).

In Case (1), from the previous discussions, we see that
$$\{t, w_0, w_1, w_{j+1}|_{1\leq j\leq 7}\}$$
forms an orthonormal basis of $\operatorname{Im} L=\mathcal{D}_{3}$.

In Case (2), we choose  $\{\td w_i|_{1\leq i\leq \td n}\}$
as an orthonormal basis of $\mathcal{D}_{3}\backslash \operatorname{Im} L$, and
such that $$\{t, w_0, w_1, w_{j+1}|_{1\leq j\leq 7}, \td w_i|_{1\leq i\leq \td n}\}$$
forms an orthonormal basis of $\mathcal{D}_{3}$,  where $\td n=n-27$.

\begin{lem}\label{lemma10.4}
The Fubini-Pick tensor satisfies
$$
\left\{\begin{array}{lll}
A(e_1,t)=0, \quad A(t,v_i)=\fr{\rho}{2}v_i, \ 1\leq i\leq 2, \quad A(t,x_j)=\fr{\rho}{2}x_j,\ 1\leq j\leq 7, \\
A(t,y_j)=\fr{\rho}{2}y_j,\ 1\leq j\leq 7,\quad A(t,t)=\rho t, \quad A(t,w_k)=\rho w_k, \quad 0\leq k\leq 8,\\
\text { if } \operatorname{Im} L \varsubsetneqq \mathcal{D}_{3}, \quad A(t,\td w_i)=0, \quad 1\leq i\leq \td n.
\end{array}
\right.
$$
\end{lem}

\textbf{Proof.}
By (iii) of Lemma \ref{lemma5.1}, we know $A(t,v_i)\in \mathcal{D}_2$. By
\begin{align*}
\lang A(t,v_i),v_j\rang&=\fr{1}{\rho}\sum_{l=1}^{2}\lang A(v_i,v_j),L(v_l,v_l)\rang \overset{(\ref{5.5})(\ref{5.6})}=\fr{\rho}{2}\delta_{ij},\quad 1\leq i,j\leq 2,\\
\lang A(t,v_i),x_j\rang&=\fr{1}{\rho}\sum_{l=1}^{2}\lang A(v_i,x_j),L(v_l,v_l)\rang \overset{Lemma\ref{lemma5.3}}=0,\quad 1\leq i\leq 2,\quad 1\leq j\leq 7, \\
\lang A(t,v_i),y_j\rang&=\fr{1}{\rho}\sum_{l=1}^{2}\lang A(v_i,y_j),L(v_l,v_l)\rang \overset{Lemma\ref{lemma5.3}}=0,\quad 1\leq i\leq 2,\quad 1\leq j\leq 7,
\end{align*}
we have $A(t,v_i)=\fr{\rho}{2}v_i, \ 1\leq i\leq 2$. Similarly, $A(t,x_j)=\fr{\rho}{2}x_j,\ A(t,y_j)=\fr{\rho}{2}y_j,\ 1\leq j\leq 7$ can be proved.

From Lemmas \ref{lemma5.14} and \ref{lemma10.2} we get
\begin{align*}
A(t,t)&=\fr{1}{\rho^2} \sum_{i=1}^{2}A(L(v_i,v_i),L(v_i,v_i))=\rho t, \\
A(t,w_0)&=\fr{\mu_1^2}{2\rho\sqrt{4\tau}}L_1=\rho w_0, \\
A(t,w_{1})&=\fr{1}{\rho\sqrt \tau}(A(L(v_1,v_1),L(v_1,v_2))+A(L(v_2,v_2),L(v_1,v_2)))=\rho w_{1}, \\
A(t,w_{j+1})&=\fr{1}{\rho\sqrt \tau}(A(L(v_1,v_1),L(v_1,y_j))+A(L(v_2,v_2),L(v_1,y_j)))   \\
&=\fr{1}{\rho\sqrt \tau}(2\tau L(v_1,y_j)+A(L(y_j,y_j),L(v_1,y_j)))=\rho w_{j+1}, \quad 1\leq j\leq 7.
\end{align*}
If $\operatorname{Im} L \varsubsetneqq \mathcal{D}_{3}$, by Lemma \ref{lemma5.4}, one get
$$A(t,\td w_i)=\fr{1}{\rho}\sum_{j=1}^{2}A(L(v_j,v_j),\td w_i)=0, \quad 1\leq i\leq \td n.$$
This completes the proof of Lemma \ref{lemma10.4}.    \hfill $\Box$

Put
\be\label{10.36}
T=\fr{1}{\sqrt{3}}e_1+\sqrt{\fr{2}{3}}t,\quad T^*=-\sqrt{\fr{2}{3}}e_1+\fr{1}{\sqrt{3}}t.
\ee
It is easily to see that if $\operatorname{Im} L=\mathcal{D}_{3}$, then $\{T, T^*, v_i|_{1\leq i\leq 2} , x_j|_{1\leq j\leq 7}, y_j|_{1\leq j\leq 7}, w_k|_{0\leq k\leq 8}\}$ (or, resp. if $\operatorname{Im} L \varsubsetneqq \mathcal{D}_{3}$, then $\{T, T^*, v_i|_{1\leq i\leq 2} , x_j|_{1\leq j\leq 7}, y_j|_{1\leq j\leq 7}, w_k|_{0\leq k\leq 8}, \td w_i|_{1\leq i\leq \td n}\}$) is an orthonormal basis of $T_pM$. By Lemmas \ref{lemma5.1} and \ref{lemma10.4} we have the following
\begin{lem}\label{lemma10.5}
With respect to the above notations, the Fubini-Pick tensor takes the following form
$$\left\{
\begin{array}{ll}
A(T,T)=\sigma T, \quad A(T,T^*)=\sigma T^*,\quad A(T,v_i)=\sigma v_i, \ (1\leq i\leq 2)  \\
A(T,x_j)=\sigma x_j, \quad A(T,y_j)=\sigma y_j, \ (1\leq j\leq 7) \quad A(T,w_k)=\sigma w_k, \ (0\leq k\leq 8)      \\
\text { if } \operatorname{Im} L \varsubsetneqq \mathcal{D}_{3}, \quad A(T,\td w_i)=0, \quad (1\leq i\leq \td n),
\end{array}
\right.
$$
where $\sigma=\fr{1}{\sqrt{3}}\mu_1$.
\end{lem}

Using the parallelism of Fubini-Pick form, Lemma 10.5, Theorems 3.3 and 3.6, we complete the proof of Theorem  \ref{thm10.1}. \hfill $\Box$

\section{Completion of the proof of theorem 1.1}

If the Fubini-Pick form $A=0$, according to Lemma 4.3 of \cite{XL1}, we have (i).

For hypersurfaces with $A\neq 0$, according to Lemma 4.2, it is necessary and sufficient to consider the cases $\left\{\mathfrak{C}_{m}\right\}_{1 \leq m \leq n}$. Firstly, by Lemma $4.4$ in \cite{XL1},  Corollaries 4.5 and  4.6, we have settled the two cases, $\mathfrak{C}_{n}$, $\mathfrak{C}_{1}$ and $\mathfrak{C}_{n-1}$, from which we have (ii) and (iii).

The remaining cases, i.e. $\mathfrak{C}_{m}$ with $2 \leq m \leq n-2$, are completely settled by Proposition $5.15$ and subsequent five theorems, i.e. Theorems $6.1,7.1,8.1,9.1$ and $10.1$. In these cases, we have (ii) and (vii).
From all of above discussions, we have completed the proof of Theorem 1.1.
\hfill $\Box$

\vskip 0.1in
\textbf{Acknowledgements}

To work for this project on Calabi hypersurfaces we have greatly benefitted from the work of \cite{HLV} and \cite{CHM}. Therefore, the authors would like to express their heartfelt thanks to Profs. Xiuxiu Cheng, Zejun Hu, Haizhong Li, Marilena Moruz and Luc Vrancken.


\end{document}